\documentclass[a4paper, reqno,12pt]{amsart}
\usepackage{amsmath, amsthm, amssymb, amstext}

\usepackage[left=0.75in,right=0.75in,top=0.75in,bottom=0.75in]{geometry}

\usepackage{xurl,hyperref,xcolor}

\usepackage[hyperpageref]{backref}
\usepackage{esint}
\hypersetup{pdfborder={0 0 0},colorlinks}
\allowdisplaybreaks
\usepackage{enumitem}
\usepackage{extarrows}
\allowdisplaybreaks
\newtheorem{knownthm}{Theorem}[section]
\newtheorem{knownlem}[knownthm]{Lemma}

\newtheorem{theorem}{Theorem}

\newtheorem{lemma}[theorem]{Lemma}
\newtheorem{proposition}[theorem]{Proposition}
\newtheorem{corollary}[theorem]{Corollary}

\newtheorem{example}[theorem]{Example}

\numberwithin{theorem}{section}
\numberwithin{equation}{section}

\usepackage{etoolbox}

\AtBeginEnvironment{thebibliography}{%
  \footnotesize
  \setlength{\itemsep}{0pt}%
  \setlength{\parsep}{0pt}%
  \setlength{\parskip}{0pt}%
}
\title[Sharp Constants and Maximizers for Moser--Trudinger Inequalities]{Sharp Constants and Maximizers for Weighted Moser--Trudinger Inequalities on the Heisenberg Group}
\author[S. Kumar]{Somnath Kumar}
\address[S. Kumar]{Department of Mathematics, Indian Institute of Technology Jodhpur, Rajasthan 342030, India}
\email{somnathkumar325@gmail.com; p25ma0012@iitj.ac.in}
\author[N. Lam]{Nguyen Lam}
\address[N. Lam]{School of Science and the Environment, Grenfell Campus, Memorial University of Newfoundland,
Corner Brook, NL A2H5G4, Canada}
\email{nlam@mun.ca}
\author[D.K. Mahanta]{Deepak Kumar Mahanta}
\address[D.K. Mahanta]{School of Mathematics, Indian Institute of Science Education and Research Thiruvananthapuram, Maruthamala P.O., Vithura, Thiruvananthapuram-695551, Kerala, India}
\email{mahanta.1@iitj.ac.in; deepakmahanta@iisertvm.ac.in}
\author[A. Sarkar]{Abhishek Sarkar}
\address[A. Sarkar]{Department of Mathematics, Indian Institute of Technology Jodhpur, Rajasthan 342030, India}
\email{abhisheks@iitj.ac.in}
\begin{document}
\begin{abstract}
This article concerns various Moser--Trudinger inequalities on the entire Heisenberg group, with particular emphasis on their sharp constants, maximizers, and applications to subelliptic partial differential equations. It resolves several questions concerning the existence of maximizers for Moser--Trudinger inequalities on the Heisenberg group in the singular case, that is, $0<\beta<Q$, where the singular weight gives compactness. More precisely, we first establish an Adimurthi--Druet-type weighted Moser--Trudinger inequality without requiring technically involved blow-up analysis, and determine its sharp constant. The resulting inequality is also shown to be equivalent to the critical weighted Moser--Trudinger inequality. In addition, for the subcritical weighted Moser--Trudinger inequality established by Lam, Lu and Tang \cite{Lam-Lu-Tang-2014}, we obtain two-sided asymptotic bounds near the sharp threshold and show the existence of maximizers in a suitable regime. The subcritical weighted Moser--Trudinger inequality also provides further insight into borderline Sobolev theory and allows us to derive a Caffarelli--Kohn--Nirenberg interpolation-type inequality in the asymptotic sense. The weighted Moser--Trudinger inequality under inhomogeneous norm constraints is also obtained, and its equivalence with the subcritical weighted Moser--Trudinger inequality is subsequently demonstrated. This equivalence provides a characterization of the associated supremum and its maximizers in the subcritical regime. Moreover, an improved subcritical weighted Moser--Trudinger inequality is derived, and the associated maximization problem is shown to admit a maximizer. The developed inequalities also allow us to treat a class of normalized subelliptic equations driven by the $Q$-sub-Laplacian and nonlinearities with singular critical exponential growth at infinity, for which we establish the existence of nontrivial ground states for every prescribed mass $a>0$, provided that the nonlinearity is large enough in a suitable sense depending on $a$.
\end{abstract}

\maketitle
\section{Introduction and Main Results}\label{sec1}
The Moser--Trudinger inequality describes the limiting behavior of the Sobolev embedding at the critical exponent. In contrast to the subcritical regime, where power-type integrability governs the embedding, the critical case exhibits exponential integrability. Now, we briefly recall its classical origins and subsequent sharp form below.
\subsection{The Moser–Trudinger Inequality: Historical Context and Classical Theory} The Sobolev embedding theorem is a cornerstone of modern analysis. In its classic form, for a bounded domain $\Omega \subset \mathbb{R}^n$ and $1 \leq p < n$, it states that $W_0^{1,p}(\Omega) \hookrightarrow L^{p^\ast}(\Omega)$, where $ p^\ast = \frac{np}{n-p}$ is the critical Sobolev exponent. This embedding is optimal, meaning that it is not possible to replace $p^\ast$ with any larger value while maintaining a universal embedding. This provides the conventional formulation of the Sobolev embedding theorem.

The transition to the borderline case $p=n$ leads to a fundamental change, like the Sobolev embedding. In this context, the critical exponent $p^\ast$ becomes infinity, but this does not yield $L^\infty$-control: functions in $ W_0^{1,n}(\Omega)$ may remain unbounded and exhibit logarithmic singularities. This prompts the natural question:
  \[ \textit{What replaces the $L^\infty$-embedding in the borderline Sobolev regime?}\]
In the 1960s, Yudovich \cite{Yudovich-1961}, Poho\v{z}aev \cite{Pokhozhaev-1965}, and Trudinger \cite{Trudinger-1967} independently addressed this question by establishing that $W_0^{1,n}(\Omega)$ embeds into an exponential-type Orlicz space. More explicitly, there are constants $ \alpha>0$ and $C(n)> 0$ such that
$$
\sup_{u\in W^{1,n}_0(\Omega),\,\|\nabla u\|_{L^n(\Omega)}\,\leq 1}~\int_{\Omega} \operatorname{\exp}\big(\alpha|u|^{\frac{n}{n-1}}\big)\,\mathrm{d}x\leq C(n)|\Omega|.
$$
The precise value of the sharp constant is given by $\alpha_n = n \omega_{n-1}^{\frac{1}{n-1}}$, where $\omega_{n-1}$ is the surface measure of the unit sphere $\mathbb{S}^{n-1}$ in $\mathbb{R}^n$. This was first highlighted in Moser's esteemed 1971 paper \cite{Moser-1970/71}. He showed that $\alpha_n$ is indeed optimal: the inequality holds for $\alpha\leq\alpha_n$ but fails for every $\alpha>\alpha_n$. This sharp inequality is commonly referred to as the Moser--Trudinger inequality.

The significance of the Moser--Trudinger inequality extends beyond the study of borderline Sobolev embeddings. It provides a natural functional framework for a variety of geometric and physical problems whose Euler-Lagrange equations involve exponential nonlinearities; see, for instance, del Pino--Musso--Ruf \cite{delPino-Musso-Ruf-2010}. 

This link between a borderline embedding theorem and the existence theory for geometrically motivated PDEs has been, since the 1980s, a major driving force for further refinements: improved versions (see, \cite{Adi-Druet2004, Chang-Yang2003}), the search for extremal functions (see, \cite{CarlesonChang1986, Flucher1992, Li2001, Struwe1988}), higher-order generalizations (see, \cite{Adams-1988, Lam-Lu-2012@@}), and extensions of these inequalities to unbounded domains (see, \cite{Cao-1992, do O-1997, Lam-Lu-2012@@, Lam-Lu-2013, Ruf-2005, Ruf-Sani-2013}).

Next, a natural question, and a deep one, is whether the supremum defining the Moser--Trudinger functional is actually attained. For different contexts, we refer to \cite{CarlesonChang1986, Flucher1992, Li2001, Struwe1988}, and to del Pino--Musso--Ruf \cite{delPino-Musso-Ruf-2012} for critical points of the Moser--Trudinger functional at energy levels above the supremum. The underlying strategy, combining concentration–compactness and blow-up analysis, shows that a maximizing sequence either converges to a genuine extremal or concentrates at a point with a universal energy profile (the standard bubble).

\subsection{Passing to the Sub-Riemannian Setting: Why the Heisenberg Group?} In a parallel development, and largely independently motivated, there has been a sustained effort to carry these embedding theorems from the Euclidean elliptic setting to the sub-elliptic world of Carnot groups, and in particular to the simplest and most studied of these, the Heisenberg group $\mathbb{H}^n$. Why $\mathbb{H}^n$? Because it is far more than a toy model. It appears naturally as the tangent space (in the Gromov sense) to any nondegenerate CR manifold, and it models the boundary geometry of strictly pseudoconvex domains in $\mathbb{C}^{n+1}$. Thus, sub-elliptic PDEs on $\mathbb{H}^n$ are exactly what we need for CR-geometric problems, the CR Yamabe problem, the Kohn Laplacian $\Box_b$, the Webster scalar curvature, in much the same way that Euclidean PDEs govern conformal geometry on the sphere.

Let us recall the setup. Denote $\xi=(z,t)=(x,y,t)$ as a generic point of $\mathbb{R}^{2n+1}$. The Heisenberg group $\mathbb{H}^{n}$ is the homogeneous non-abelian Lie group $(\mathbb{R}^{2n+1},\circ)$, whose group law is given by
$$
\xi \circ \xi^\prime:= \big(z + z',\; t + t' + 2\big((y, x') - (x, y')\big)\big)\quad\text{for all}\quad \xi,\,\xi^\prime\in \mathbb{H}^{n},$$
where $\xi=(x,y,t)$, $\xi^\prime=(x^\prime,y^\prime,t^\prime)$ and $(\cdot,\cdot)$ denotes the inner product on $\mathbb{R}^{n}$. In addition, there hold $\xi^{-1}=-\xi$ and $(\xi\circ \xi^\prime)^{-1}=(\xi^\prime)^{-1} \circ\,\xi^{-1}$. For a fixed $\eta\in\mathbb{H}^{n}$, the left translation by $\eta$ is the mapping $\tau_{\eta}\colon\mathbb{H}^n\to\mathbb{H}^n$ which, according to the group law $\circ$, is defined by $\tau_\eta(\xi)=\eta\circ\xi$ for all $\xi\in\mathbb{H}^n$. Its Lie algebra $\mathfrak{h} = V_1 \oplus V_2$ is generated by the left-invariant vector fields
$$X_j := \partial_{x_j} + 2y_j \partial_t\quad\text{and}\quad Y_j := \partial_{y_j} - 2x_j \partial_t\quad\text{for}\quad j=1,\dots,n,$$
where $V_1=\operatorname{Span}\{X_j,Y_j\}_{j=1}^n$ and $V_2 = \mathbb{R} T$, with $T = \partial_t$. The only nontrivial commutator is $[X_j, Y_j] = -4T$. The vanishing of all commutators of length greater than two shows that $\mathbb{H}^n$ is a nilpotent, graded, and stratified Lie group of step two. The left-invariant vector fields in $\operatorname{Span}\{X_j,Y_j\}_{j=1}^n$ are called horizontal vector fields. The natural dilation $\delta_R\colon\mathbb{H}^n\to\mathbb{H}^n$, defined by
$ \delta_R(\xi)=(Rz,R^2t)$ for all $R>0$, have Jacobian determinant $R^{2n+2}$. Thus, the natural number $Q=2n+2$ is called the homogeneous dimension of $\mathbb{H}^n$. It plays a role analogous to that of the topological dimension in the Euclidean setting and is the relevant dimension for analysis on $\mathbb{H}^n$, rather than the topological dimension $2n+1$. The function $r\colon\mathbb{H}^n\to\mathbb{R}$ is called the Kor\'{a}nyi norm, which is defined by $r(\xi):=(|z|^4+t^2)^{\frac{1}{4}}$ for all $\xi\in \mathbb{H}^n$, where $|z|$ is the natural Euclidean norm in $\mathbb{R}^{2n}$. It satisfies $r(\delta_R(\xi))=R r(\xi)$ for all $R>0$. The associated Kor\'{a}nyi distance function $d_K\colon \mathbb{H}^n\times \mathbb{H}^n\to\mathbb{R}$ is given by $d_K(\xi,\xi'):=r(\xi^{-1}\circ\xi')$ for all $\xi,\,\xi^\prime\in\mathbb{H}^n$. The Kor\'{a}nyi open ball of radius $R$ and centered at $\xi_0$ is defined as $B_R(\xi_0):=\{\xi\in\mathbb{H}^n\colon\,d_K(\xi,\xi_0)<R\}$. In the sequel, $B_R$ will denote $B_R(O)$, where $O=(0,0)$ denotes the identity element of the Heisenberg group $\mathbb{H}^n$.

The group $\mathbb{H}^n$ is naturally equipped with the Haar measure induced by the Lebesgue measure on $\mathbb{R}^{2n+1}$. Indeed, the latter is invariant under left translations, and we denote this measure by $\mathrm{d}\xi$. For every measurable set $\mathcal{E}\subset\mathbb{H}^n$, we write $|\mathcal{E}|$ for its $(2n+1)$-dimensional Lebesgue measure. The Haar measure is homogeneous of degree $Q=2n+2$ under the dilations $\delta_R$. It follows that
$$
|\delta_R(\mathcal{E})|=R^Q|\mathcal{E}|,
\quad
\mathrm{d}(\delta_R(\xi))=R^Q\mathrm{d}\xi,\quad \text{and}\quad |B_R|=R^Q|B_1|.
$$

Let $u\in C^1(\mathbb{H}^n,\mathbb{R})$. The subelliptic (horizontal) gradient of $u$ is the $2n$-dimensional vector field defined by $D_H u := (X_1 u, \dots, X_n u, Y_1 u, \dots, Y_n u),$
Moreover, we have
$$|D_H u|_H=\sqrt{(D_H u,D_H u)_H}=\bigg(\sum_{j=1}^{n} (X_ju)^2+(Y_ju)^2\bigg)^{\frac{1}{2}},$$
where $(\cdot,\cdot)_H$ is the inner product and $|\cdot|_H$ is the corresponding norm on $\mathbb{R}^{2n}$. Note that, if $g\in C^1(\mathbb{R)}$, then there holds $D_H g(u) =g^\prime(u) D_H u.$
Let $A\colon\xi\in\mathbb{H}^n\to A(\xi)=(A_1\xi,\cdots,A_{2n}\xi)$ be a $C^1$ vector field. Then the horizontal divergence of $A$ is given by
$$\text{div}_H A(\xi):=\sum_{j=1}^{n} X_j(A_j\xi)+\sum_{j=1}^{n} Y_j(A_{n+j}\xi)\quad\text{for all}\quad\xi \in\mathbb{H}^n.$$ The sub-Laplacian is the second-order differential operator defined by $\Delta_{H}:= \sum_{j=1}^n (X_j^2 + Y_j^2),$ which is the natural replacement for the Laplacian; it is hypoelliptic, as established by H\"{o}rmander \cite{Hormander-1967}, but not elliptic, since it omits the $T$ direction, which is recovered only through the non-commutativity. Recall that $X_j$ and $Y_j$ are left-invariant under the group operation and are homogeneous of degree one with respect to the dilation $\delta_R$. It follows that $\Delta_{H}$ is left translation invariant and $\delta_R$-homogeneous of degree two. To define the notation of horizontal Sobolev spaces, we fixed $p\geq1$.
Then the first order horizontal Sobolev space $HW^{1,p}(\mathbb{H}^n)$ is given by
$HW^{1,p}(\mathbb{H}^n):=\big\{ u\in L^p(\mathbb{H}^n)\colon\,|D_H u|_H\in L^p(\mathbb{H}^n)\big\},$
with its natural norm
$$ \|u\|:=\bigg(\int_{\mathbb{H}^n}|D_H u|_H^p\,\mathrm{d}\xi+\int_{\mathbb{H}^n}|u|^p\,\mathrm{d}\xi\bigg)^{\frac{1}{p}}.$$
In particular, $HW^{1,p}(\mathbb{H}^n)$ is the completion of $C^\infty_0(\mathbb{H}^n)$ with respect to this norm. For later use, given $\tau>0$, we introduce the equivalent norm
$$ \|u\|_\tau:=\bigg(\int_{\mathbb{H}^n}|D_H u|_H^p\,\mathrm{d}\xi+\tau\int_{\mathbb{H}^n}|u|^p\,\mathrm{d}\xi\bigg)^{\frac{1}{p}}.$$
It is well-known that the space $HW^{1,p}(\mathbb{H}^n)$ is a uniformly convex and separable Banach space for all $p\in(1,+\infty)$. Next, let $\Omega\subset\mathbb{H}^n$ be any smooth bounded domain. Then the horizontal Sobolev space $HW^{1,p}_0(\Omega)$ is defined as the completion of $C^\infty_0(\Omega)$ with respect to the Dirichlet norm given by
$$\|D_H u\|_{L^p(\Omega)}:= \left(\int_\Omega |D_H u|_H^p\,\mathrm{d}\xi\right)^{\frac{1}{p}}.$$

In what follows, $\|\cdot\|_p$ denotes the usual $L^{p}(\mathbb{H}^n)$-norm for $p\in[1,+\infty]$, while $\|\cdot\|_{L^p(\Omega)}$ denotes the corresponding norm on $L^{p}(\Omega)$ for any measurable set $\Omega\subset \mathbb{H}^n$. In addition, for $p\in[1,+\infty)$ and $0\leq\beta<Q$, the weighted Lebesgue space $L^p(\mathbb{H}^n,r(\xi)^{-\beta}\,\mathrm{d}\xi)$ consists of real-valued measurable function $u$ such that $|u|^p r(\xi)^{-\beta}\in L^1(\mathbb{H}^n)$. Its norm is defined by
$$ \|u\|_{p,\beta}:=\left(\int_{\mathbb{H}^n} |u|^p{r(\xi)^{-\beta}}\,\mathrm{d}\xi\right)^{\frac{1}{p}}.$$
\subsubsection{\textnormal{\bf{Literature Survey on the Moser--Trudinger Inequality on the Heisenberg Group
}}} The celebrated Folland--Stein \cite{Folland-Stein} theory of Sobolev spaces on stratified groups asserts that, for any bounded domain $\Omega\subset\mathbb{H}^n$ and $1\leq p<Q$, the embedding $HW_0^{1,p}(\Omega) \hookrightarrow L^{p^*}(\Omega)$ holds, where $p^\ast:=\frac{Qp}{Q-p}$ denotes the critical Sobolev exponent. 
The limiting case $p=Q$ marks the point at which the usual Sobolev embedding gives way to a genuinely critical phenomenon, and naturally leads to the question:
  \[ \textit{What happens at the Sobolev borderline case $p = Q$? Does the sharp exponential}\]
  \[\textit{integrability have a perfect analogue in the sub-Riemannian setting?}\]
This question was fully answered by Cohn and Lu \cite{Cohn-Lu-2001}. Indeed, they established the sharp exponential integrability result on domains of finite measure, with the sharp constant characterized as follows.
\begin{knownthm}[Cohn--Lu \cite{Cohn-Lu-2001}]\label{knownthm}
For every bounded domain $\Omega\subset\mathbb{H}^n$, there exists a constant $C_0>0$ such that
$$
\sup_{\substack{u \in HW_0^{1,Q}(\Omega),\, \|D_H u\|_{L^Q(\Omega)} \le 1}} \; \frac{1}{|\Omega|} \int_\Omega \exp\big( \alpha_Q |u|^{\frac{Q}{Q-1}} \big) \mathrm{d}\xi\;\le\; C_0,
$$
where the sharp constant is $\alpha_Q = Q\, c_Q^{\frac{1}{Q-1}}$ with
$c_Q=\int_{\{\xi\in \mathbb{H}^n\colon~r(\xi)=1\}}|z|^Q\,\mathrm{d}\mu
$. Moreover, the exponent $\alpha_Q$ is sharp: replacing $\alpha_Q$ by any $\alpha>\alpha_Q$ makes the supremum infinite.
\end{knownthm}

This theorem is the exact sub-Riemannian analogue of Moser's classical inequality, but its proof reveals a genuinely different mechanism dictated by the geometry of $\mathbb{H}^n$. The failure of the classical P\'{o}lya--Szeg\"{o} inequality for the horizontal gradient motivated Cohn--Lu \cite{Cohn-Lu-2001}, following Adams \cite{Adams-1988}, to develop a potential-theoretic method for proving the sharp Moser--Trudinger inequality on the Heisenberg group. Balogh et al. \cite{BMT-2003} subsequently extended this result to general Carnot groups, while sharp weighted Moser--Trudinger inequalities on Heisenberg-type groups, both on bounded domains and on the whole space, were obtained in \cite{Lam-Tang-2013}.

The weighted counterpart of the Moser--Trudinger inequality was motivated by the seminal contribution of Adimurthi--Sandeep \cite{Adi-Sandeep-2007} and further generalized by Lam--Lu--Tang \cite{Lam-Lu-Tang-2012} into the Heisenberg setting. If $\Omega\subset\mathbb{H}^n$ is bounded and $0\leq\beta<Q$, then there exists $C=C(Q,\beta)>0$ such that
 \begin{align*}
    \sup_{u\in HW^{1,Q}_0(\Omega),\,\|D_H u\|_{L^Q(\Omega)\leq1}}\frac{1}{|\Omega|^{1-\frac{\beta}{Q}}}\int_\Omega\frac{\operatorname{\exp}\Big(\big(1-\frac{\beta}{Q}\big)\alpha_Q|u|^{\frac{Q}{Q-1}}\Big)}{r(\xi)^\beta}\,\mathrm{d}\xi\leq C.
\end{align*}
The exponent $\alpha_Q$ is sharp, as the above supremum ceases to be finite once $\alpha_Q$ is replaced by any $\alpha>\alpha_Q$.

The inequality of Adimurthi--Yang \cite{Adi-Yang-2010} admits a weighted extension to $\mathbb{H}^n$, established by Lam--Lu \cite{Lam-Lu-2012} without recourse to symmetrization. For any $\tau>0$, $0\leq\beta<Q$, and $0<\alpha\leq\alpha_Q$, they proved that
\begin{align*}
 \sup_{u\in HW^{1,Q}(\mathbb{H}^n),~\|u\|_{\tau}\leq 1} ~\int_{\mathbb{H}^n}\cfrac{\Phi\Big(\big(1-\frac{\beta}{Q}\big)\alpha|u|^{\frac{Q}{Q-1}}\Big)}{r(\xi)^\beta}\,\mathrm{d}\xi<+\infty,
 \end{align*}
where $\alpha_Q=Q\left(2\pi^n\Gamma\left(\frac{1}{2}\right)\Gamma\left(\frac{Q-1}{2}\right){\Gamma\left(\frac{Q}{2}\right)}^{-1}{\Gamma(n)}^{-1}\right)^{\frac{1}{Q-1}}$ and $\Phi\colon[0,+\infty)\to [0,+\infty)$ is given by 
\begin{align*}
   \Phi(s):=\exp(s)-\sum_{j=0}^{Q-2}\frac{s^j}{j!}\quad\text{for all}\quad s\geq0,
\end{align*} which is also known as the Young function.
The function defined above will be used below in the statements and proofs of the relevant weighted Moser--Trudinger inequalities on the Heisenberg group.
In particular, if $\alpha>\alpha_{Q}$, then the integral is still finite, but the supremum diverges to infinity.

Lam--Lu--Tang \cite{Lam-Lu-Tang-2014} provided a detailed account of the subcritical weighted Moser--Trudinger inequality. The potential-dependent extension of this inequality was later established by Chen--Lu--Zhu \cite{Chen-Lu-Zhu-2023}. Cohn--Lam--Lu--Yang \cite{Cohn-Lam-Lu-Yang-2012} pursued a symmetrization approach on $\mathbb{H}^n$ and, with the aid of the modified P\'{o}lya--Szeg\"{o} inequality in \cite[Theorem 5.3]{JJM2019} as a key ingredient, obtained an alternative form. Their result states that, for $0<\alpha\leq\alpha^\ast$ and $0\leq\beta<Q$, one has
 \begin{align*}
 \sup_{u\in HW^{1,Q}(\mathbb{H}^n),~\|u\|\leq 1} ~\int_{\mathbb{H}^n}\cfrac{\Phi\big(\alpha|u|^{\frac{Q}{Q-1}}\big)}{r(\xi)^\beta}\,\mathrm{d}\xi<+\infty\iff \frac{\alpha}{\alpha^\ast}+\frac{\beta}{Q}\leq 1,
 \end{align*}
 where $\alpha^{\ast}=\frac{\alpha_Q}{c^\ast}$, $u^\ast$ is the decreasing rearrangement of the nonnegative function $u$, and $c^\ast$ is given by $$c^\ast=\inf\bigg\{ c^{\frac{1}{Q-1}}\colon\,\displaystyle\int_{\mathbb{H}^n} |D_H u^\ast|_H^Q\,\mathrm{d}\xi\leq c \displaystyle\int_{\mathbb{H}^n} |D_H u|_H^Q\,\mathrm{d}\xi,\, u\in HW^{1,Q}(\mathbb{H}^n)\bigg\}.$$

Later, Li--Lu--Zhu \cite{Li-Lu-Zhu-2018} established, by a completely symmetrization-free method, an analog of the Lions concentration-compactness principle \cite{Lions-1985} on both bounded domains and the whole space $\mathbb{H}^n$. In the bounded-domain case, let $\Omega\subset\mathbb{H}^n$ be bounded and let $0\leq\beta<Q$. If there exists a sequence $\{u_k\}_k\subset HW^{1,Q}_0(\Omega)$ with $\|D_H u_k\|_{L^Q(\Omega)}=1$ and $u_k\rightharpoonup u\neq 0$ in $HW^{1,Q}_0(\Omega)$ as $k\to\infty$, then
$$\sup_{k\in\mathbb{N}}\int_\Omega\frac{\operatorname{\exp}\Big(\big(1-\frac{\beta}{Q}\big)\alpha_{Q}\ell|u_k|^{\frac{Q}{Q-1}}\Big)}{r(\xi)^\beta}\,\mathrm{d}\xi<+\infty,$$
where $0<\ell<\zeta_{Q}(u):=\big(1-\|D_H u\|^Q_{L^Q(\Omega)}\big)^{-\frac{1}{Q-1}}$. In addition, the value $\zeta_Q(u)$ is the sharp threshold: the above supremum is infinite for all $\ell\geq\zeta_Q(u)$. The whole-space counterpart is obtained by considering a sequence $\{u_k\}_k\subset HW^{1,Q}(\mathbb{H}^n)$ with $\| u_k\|=1$ and $u_k\rightharpoonup u\neq 0$  in $HW^{1,Q}(\mathbb{H}^n)$ as $k\to\infty$. It states that for $0\leq\beta<Q$, there holds
$$\sup_{k\in\mathbb{N}}\int_{\mathbb{H}^n}\cfrac{\Phi\Big(\big(1-\frac{\beta}{Q}\big)\alpha_{Q}\ell|u_k|^{\frac{Q}{Q-1}}\Big)}{r(\xi)^\beta}\,\mathrm{d}\xi<+\infty,$$
where $0<\ell<\tilde{\zeta}_{Q}(u):=\big(1-\| u\|^Q\big)^{-\frac{1}{Q-1}}$. This threshold is sharp, as the supremum is infinite for all $\ell>\tilde{\zeta}_Q(u)$. The subsequent literature includes the contributions of Chen--Lu--Zhu \cite{Chen-Lu-Zhu-2021}, Figueiredo et al. \cite{Figueiredo-Mahanta-Sarkar-2027} and Li--Lu--Zhu \cite{Li-Lu-Zhu-2021}.

In addition, the maximizers for Moser--Trudinger inequalities in the subcritical case have been investigated by Ishiwata et al. \cite{Ishiwata-Nakamura-Wadade-2014} and Cassani et al. \cite{Cassani-Sani-Tarsi-2014}; related problems subject to inhomogeneous norm constraints were treated by Ikoma et al. \cite{Ikoma-Ishiwata-Wadade-2019}, Lam--Lu--Zhang \cite{Lam-Lu-Zhang-2019} and Nguyen \cite{Nguyen-2021}. The equivalence between critical and subcritical Moser--Trudinger inequalities was first established in the planar case by Cassani et al. \cite{Cassani-Sani-Tarsi-2014}, which was later extended to higher dimensions by Lam--Lu--Zhang \cite{Lam-Lu-Zhang-2017}. This equivalence further facilitated the study of maximizers for Moser--Trudinger inequalities under inhomogeneous norm constraints, as pursued in \cite{Lam-2017, Lam-Lu-Zhang-2019}. Although Lam--Lu \cite{Lam-Lu-2012} established the sharp critical Moser--Trudinger inequality on the Heisenberg group, the corresponding maximizer problem is complicated by the absence of a suitable P\'{o}lya--Szeg\"{o} inequality for the subelliptic gradient; see Jerison--Lee \cite{JL1988} and Manfredi et al. \cite{JJM2019}. In a recent contribution, Figueiredo et al. \cite{Figueiredo-Mahanta-Sarkar-2027} established the equivalence between critical and subcritical Moser--Trudinger-type inequalities in the subspace of the borderline Sobolev space on the Heisenberg group.

The Moser--Trudinger theory has also been extensively developed in complex, CR, and conformal geometries. Cohn--Lu \cite{Cohn-Lu-2004} established sharp Moser--Trudinger inequalities on complex spheres in $\mathbb{C}^n$ for smooth and holomorphic functions. Gui--Moradifam \cite{Changfeng-Moradifam-2018} established the optimal constant in a Moser--Trudinger-type inequality on the sphere, confirming a conjecture that was proposed by Chang--Yang \cite{Chang-Yang-1987}. In addition, De Marchis et al. \cite{DMMT-2022} studied the critical points of the Moser--Trudinger functional on closed Riemann surfaces, proving the existence of positive critical points through a min--max approach combined with compactness, quantization, and sharp energy estimates. 

Much of the Moser--Trudinger theory on the Heisenberg group described above is due to Lu and his collaborators. This includes the sharp inequality on bounded domains \cite{Cohn-Lu-2001} and its weighted version \cite{Lam-Lu-Tang-2012}, the sharp critical and subcritical inequalities on the whole group \cite{Cohn-Lam-Lu-Yang-2012, Lam-Lu-2012, Lam-Lu-Tang-2014}, the concentration-compactness principle \cite{Li-Lu-Zhu-2018, Li-Lu-Zhu-2021}, the sharpened inequality of Adimurthi--Druet type \cite{Chen-Lu-Zhu-2021}, and the existence of least energy solutions to the corresponding subelliptic equations \cite{Chen-Lu-Zhu-2023}. Lu and his collaborators also introduced the rearrangement-free argument for sharp Moser--Trudinger and Adams inequalities on unbounded domains \cite{Lam-Lu-2013} and proved the equivalence between the critical and subcritical inequalities \cite{Lam-Lu-Zhang-2017}. These two ideas are the starting point of the present work.

\subsection{Notations}
 In the entire article, $(X^\ast,\|\cdot\|_{\ast})$ stands for the continuous dual of the Banach space $(X,\|\cdot\|)$, and $\langle\cdot,\cdot\rangle$ the duality pairing between $X^\ast$ and $X$. We write $\chi_{\mathcal{E}}$ for the characteristic function of $\mathcal{E}\subset\mathbb{H}^n$, $Q^\prime=\frac{Q}{Q-1}$ for the H\"{o}lder conjugate of $Q$, and $o_k(1)\to0$ as $k\to\infty$. The symbols $\rightharpoonup$ and $\to$ denote weak and strong convergences, respectively. Moreover, $\lesssim$ and $\gtrsim$ denote inequalities up to a positive constant, $O(f)$ a quantity bounded by a constant multiple of $f$, and $\underset{\approx}{<}$ and $\thickapprox$ indicate approximate closeness and equality, respectively.
\subsection{{Highlights of the Main Results}} Now, the results established in this article are presented in the sequel. The discussion begins with a weighted Moser--Trudinger inequality of Adimurthi--Druet-type, together with the determination of its sharp constant.
\begin{theorem}[\textbf{The Adimurthi--Druet-Type Weighted Moser--Trudinger Inequality}]\label{thm1.1}
Let $0\leq \alpha<1$, $0\leq\beta<Q$ and $\tau>0$ be real numbers. Then there holds
\begin{equation}
\begin{aligned}\label{eq1.2}
    \textnormal{ADMT}(\alpha,\beta,\tau):&=\sup_{u\in HW^{1,Q}(\mathbb{H}^n),\|u\|_\tau\leq 1}\int_{\mathbb{H}^n}\frac{\Phi\big(\alpha_Q \big(1-\frac{\beta}{Q}\big)\big(1+\alpha \tau \|u\|_Q^Q \big)^{\frac{1}{Q-1}}|u|^{Q^\prime}\big)}{r(\xi)^\beta}\,\mathrm{d}\xi\\&<+\infty. 
\end{aligned}
\end{equation}
In addition, if $\alpha\geq1$, then the above supremum $\textnormal{ADMT}(\alpha,\beta,\tau)$ ceases to be finite.
\end{theorem}
In what follows, we revisit the subcritical weighted Moser--Trudinger inequality initially established by Lam--Lu--Tang \cite{Lam-Lu-Tang-2014} and prove precise asymptotic bounds for the corresponding supremum near the sharp threshold.
\begin{theorem}[\textbf{The Subcritical Weighted Moser--Trudinger Inequality}]\label{thm1.3}
Assume that $0<\alpha<\alpha_Q$ and $0\leq\beta<Q$. Then there holds
$$\textnormal{SCMT}(\alpha,\beta):=\sup_{u\in HW^{1,Q}(\mathbb{H}^n),~\|D_H u\|_{Q}\leq 1}\cfrac{1}{\|u\|_Q^{Q-\beta}} \int_{\mathbb{H}^n}\frac{\Phi\big(\alpha \big(1-\frac{\beta}{Q}\big)|u|^{Q^\prime}\big)}{r(\xi)^\beta}\,\mathrm{d}\xi<+\infty.$$
Moreover, the constant $\alpha_Q$ is sharp in the sense that $\textnormal{SCMT}(\alpha_Q,\beta)=+\infty$. Also, there exist two  positive constants $c(Q,\beta)$ and $C(Q, \beta)$, which depend only on $Q$ and $\beta$ such that for $\alpha\underset{\approx}{<}\alpha_Q$, the following inequality holds
$$ \cfrac{c(Q,\beta)}{\Big(1-\big(\frac{\alpha}{\alpha_Q}\big)^{Q-1}\Big)^{1-\frac{\beta}{Q}}}\leq \textnormal{SCMT}(\alpha,\beta)\leq \cfrac{C(Q,\beta)}{\Big(1-\big(\frac{\alpha}{\alpha_Q}\big)^{Q-1}\Big)^{1-\frac{\beta}{Q}}}.$$
\end{theorem}
The following result establishes the existence of a maximizer for the weighted subcritical Moser--Trudinger inequality under a suitable restriction on the parameter.
\begin{theorem}\label{thm1.4}
Let $0<\alpha<\alpha_Q$ and $0<\beta<Q$ be real numbers. Then there exists a maximizer for $\textnormal{SCMT}(\alpha,\beta)$.
\end{theorem}
The preceding subcritical weighted Moser--Trudinger inequality also yields a Caffarelli--Kohn--Nirenberg interpolation-type inequality in the asymptotic sense, providing a further consequence within the borderline Sobolev framework.
\begin{theorem}[\textbf{The Caffarelli--Kohn--Nirenberg Interpolation-Type Inequality}]\label{thm1.5}
Let $0\leq\beta<Q$ and $\vartheta\geq Q$ be real numbers. Then there holds
\begin{align}\label{eq1.3}
  \|u\|_{\vartheta,\beta}\leq C(Q,\beta,\vartheta) \|u\|_Q^{\frac{Q-\beta}{\vartheta}}\|D_Hu\|_Q^{1-\frac{Q-\beta}{\vartheta}}\quad\text{for all}\quad u\in HW^{1,Q}(\mathbb{H}^n),
\end{align}
where $C(Q,\beta,\vartheta)>0$ is a constant  depending only on $Q$, $\beta$ and $\vartheta$. Define the asymptotic best constant $\boldsymbol{\Lambda}_{Q,\beta}$ corresponding to the above inequality by
\begin{align*}
\boldsymbol{\Lambda}_{Q,\beta}:=\limsup_{\vartheta\to\infty} \sup_{u\in HW^{1,Q}(\mathbb{H}^n)\setminus\{0\}} \Bigg(\cfrac{\|u\|_{\vartheta,\beta}}{\vartheta^{\frac{1}{Q^\prime}}\|u\|_Q^{\frac{Q-\beta}{\vartheta}}\|D_Hu\|_Q^{1-\frac{Q-\beta}{\vartheta}}}\Bigg).
\end{align*}
Notably, if we assume $\lambda>\boldsymbol{\Lambda}_{Q,\beta}$, then there exists a positive constant $\vartheta_0=\vartheta_0(Q,\beta,\lambda)\geq Q$ depending only on $Q$, $\beta$ and $\lambda$ such that for any $\vartheta\geq\vartheta_0$, there holds
\begin{align*}
  \|u\|_{\vartheta,\beta}\leq \lambda\vartheta^{\frac{1}{Q^\prime}} \|u\|_Q^{\frac{Q-\beta}{\vartheta}}\|D_Hu\|_Q^{1-\frac{Q-\beta}{\vartheta}}\quad\text{for all}\quad u\in HW^{1,Q}(\mathbb{H}^n).
\end{align*}
Moreover, the asymptotic best constant $\boldsymbol{\Lambda}_{Q,\beta}$ admits an explicit expression, which is given by $\boldsymbol{\Lambda}_{Q,\beta}=\big(\alpha_Q \big(1-\frac{\beta}{Q}\big)Q^\prime\,\operatorname{\exp}(1)\big)^{-\frac{1}{Q^\prime}}$.
\end{theorem}
In a complementary direction, a weighted Moser--Trudinger inequality under an inhomogeneous norm constraint is established and shown to be equivalent to the subcritical weighted Moser--Trudinger inequality. The precise statement is as follows.
\begin{theorem}[\textbf{The Weighted Moser--Trudinger Inequality under Inhomogeneous Norm Constraints}]\label{thm1.6}
Let $\alpha>0$, $0\leq\beta<Q$ and $a,b>0$ be real numbers. Then there holds
$$\textnormal{CMT}_{a,b}(\alpha,\beta):=\sup_{u\in HW^{1,Q}(\mathbb{H}^n),~\|D_H u\|^a_{Q}+\| u\|^b_Q\leq 1}\int_{\mathbb{H}^n}\frac{\Phi\big(\alpha \big(1-\frac{\beta}{Q}\big)|u|^{Q^\prime}\big)}{r(\xi)^\beta}\,\mathrm{d}\xi<+\infty $$
if and only if either $\alpha<\alpha_Q$ or $\alpha=\alpha_Q$ with $b\leq Q$. The constant $\alpha_Q$ is sharp. In addition, the following identity
\begin{align}\label{eq1.4}
  \textnormal{CMT}_{a,b}(\alpha,\beta)=\sup_{\zeta\in(0,\alpha)}\Bigg(\cfrac{1-\big(\frac{\zeta}{\alpha}\big)^{(\frac{Q-1}{Q})a}}{\big(\frac{\zeta}{\alpha}\big)^{(\frac{Q-1}{Q})b}}\Bigg)^{\frac{Q}{b}(1-\frac{\beta}{Q})}\textnormal{SCMT}(\zeta,\beta)
\end{align}
holds in both cases $\alpha<\alpha_Q$ and $\alpha=\alpha_Q$ with $b\leq Q$. More precisely, there hold $\textnormal{CMT}_{Q,Q}(\alpha_Q,\beta)<+\infty $
and
\begin{align*}\textnormal{CMT}_{Q,Q}(\alpha_Q,\beta)=\sup_{\zeta\in(0,\alpha_Q)}\Bigg(\cfrac{1-\big(\frac{\zeta}{\alpha_Q}\big)^{Q-1}}{\big(\frac{\zeta}{\alpha_Q}\big)^{Q-1}}\Bigg)^{1-\frac{\beta}{Q}}\textnormal{SCMT}(\zeta,\beta).
\end{align*}
In particular, the above identity illustrates the equivalence between the critical and subcritical weighted Moser--Trudinger inequalities.
\end{theorem}
The established equivalence \eqref{eq1.4} naturally raises the question of whether the corresponding supremum is attained in the subcritical regime. The answer is given in the following theorem, which characterizes the existence of a maximizer.
\begin{theorem}\label{thm1.7}
For all $0<\beta<Q$ and $a,b>0$, there exists a maximizer for $\textnormal{CMT}_{a,b}(\alpha,\beta)$ in the subcritical regime: either $(\alpha<\alpha_Q,\, b>0)$ or $(\alpha=\alpha_Q,\, b<Q)$.
\end{theorem}
In the same spirit, the subcritical weighted Moser--Trudinger inequality can be further refined, leading to the following result.
\begin{theorem}[\textbf{The Improved Subcritical Weighted Moser--Trudinger Inequality}]\label{thm1.8}
 Assume that $0<\alpha<\alpha_Q$ and $0\leq\beta<Q$. Then there holds
 \begin{align*}
  \textnormal{SCMT}_{\operatorname{e}}(\alpha,\beta):=\sup_{u\in HW^{1,Q}(\mathbb{H}^n),~\|D_H u\|_{Q}\leq 1}\cfrac{1}{\|u\|_Q^{Q-\beta}} \int_{\mathbb{H}^n}\frac{\operatorname{\exp}\big(\alpha \big(1-\frac{\beta}{Q}\big)|u|^{Q^\prime}\big)|u|^Q}{r(\xi)^\beta}\,\mathrm{d}\xi<+\infty.
 \end{align*}
In addition, the constant $\alpha_Q$ is sharp in the sense that the above supremum is infinite for all $\alpha\geq \alpha_Q$.
\end{theorem}
Next, turning to the associated variational problem, the preceding inequality yields a maximization problem whose attainment in a suitable regime is characterized by the following theorem.
\begin{theorem}\label{thm1.9}
For all $0<\alpha<\alpha_Q$ and $0<\beta<Q$, there exists a maximizer for $\textnormal{SCMT}_{\operatorname{e}}(\alpha,\beta)$.
\end{theorem}
\subsection{Novelty and Contributions of the Present Work} The results obtained here enrich the existing theory of Moser--Trudinger inequalities on the Heisenberg group in several directions. The key novelties and technical challenges are summarized below.
\begin{itemize}
 \item[$\triangleright$] The Adimurthi--Druet-type inequality, rooted in the concentration--compactness principle of P.-L. Lions \cite{Lions-1985}, has evolved differently in Euclidean and Heisenberg settings. In Euclidean spaces, the results of Adimurthi--Druet \cite{Adi-Druet2004} and Yang \cite{Yang-2006} rely on delicate blow-up analysis. On the Heisenberg group, Li--Lu--Zhu \cite{Li-Lu-Zhu-2018} established the concentration--compactness principle through rearrangement-free techniques, followed by the Adimurthi--Druet-type improvement of Chen--Lu--Zhu \cite{Chen-Lu-Zhu-2021}. Motivated by these works, we establish its weighted counterpart under a weighted-norm constraint.
\item[$\triangleright$] The existence and nonexistence of maximizers for variational problems associated with Moser--Trudinger inequalities has been extensively studied in the Euclidean setting, with particular emphasis on the critical case; see, among others, Ishiwata \cite{Ishiwata-2011} and Nguyen \cite{Nguyen-2021,Nguyen-2021@}. The classical analysis often relies on symmetric decreasing rearrangement, together with the P\'{o}lya--Szeg\"{o} inequality, which provides a fundamental tool for reducing the variational problem to radially symmetric functions while preserving the Sobolev norm. On the Heisenberg group, the maximization problems associated with Moser--Trudinger inequalities are complicated by the absence of a suitable P\'{o}lya--Szeg\"{o} inequality for the subelliptic gradient. The resulting lack of radial reduction prevents a direct application of the concentration-compactness and blow-up machinery. In this work, we partially address this problem by proving the existence of a maximizer for weighted Moser--Trudinger inequalities, both in the subcritical regime and under inhomogeneous norm constraints on $\mathbb H^n$, without relying on radial rearrangement. The key ingredient is a compact embedding of the borderline Sobolev space into a suitable weighted Lebesgue space on $\mathbb H^n$.
\item[$\triangleright$] In the borderline Sobolev regime, we establish new Caffarelli--Kohn--Nirenberg interpolation-type inequalities on the Heisenberg group, with a parameter range extending the corresponding Euclidean results of Caffarelli et al. \cite{CKN-1984}. However, for homogeneous groups, related inequalities have been established by Ozawa et al. \cite{Ozawa-Ruzhansky-Suragan-2019}, and our results further extend the Caffarelli--Kohn--Nirenberg theory on the Heisenberg group. This provides a new extension of the Caffarelli--Kohn--Nirenberg framework to the context of the Heisenberg group (for the special case in \cite{Ozawa-Ruzhansky-Suragan-2019} with $\alpha=0=\beta$, the CKN inequality only holds for $p=Q$). In addition, we identify the relationship between the sharp constants of the subcritical weighted Moser--Trudinger inequalities and those of the Caffarelli--Kohn--Nirenberg interpolation-type inequalities in the asymptotic sense.
\item[$\triangleright$] The present results extend beyond the weighted Moser–Trudinger theory to address normalized $Q$-sub-Laplacian problems with singular critical exponential nonlinearities, a previously open area in the Heisenberg setting due to the absence of suitable Caffarelli–Kohn–Nirenberg inequalities. The challenges, especially under prescribed-mass constraints, arise from balancing $Q$-Dirichlet energy, lower-order contributions, and singular nonlinearity. The newly established Caffarelli–Kohn–Nirenberg interpolation-type inequalities, along with subcritical weighted Moser–Trudinger inequalities, provide the analytical framework to tackle these issues.
\end{itemize}
\subsection{From the Derived Inequalities to Subelliptic PDEs}
The study of subelliptic equations with critical exponential growth on the Heisenberg group has developed in several directions, including critical $Q$-sub-Laplacian equations, singular $(p,Q)$-sub-Laplacian equations, and their nonlocal counterparts with Stein--Weiss-type interactions; see, for instance, \cite{Chen-Lu-Zhu-2023, Figueiredo-Mahanta-Sarkar-2027, Lam-Lu-2012, Lam-Lu-Tang-2012, Li-Lu-Zhu-2018, DTAL-2026, PT-2022} and the references therein. However, normalized problems driven by the $Q$-sub-Laplacian have not yet been studied. The present work addresses this problem.

An application of the established inequalities is illustrated by the study of nontrivial normalized ground states for a class of quasilinear subelliptic equations on the entire Heisenberg group, driven by the $Q$-sub-Laplacian, which brings a more pronounced “PDE flavor" to the analysis. To this end, we consider the following problem:
\begin{equation}\label{main problem}
-\Delta_{H,Q} u+\lambda|u|^{Q-2}u
=\frac{g(u)}{r(\xi)^\beta}
\quad\text{in}\quad \mathbb{H}^n,
\quad u\in HW^{1,Q}(\mathbb{H}^n),
\tag{$\mathrm{P}_a$}
\end{equation}
subject to the prescribed mass constraint
\begin{align*}
\int_{\mathbb{H}^n}|u|^Q\,\mathrm{d}\xi=a^Q,
\end{align*}
where $0<\beta<Q$, $a>0$ is prescribed, and $\lambda\in\mathbb{R}$ acts as the Lagrange multiplier. The nonlinearity $g\colon \mathbb{R}\to \mathbb{R}$ is a continuous function, which has singular critical exponential growth at infinity in the sense of the Moser--Trudinger inequality. In addition, the operator $\Delta_{H,t} \cdot=\text{div}_H(|D_H \cdot|_H^{t-2}D_H \cdot)$ for all $t>1$ is known as the $t$-sub-Laplacian, where $D_H \cdot$ is the subelliptic gradient.

Throughout this article, we assume that $g\colon \mathbb{R}\to\mathbb{R}$ satisfies the following hypotheses:
\begin{enumerate}[label=\textnormal{(\bf g\arabic*)},ref=\textnormal{\bf g\arabic*}]
	\item\label{g1} $g\in C(\mathbb{R},\mathbb{R})$ and there exists $\tau>2Q-\beta$ such that $g(s)=o(|s|^{\tau-1})$ as $s\to 0$.
\item\label{g2} \textit{Singular critical exponential growth at infinity}: There holds
    \begin{align*}
	\lim_{|s|\to +\infty} \frac{|g(s)|}{\operatorname{exp}\left(\alpha\big(1-\frac{\beta}{Q}\big) |s|^{Q^\prime}\right)}=
	\begin{cases}
		0       & \text{if  }\quad \alpha>\alpha_Q, \\
		+\infty & \text{if }\quad0<\alpha<\alpha_Q.
	\end{cases}
\end{align*}
\item\label{g3} \textit{Ambrosetti--Rabinowitz (A--R) condition}: There exists a constant $\theta>2Q-\beta$ such that
      $$0<\theta G(s):=\theta\int_{0}^{s} g(t)\,\mathrm{d}t\leq g(s)s\quad{\text{for all}}\quad s\in\mathbb{R}\setminus\{0\}. $$
\item\label{g4} There exist constants $p>2Q-\beta$ and $\mu>0$ such that
       $$\operatorname{sgn}(s)g(s)\geq\mu|s|^{p-1}\quad{\text{for all}}\quad s\in\mathbb{R}\setminus\{0\},$$
       where $\operatorname{sgn}\colon \mathbb{R}\setminus\{0\}\to \{-1,1\}$ denotes the signum function, given by   $\operatorname{sgn}(s)=-1$ for all $s<0$ and $\operatorname{sgn}(s)=1$ for all $s>0$.
\item\label{g5} The map $s\mapsto \displaystyle\frac{\widetilde{G}(s)}{|s|^{2Q-\beta-1}s}$ is strictly increasing on $(-\infty,0)$ and $(0,+\infty)$, where $\widetilde{G}\colon \mathbb{R}\to\mathbb{R}$ is given by
\begin{align*}
 \widetilde{G}(s)=g(s)s-(Q-\beta) G(s)\quad{\text{for all}}\quad s\in \mathbb{R}.
\end{align*}
\end{enumerate}
In virtue of \eqref{g1} and \eqref{g2}, for any $\varepsilon>0$, $\alpha>\alpha_Q$ and $q>Q$, there exists a constant $C_\varepsilon>0$ depending on $\varepsilon$, $q$, $\alpha$, $\beta$ and $Q$ such that
\begin{align}\label{eq1.5}
    |g(s)|\leq \varepsilon |s|^{\tau-1}+C_\varepsilon|s|^{q-1}\Phi\Big(\alpha\Big(1-\frac{\beta}{Q}\Big)|s|^{Q^\prime}\Big)\quad{\text{for all}}\quad s\in\mathbb{R},
\end{align}
and we also obtain from \eqref{g3} and \eqref{eq1.5} that
\begin{align}\label{eq1.6}
    \max\left\{|g(s)s|,\,|G(s)|\right\}\leq \varepsilon |s|^{\tau}+C_\varepsilon|s|^{q}\Phi\Big(\alpha\Big(1-\frac{\beta}{Q}\Big)|s|^{Q^\prime}\Big)\quad{\text{for all}}\quad s\in\mathbb{R}.
\end{align}

To illustrate the preceding hypotheses, we consider the following example of a nonlinear term satisfying all of them.
\begin{example} The function $G\colon \mathbb{R}\to \mathbb{R}$ given by $G(s):=\mu|s|^p \operatorname{exp}\big(\alpha_Q
\big(1-\frac{\beta}{Q}\big)|s|^{Q^\prime}\big)$ for all $s\in\mathbb{R}$, where $Q=2n+2$, $0<\beta<Q$, $\mu>0$, $p>2Q-\beta$, $g(s)=G^\prime(s)$ for all $s\in \mathbb{R}\setminus\{0\}$ and $g(0)=0$. Then $g$ satisfies all the hypotheses \eqref{g1}--\eqref{g5}.
\end{example}
The existence of normalized solutions to nonlinear Schr\"{o}dinger equations has received considerable attention in recent years. Such solutions arise naturally in the study of standing waves for the nonlinear Schr\"{o}dinger equation of the form
\begin{align*}
 i\partial_t\psi-\Delta\psi+f(\psi)=0
\quad\text{in }\quad\mathbb{R}^n\times (0,+\infty),
\end{align*}
through the standing-wave ansatz $\psi(x,t)=\operatorname{\exp}\left(-i\lambda t\right)u(x)$, which leads to the stationary equation as follows
\begin{align*}
 -\Delta u+\lambda u=f(u)\quad \text{in}\quad\mathbb{R}^n, \qquad\|u\|_2=a.
\end{align*}
The prescribed mass $a>0$ is fixed, while $\lambda$ serves as the associated Lagrange multiplier. The variational framework for normalized problems was developed by Jeanjean \cite{Jeanjean-1997} and Bartsch--Soave \cite{BS-2017}, and extended to the Sobolev-critical setting by Soave \cite{Soave-2020}. In the planar case, where the critical growth is governed by the Moser--Trudinger inequality, Alves--Ji--Miyagaki \cite{AJM-2022} investigated the existence of solutions for normalized, while Chen et al. \cite{CRTY-2023} subsequently studied normalized and ground state solutions for Schr\"{o}dinger equations in $\mathbb{R}^2$ with critical exponential nonlinearities. The theory has since been extended in several directions. Recently, Cassani et al. \cite{CHTZ-2026} established the existence of two positive normalized solutions to a mass-mixed elliptic problem in $\mathbb{R}^2$ for sufficiently small prescribed mass, one of local-minimum type and the other of mountain-pass type. For coupled Schr\"{o}dinger systems, Deng et al. \cite{DHZZ-2024} studied normalized solutions with critical exponential growth in the mass-supercritical regime. These works highlight the challenges arising from the prescribed mass constraint and the loss of compactness at the Moser--Trudinger critical threshold.

It is important to notice that \eqref{main problem} can be viewed as a constrained variational problem. In fact, its solutions are precisely the constrained critical points of the $C^1$-functional $\mathcal{J}\colon HW^{1,Q}(\mathbb{H}^n)\to\mathbb{R}$, associated with the corresponding Euler-Lagrange equation and defined by
\begin{align}\label{eq8.9990}
 \mathcal{J}(u):=\frac{1}{Q}\|D_Hu\|_Q^Q- \int_{\mathbb{H}^n}\frac{G(u)}{r(\xi)^\beta} \,\mathrm{d}\xi\quad\text{for all}\quad u\in HW^{1,Q}(\mathbb{H}^n),
\end{align}
constrained on the sphere
\begin{align*}
   S(a):=\left\{u\in HW^{1,Q}(\mathbb{H}^n)\colon\, \|u\|_Q=a\right \}.
\end{align*}
Moreover, every critical point $v_a\in {S}(a)$ of $\left.\mathcal{J}\right|_{{S}(a)}$ corresponds to a Lagrange multiplier $\lambda_a\in\mathbb{R}$ such that the pair $(v_a,\lambda_a)\in HW^{1,Q}(\mathbb{H}^n)\times\mathbb{R}$ solves \eqref{main problem}. Let $\mathcal{P}\colon HW^{1,Q}(\mathbb{H}^n)\to \mathbb{R}$ be the Poho\v{z}aev functional, which is given by
\begin{align}\label{eq9.999}
\mathcal{P}(u):=\|D_Hu\|_Q^Q+(Q-\beta) \int_{\mathbb{H}^n}\frac{G(u)}{r(\xi)^\beta}\,\mathrm{d}\xi-\int_{\mathbb{H}^n}\frac{g(u)u}{r(\xi)^\beta}\,\mathrm{d}\xi
\end{align}
for all $u\in HW^{1,Q}(\mathbb{H}^n)$. Now, we define the Poho\v{z}aev manifold by
\begin{align*}
 \mathcal{P}(a):=\left\{u\in S(a)\colon\,\mathcal{P}(u)=0\right\}.
\end{align*}
This naturally leads to the minimization problem for $\mathcal{J}$ over the Poho\v{z}aev manifold $\mathcal{P}(a)$, given by
\begin{align*}
    m(a):=\inf\left\{\mathcal{J}(u)\colon\, u\in \mathcal{P}(a)\right\}.
\end{align*}
In what follows, the critical points of $\left.\mathcal{J}\right|_{S(a)}$ that lie in $\mathcal{P}(a)$ and attain the minimum level $m(a)$ are referred to as ground state solutions of \eqref{main problem}. Formally, every critical point of $\left.\mathcal{J}\right|_{S(a)}$ lies in $\mathcal{P}(a)$. To prove this, one needs a Poho\v{z}aev identity for weak solutions of \eqref{main problem}, and we do not know such an identity for the $Q$-sub-Laplacian with the singular weight and critical exponential growth. Thus, a ground state solution has the least energy among all solutions of \eqref{main problem} lying in $\mathcal{P}(a)$.

The result below concerns the existence of ground state solutions to \eqref{main problem}.
\begin{theorem}\label{thm1.11}
Let $0<\beta<Q$ and $a>0$ be satisfied. In addition, suppose that $g$ satisfies \eqref{g1}--\eqref{g4}. Let $\mu^\ast=\mu^\ast(Q,\beta,a,p,\theta)>0$ be given by \eqref{eq8.mustar}. If \eqref{g4} holds with $\mu>\mu^\ast$, then \eqref{main problem} admits a weak solution $(u_a,\lambda_a)\in HW^{1,Q}(\mathbb{H}^n)\times\mathbb{R}$ with $\|u_a\|_Q=a$, $\lambda_a>0$, $u_a\in\mathcal{P}(a)$ and $\mathcal{J}(u_a)=\gamma_\mu(a,\beta)$, where $\gamma_\mu(a,\beta)$ is the mountain pass level defined in Section \ref{sec8}. Moreover, if \eqref{g5} holds, then $u_a$ is a nontrivial ground state solution of \eqref{main problem}, that is,  $$\mathcal{J}(u_a)=m(a)=\inf_{u\in \mathcal{P}(a)} \mathcal{J}(u).$$
\end{theorem}
We note that the constant $\mu^\ast=\mu^\ast(Q,\beta,a,p,\theta)>0$ defined in \eqref{eq8.mustar} satisfies the scaling relation $\mu^\ast(Q,\beta,a,p,\theta)=a^{-(Q-\beta)}\mu^\ast(Q,\beta,1,p,\theta)$. Thus, the condition $\mu>\mu^\ast$ in Theorem \ref{thm1.11} is equivalent to $a>a_\ast(\mu):=\left(\mu^\ast(Q,\beta,1,p,\theta)/\mu\right)^{\frac{1}{Q-\beta}}$. Hence, if $g$ satisfies \eqref{g1}--\eqref{g5}, then for every $\mu>0$, problem \eqref{main problem} admits a nontrivial ground state solution for every $a>a_\ast(\mu)$.

The rest of this article is organized as follows. Section \ref{sec3} contains the proof of Theorem \ref{thm1.1}, which is obtained by means of change-of-variables arguments. The proofs of Theorems \ref{thm1.3} and \ref{thm1.4} are given in Section \ref{sec4}, while Section \ref{sec5} is devoted to the proof of Theorem \ref{thm1.5}. Section \ref{sec6} establishes Theorems \ref{thm1.6} and \ref{thm1.7}, and Section \ref{sec7} contains the proofs of Theorems \ref{thm1.8} and \ref{thm1.9}. Finally, Section \ref{sec8} is devoted to the proof of Theorem \ref{thm1.11}, concerning the existence of nontrivial ground state solutions to the variational problem under consideration. Section \ref{sec9} explores several open questions arising from the present study and outlines potential directions for future research.
\section{The Adimurthi--Druet-Type Weighted Inequality and Its Consequence}\label{sec3}
This section is devoted to the proof of Theorem \ref{thm1.1}. The proof of Theorem \ref{thm1.1} uses a suitable change-of-variables technique to establish a weighted version of the Adimurthi--Druet-type Moser--Trudinger inequality on the Heisenberg group $\mathbb{H}^n$, for which the unweighted counterpart was established by Chen--Lu--Zhu \cite{Chen-Lu-Zhu-2021}. Moreover, we also derive the equivalence of this inequality to the critical weighted Moser--Trudinger inequality, as outlined in Corollary \ref{cor1.2}.
\subsection{{Proof of Theorem \ref{thm1.1}}}
Take $u\in HW^{1,Q}(\mathbb{H}^n)$ such that $\|u\|_\tau\leq1$. By setting $\theta=\tau^{\frac{1}{Q}}$, we introduce the function $v\colon \mathbb{H}^n\to\mathbb{R}$ through $u(\xi)=v(\delta_{\theta}(\xi))$ for all $\xi\in \mathbb{H}^n$. Then, straightforward calculations give $\|D_H u\|_Q= \|D_H v\|_Q$ and $\|u\|_Q=\theta^{-1}\|v\|_Q$. This shows that $\|v\|\leq 1$ and thus, we have
\begin{align*}
 \int_{\mathbb{H}^n}&\frac{\Phi\big(\alpha_Q \big(1 -\frac{\beta}{Q}\big)\big(1+\alpha \tau \|u\|_Q^Q \big)^{\frac{1}{Q-1}}|u|^{Q^\prime}\big)}{r(\xi)^\beta}\,\mathrm{d}\xi\\ &\qquad=\frac{1}{\tau^{1-\frac{\beta}{Q}}}\int_{\mathbb{H}^n}\frac{\Phi\big(\alpha_Q \big(1-\frac{\beta}{Q}\big)\big(1+\alpha \|v\|_Q^Q \big)^{\frac{1}{Q-1}}|v|^{Q^\prime}\big)}{r(\xi)^\beta}\,\mathrm{d}\xi.
\end{align*}
In particular, it implies
\begin{align}\label{eq3.1}
 \textnormal{ADMT}(\alpha,\beta,\tau)\leq \frac{1}{\tau^{1-\frac{\beta}{Q}}} \textnormal{ADMT}(\alpha,\beta,1).
\end{align}
Next, choose $v\in HW^{1,Q}(\mathbb{H}^n)$ with $\|v\|\leq1$. Define the function $w\colon \mathbb{H}^n\to\mathbb{R}$ by $w(\xi):=\big(1+\alpha \|v\|_Q^Q \big)^{\frac{1}{Q}}v(\xi)$ for all $\xi\in\mathbb{H}^n$. It follows that
\begin{align*}
    \|D_H w\|_Q^Q+(1-\alpha)\|w\|_Q^Q&=\big(1+\alpha \|v\|_Q^Q \big)\|D_H v\|_Q^Q+(1-\alpha)\big(1+\alpha \|v\|_Q^Q \big)\| v\|_Q^Q\\
    &=\big(1+\alpha \|v\|_Q^Q \big)(\|D_H v\|_Q^Q+\|v\|_Q^Q)-\alpha\big(1+\alpha \|v\|_Q^Q \big)\| v\|_Q^Q\\ &
    \leq 1+\alpha \|v\|_Q^Q -\alpha\big(1+\alpha \|v\|_Q^Q \big)\| v\|_Q^Q\leq 1
\end{align*}
and
\begin{align*}
   \int_{\mathbb{H}^n}\frac{\Phi\big(\alpha_Q \big(1-\frac{\beta}{Q}\big)\big(1+\alpha \|v\|_Q^Q \big)^{\frac{1}{Q-1}}|v|^{Q^\prime}\big)}{r(\xi)^\beta}\,\mathrm{d}\xi= \int_{\mathbb{H}^n}\frac{\Phi\big(\alpha_Q \big(1-\frac{\beta}{Q}\big)|w|^{Q^\prime}\big)}{r(\xi)^\beta}\,\mathrm{d}\xi.
\end{align*}
Hence, we conclude that $\|w\|_{1-\alpha}\leq 1$. Now, by applying the weighted  Moser--Trudinger inequality established by Lam--Lu \cite[Theorem 1.6]{Lam-Lu-2012}, one obtains from the above estimate that
\begin{align}\label{eq3.2}
 \textnormal{ADMT}(\alpha,\beta,1)\leq  \sup_{w\in HW^{1,Q}(\mathbb{H}^n),~\|w\|_{1-\alpha}\leq 1} \int_{\mathbb{H}^n}\frac{\Phi\big(\alpha_Q\big(1-\frac{\beta}{Q}\big) |w|^{Q^\prime}\big)}{r(\xi)^\beta}\,\mathrm{d}\xi<+\infty.
\end{align}
Together, \eqref{eq3.1} and \eqref{eq3.2} imply the desired inequality \eqref{eq1.2}. The proof is completed by showing that $\textnormal{ADMT}(1,\beta,\tau)=+\infty$. To establish this, we introduce
$$\xi=(z,t)\in\mathbb{H}^n\setminus\{O\},~z^\ast=\frac{z}{r(\xi)},~t^\ast=\frac{t}{r(\xi)^2}\quad\text{and}\quad\xi^\ast=(z^\ast,t^\ast). $$
By construction, the point $\xi^\ast$ lies in $\Sigma=\{\xi\in \mathbb{H}^n\colon\,r(\xi)=1\}$, the unit sphere on the Heisenberg group. In accordance with Cohn--Lu \cite{Cohn-Lu-2001}, we define the constant
\begin{align*}
c_m:=\int_\Sigma|z^\ast|^m\,\mathrm{d}\mu(\xi^\ast)=\omega_{2n-1}\Gamma\bigg(\frac{1}{2}\bigg)\Gamma\bigg(\frac{Q+m-2}{4}\bigg)\Gamma\bigg(\frac{Q+m}{4}\bigg)^{-1}
\end{align*}
for all $m>-2n$, where $\mathrm{d}\mu(\xi^\ast)$ denotes the Radon measure on $\Sigma$, and $\omega_{2n-1}=2\pi^n\Gamma(n)^{-1}$ denotes the surface measure of the unit sphere $\mathbb{S}^{2n-1}$ in $\mathbb{R}^{2n}$. In terms of $c_Q$, the critical constant $\alpha_Q$ takes the form
\begin{align*}
 \alpha_Q:=Qc_Q^{\frac{1}{Q-1}}=Q\bigg(2\pi^n\Gamma\bigg(\frac{1}{2}\bigg)\Gamma\bigg(\frac{Q-1}{2}\bigg){\Gamma\bigg(\frac{Q}{2}\bigg)}^{-1}{\Gamma(n)}^{-1}\bigg)^{\frac{1}{Q-1}}.
\end{align*}
Now, let $\{R_k\}_k$ be the sequence of real numbers given by $R_k:=\frac{(\ln k)^{\frac{1}{2Q}}}{\ln(\ln k)}$ for all integers $k\geq3$. Note that $R_k\to+\infty$ as $k\to\infty$. Define the Moser's function sequence $\{v_k\}_k$ given by
\begin{align}\label{eq3.3}
  v_k(\xi):=\frac{1}{c^{\frac{1}{Q}}_Q}\begin{cases}
    (\ln k)^{\frac{Q-1}{Q}}\quad\quad\quad\quad\quad\text{if}\quad 0\leq r(\xi)\leq \frac{R_k}{k},\\
    (\ln k)^{-\frac{1}{Q}}\ln\big(\frac{R_k}{r(\xi)}\big)\quad\quad\text{if}\quad \frac{R_k}{k}\leq r(\xi)\leq R_k,\\
    0 \quad\quad\quad\quad\quad\quad\quad\quad\quad\text{if}\quad r(\xi)\geq R_k,
\end{cases}
\end{align}
for all integers $k\geq3$. By applying the result of Cohn--Lu \cite[Lemma 2.1]{Cohn-Lu-2001}, we obtain $|D_H r(\xi)|_H=\frac{|z|}{r(\xi)}$, which consequently yields
$$ |D_H v_k|_H=c^{-\frac{1}{Q}}_Q(\ln k)^{-\frac{1}{Q}}\frac{|z|}{r(\xi)^2}\chi _{B_{R_k}\setminus B_{k^{-1}R_k}}(\xi)\quad\text{for all}\quad k\in\mathbb{N}.$$
Invoking the polar coordinate formula of Folland--Stein \cite[Proposition 1.15, p. 14]{FS-82}, we obtain
   \begin{align*}
    \|D_H v_k\|^Q_Q&=\int_\Sigma\int_{k^{-1}R_k}^{R_k}\bigg(c^{-\frac{1}{Q}}_Q(\ln k)^{-\frac{1}{Q}}\frac{|z^\ast|}{\eta}\bigg)^Q\eta^{Q-1}\,\mathrm{d}\eta\,\mathrm{d}\mu(\xi^\ast)\\
    &=c^{-1}_Q(\ln k)^{-1}c_Q\int_{k^{-1}R_k}^{R_k}\eta^{-1}\,\mathrm{d}\eta=c^{-1}_Q(\ln k)^{-1}c_Q\ln\bigg(\frac{R_k}{k^{-1}R_k}\bigg)=1
\end{align*}
and
 \begin{align*}
    \| v_k\|^Q_Q&=\int_\Sigma\int_{0}^{k^{-1}R_k}\bigg(c^{-\frac{1}{Q}}_Q(\ln k)^{\frac{Q-1}{Q}}\bigg)^Q\eta^{Q-1}\,\mathrm{d}\eta\,\mathrm{d}\mu(\xi^\ast)\\
    &\qquad+\int_\Sigma\int_{k^{-1}R_k}^{R_k}\bigg(c^{-\frac{1}{Q}}_Q(\ln k)^{-\frac{1}{Q}}\bigg|\ln{\bigg(\frac{R_k}{\eta}}\bigg)\bigg|\bigg)^Q\eta^{Q-1}\,\mathrm{d}\eta\,\mathrm{d}\mu(\xi^\ast)\\
    &=c^{-1}_Q\bigg(\int_\Sigma\,\mathrm{d}\mu(\xi^\ast)\bigg)(\ln k)^{Q-1}\int_{0}^{k^{-1}R_k}\eta^{Q-1}\,\mathrm{d}\eta\\
&\qquad+\frac{c^{-1}_Q\big(\int_\Sigma\,\mathrm{d}\mu(\xi^\ast)\big)}{\ln k}\int_{k^{-1}R_k}^{R_k}\bigg|\ln{\bigg(\frac{R_k}{\eta}}\bigg)\bigg|^Q\eta^{Q-1}\,\mathrm{d}\eta\\
&=c^{-1}_Q\bigg(\int_\Sigma\,\mathrm{d}\mu(\xi^\ast)\bigg)\bigg(\frac{(\ln k)^{Q-1}}{Q}(k^{-1}R_k)^Q+\frac{R_k^Q}{\ln k}\int_{0}^{\ln k}s^Q \operatorname{exp}(-Qs)\,\mathrm{d}s\bigg)\\
&=\frac{\widetilde{\Pi}_{k} R_k^Q}{\ln k}\bigg(1+\frac{1}{Q\Pi_{k}}\bigg(\frac{\ln k}{k}\bigg)^Q\bigg)=\frac{\widetilde{\Pi}_{k} R_k^Q}{\ln k}\big(1+o_k(1)\big)\quad\text{as}\quad k\to\infty,
\end{align*}
where the estimates above are obtained using
\begin{align*}
\widetilde{\Pi}_{k}=c^{-1}_Q\bigg(\int_\Sigma\,\mathrm{d}\mu(\xi^\ast)\bigg)\Pi_{k}\qquad\text{and}\qquad  \Pi_{k}=\int_{0}^{\ln k}s^Q \operatorname{exp}(-Qs)\,\mathrm{d}s
\end{align*}
as well as the limit
\begin{align*}
    \lim_{k\to\infty}\Pi_{k}=\lim_{k\to\infty}\int_{0}^{\infty}s^Q \operatorname{exp}(-Qs)\chi_{[0,\ln k]}(s)\mathrm{d}s=\int_{0}^{\infty}s^Q \operatorname{exp}(-Qs)\mathrm{d}s=\frac{\Gamma(Q+1)}{Q^{Q+1}},
\end{align*}
where the monotone convergence theorem is applied; see, e.g., Royden--Fitzpatrick \cite[{\S}4.3, p.~83]{Royden-Fitzpatrick-1988}. It follows directly that
\begin{align*}
    \|v_k\|_\tau^Q=1+\tau\bigg(\frac{\widetilde{\Pi}_{k} R_k^Q}{\ln k}\bigg)\big(1+o_k(1)\big)=1+o_k(1)\quad\text{as}\quad k\to\infty.
\end{align*}
For each $k\in\mathbb{N}$, let $w_k$ be defined by $ w_k(\xi)=\frac{v_k(\xi)}{\|v_k\|_\tau}$ for all $\xi\in\mathbb{H}^n$. This shows that the sequence $\{w_k\}_k$ satisfies $\|w_k\|_\tau=1$. Now, taking $\alpha=1$, we obtain the following estimate on the Kor\'anyi open ball $B_{k^{-1}R_k}$, given by
\begin{align*}
   \alpha_Q \Big(1-\frac{\beta}{Q}\Big)&\big(1+\tau\|w_k\|^Q_Q\big)^{\frac{1}{Q-1}} |w_k|^{Q^\prime}=Qc_Q^{\frac{1}{Q-1}}\Big(1-\frac{\beta}{Q}\Big)|v_k|^{\frac{Q}{Q-1}}\frac{(1+2\tau\|v_k\|_Q^Q)^{\frac{1}{Q-1}}}{(1+\tau\|v_k\|_Q^Q)^{\frac{2}{Q-1}}}\\
   &=(Q-\beta)\ln k\bigg(1-\frac{\tau^2}{Q-1}\|v_k\|_Q^{2Q}+\frac{2\tau^3}{Q-1}\|v_k\|_Q^{3Q}(1+o_k(1))\bigg)\\
   &=(Q-\beta)\ln k-(Q-\beta)\ln k\bigg(\frac{\tau^2}{Q-1}\frac{\widetilde{\Pi}_{k}^2 R_k^{2Q}}{(\ln k)^2}(1+o_k(1))\bigg)\\
   &=(Q-\beta)\ln k-\tau^2\bigg(\frac{Q-\beta}{Q-1}\bigg)\frac{ \widetilde{\Pi}_{k}^2}{(\ln(\ln k))^{2Q}}(1+o_k(1))\quad\text{as}\quad k\to\infty.
\end{align*}
Also, the weighted measure of the Kor\'anyi open ball $B_{k^{-1}R_k}$ is computed as
\begin{align*}
    \int_{B_{\frac{R_k}{k}}}\frac{\mathrm{d}\xi}{r(\xi)^\beta}&=\int_\Sigma\int_{0}^{\frac{R_k}{k}}\eta^{Q-\beta-1}\,\mathrm{d}\eta\,\mathrm{d}\mu(\xi^\ast)=\bigg(\int_\Sigma\,\mathrm{d}\mu(\xi^\ast)\bigg)\frac{\operatorname{\exp}\big((Q-\beta)(\ln R_k-\ln k)\big)}{Q-\beta}\\
    &=\bigg(\int_\Sigma\,\mathrm{d}\mu(\xi^\ast)\bigg)\frac{\operatorname{\exp}\Big((Q-\beta)\big(\frac{1}{2Q}\ln(\ln k)-\ln(\ln(\ln k))-\ln k\big)\Big)}{Q-\beta}.
\end{align*}
Taken together, these estimates give
\begin{align*}
    \textnormal{ADMT}&(1,\beta,\tau)\geq\int_{\mathbb{H}^n}\frac{\Phi\big(\alpha_Q \big(1-\frac{\beta}{Q}\big)\big(1+ \tau \|w_k\|_Q^Q \big)^{\frac{1}{Q-1}}|w_k|^{Q^\prime}\big)}{r(\xi)^\beta}\,\mathrm{d}\xi\\
    &\geq \int_{B_{\frac{R_k}{k}}}\frac{\Phi\Big((Q-\beta)\ln k-\tau^2\big(\frac{Q-\beta}{Q-1}\big)\frac{ \widetilde{\Pi}_{k}^2}{(\ln(\ln k))^{2Q}}(1+o_k(1))\Big)}{r(\xi)^\beta}\,\mathrm{d}\xi\\
    &\gtrsim \operatorname{\exp}\bigg((Q-\beta)\ln k-\tau^2\bigg(\frac{Q-\beta}{Q-1}\bigg)\frac{ \widetilde{\Pi}_{k}^2}{(\ln(\ln k))^{2Q}}(1+o_k(1))\bigg)\int_{B_{\frac{R_k}{k}}}\frac{\mathrm{d}\xi}{r(\xi)^\beta}\\
    &\gtrsim \Theta_k\to+\infty\quad\text{as}\quad k\to\infty,
\end{align*}
where the sequence $\{\Theta_k\}_k$ is given by  $$\Theta_k=\operatorname{\exp}\Bigg(\bigg(\frac{Q-\beta}{2Q}\bigg)\ln(\ln k)\bigg(1-2Q\bigg(\frac{\ln(\ln(\ln k))}{\ln(\ln k)}\bigg)-\frac{2Q\tau^2}{Q-1}\bigg(\frac{ \widetilde{\Pi}_{k}^2}{(\ln(\ln k))^{2Q+1}}\bigg)(1+o_k(1))\bigg)\Bigg)$$ for each integer $k\geq3$.
This proves the desired assertion, namely, $\textnormal{ADMT}(1,\beta,\tau)=+\infty$, which completes the proof. \qed


Now, by taking $\alpha=0$ in \eqref{eq1.2}, and by the proof of Theorem \ref{thm1.1}, we obtain the following result, where the critical weighted Moser--Trudinger inequality is the one of Lam--Lu \cite[Theorem 1.6]{Lam-Lu-2012}.
\begin{corollary}\label{cor1.2}
The Adimurthi--Druet-type weighted Moser--Trudinger inequality is equivalent to the critical weighted Moser--Trudinger inequality on the Heisenberg group $\mathbb{H}^n$.
\end{corollary}
To end this section, we recall the following particular case of the weighted Moser--Trudinger inequality on the Heisenberg group; see, e.g., Cohn--Lam--Lu--Yang \cite[Theorem 1.1]{Cohn-Lam-Lu-Yang-2012}.
\begin{knownlem}\label{lem3.1}
    For any $\alpha>0$, $0\leq\beta<Q$ and $u\in HW^{1,Q}(\mathbb{H}^n)$, there holds
    \begin{align*}
        \frac{\Phi\big(\alpha|u|^{Q^\prime}\big)}{r(\xi)^\beta}\in L^1(\mathbb{H}^n).
    \end{align*}
\end{knownlem}
\section{The Subcritical Weighted Moser--Trudinger Inequality and Its Maximizer}\label{sec4}
The discussion in this section is twofold. Firstly, we derive new two-sided asymptotic bounds for the supremum associated with the subcritical weighted Moser--Trudinger inequality near the critical threshold in Theorem \ref{thm1.3}, building on the rearrangement-free approach of the level-set method of Lam--Lu--Zhang \cite{Lam-Lu-Zhang-2017}, while the underlying inequality was already established by Lam--Lu--Tang \cite{Lam-Lu-Tang-2014}. Secondly, we establish the existence of maximizers for this inequality in Theorem \ref{thm1.4}, following Ishiwata et al. \cite{Ishiwata-Nakamura-Wadade-2014}, with compactness obtained through a weighted Sobolev embedding in the full borderline Sobolev space on the Heisenberg group, without radial restrictions. The result below plays a central role in our argument.
\begin{lemma}\label{lem4.1}
    Let $0<\alpha<\alpha_Q$ and $0\leq\beta<Q$ be real numbers. Then there holds
$$\textnormal{SCMT}(\alpha,\beta)=\sup_{u\in HW^{1,Q}(\mathbb{H}^n),\,\|D_H u\|_{Q}\leq 1,\,\|u\|_Q=1} \int_{\mathbb{H}^n}\frac{\Phi\big(\alpha \big(1-\frac{\beta}{Q}\big)|u|^{Q^\prime}\big)}{r(\xi)^\beta}\,\mathrm{d}\xi.$$
\end{lemma}
\begin{proof}
 Let $u\in HW^{1,Q}(\mathbb{H}^n)$ be such that $\|D_H u\|_{Q}\leq 1$. Define $w\colon \mathbb{H}^n\to \mathbb{R}$ by $w(\xi)=u(\delta_\theta(\xi))$ for all $\xi\in \mathbb{H}^n$, where $\theta=\|u\|_Q>0$. By direct computations, we obtain 
 \begin{align*}
    \|w\|_Q=\theta^{-1}\| u\|_Q =1,\quad\|D_H w\|_Q=\|D_H u\|_Q\leq 1
 \end{align*}
and
 \begin{align*}
     \int_{\mathbb{H}^n}\frac{\Phi\big(\alpha \big(1-\frac{\beta}{Q}\big)|w|^{Q^\prime}\big)}{r(\xi)^\beta}\,\mathrm{d}\xi=\|u\|_Q^{-(Q-\beta)}\int_{\mathbb{H}^n}\frac{\Phi\big(\alpha \big(1-\frac{\beta}{Q}\big)|u|^{Q^\prime}\big)}{r(\xi)^\beta}\,\mathrm{d}\xi.
 \end{align*}
Hence, the preceding relation leads to the desired conclusion. This finishes the proof.
\end{proof}
In light of Lemma \ref{lem4.1}, we may assume, without loss of generality, that $\|u\|_Q=1$ when proving the subcritical weighted Moser--Trudinger inequality. Next, by having established the preceding result at hand, we now prove Theorem \ref{thm1.3}.
\subsection{{Proof of Theorem \ref{thm1.3}}} In view of Lemma \ref{lem4.1} and the fact that $C_0^\infty(\mathbb{H}^n)$ is dense in $HW^{1,Q}(\mathbb{H}^n)$, it suffices to prove Theorem \ref{thm1.3} for all $u\in C_0^\infty(\mathbb{H}^n)\setminus\{0\}$, $u\geq 0$, $\|u\|_Q= 1$ and $\|D_Hu\|_Q\leq 1$. Now, we define the set
$$ S:=\Bigg\{\xi\in\mathbb{H}^n\colon\,u(\xi)>\bigg(1-\Big(\frac{\alpha}{\alpha_Q}\Big)^{Q-1}\bigg)^{\frac{1}{Q}}\Bigg\}. $$
Note that $S$ is a bounded set in $\mathbb{H}^n$, and we have
\begin{align}\label{eq4.1}
 |S|= \int_{S} \,\mathrm{d}\xi\leq \bigg(1-\Big(\frac{\alpha}{\alpha_Q}\Big)^{Q-1}\bigg)^{-1}\int_{S}|u|^Q\,\mathrm{d}\xi\leq \bigg(1-\Big(\frac{\alpha}{\alpha_Q}\Big)^{Q-1}\bigg)^{-1}.
\end{align}
By decomposing the integral in the following manner, we arrive at
\begin{align*}
 \int_{\mathbb{H}^n}\frac{\Phi\big(\alpha \big(1-\frac{\beta}{Q}\big)|u|^{Q^\prime}\big)}{r(\xi)^\beta}\,\mathrm{d}\xi=\bigg(\int_{\mathbb{H}^n\setminus S}+\int_{S}\bigg)\frac{\Phi\big(\alpha \big(1-\frac{\beta}{Q}\big)|u|^{Q^\prime}\big)}{r(\xi)^\beta}\,\mathrm{d}\xi=:I+J.
\end{align*}
We first estimate $I$. By the inclusion $\mathbb{H}^n\setminus S\subset \{\xi\in \mathbb{H}^n\colon\,u(\xi)\leq 1\}$ and the choice $\alpha\in(0,\alpha_Q)$, we infer that
\begin{equation}\label{eq4.2}
    \begin{aligned}
      I & \leq \int_{\{\xi\in \mathbb{H}^n\colon\,u(\xi)\leq 1\}}\frac{1}{r(\xi)^\beta}\sum_{j=Q-1}^{\infty} \frac{\big(\alpha\big(1-\frac{\beta}{Q}\big)\big)^j}{j!}\,|u|^{Q^\prime j}\,\mathrm{d}\xi  \\
&\leq \int_{\{\xi\in \mathbb{H}^n\colon\,u(\xi)\leq 1\}}\frac{1}{r(\xi)^\beta}\sum_{j=Q-1}^{\infty} \frac{\big(\alpha\big(1-\frac{\beta}{Q}\big)\big)^j}{j!}\,|u|^{Q}\,\mathrm{d}\xi  \\
 & \leq \operatorname{\exp}\bigg(\alpha\bigg(1-\frac{\beta}{Q}\bigg)\bigg) \Bigg(\int_{\{\xi\in \mathbb{H}^n\colon\, u(\xi)\leq 1,\,r(\xi)<1\}}+\int_{\{\xi\in \mathbb{H}^n\colon\,u(\xi)\leq1,\,r(\xi)\geq1\}}\Bigg)\frac{|u|^{Q}}{r(\xi)^\beta}\,\mathrm{d}\xi  \\
& \leq \operatorname{\exp}\bigg(\alpha_Q\bigg(1-\frac{\beta}{Q}\bigg)\bigg)\Bigg(\int_{\{ \xi\in \mathbb{H}^n\colon\,r(\xi)<1\}}\frac{1}{r(\xi)^\beta}\,\mathrm{d}\xi+\int_{\{ \xi\in \mathbb{H}^n\colon\,r(\xi)\geq1\}}|u|^{Q}\,\mathrm{d}\xi\Bigg) \\
&\leq  C(Q,\beta),   
    \end{aligned}
\end{equation}
where $C(Q,\beta)>0$ is a constant depending only on $Q$ and $\beta$. Finally, we compute $J$. Define the function $w\colon S\to\mathbb{R}$ by $$w(\xi):=u(\xi)-\bigg(1-\Big(\frac{\alpha}{\alpha_Q}\Big)^{Q-1}\bigg)^{\frac{1}{Q}}\quad\text{for all}\quad\xi\in S.$$
It is evident that $w\in HW^{1,Q}_0(S)$ with $\|D_H w\|_{L^Q(S)}\leq 1$. Recall the following elementary inequality: for every $\varepsilon>0$ and every $p>1$, one has
 \begin{align}\label{eq4.999}
     (a+b)^p\leq (1+\varepsilon)a^p+C_\varepsilon b^p\quad\text{for all}\quad a,b>0,
\end{align}
where $C_\varepsilon=(1-(1+\varepsilon)^{-\frac{1}{p-1}})^{1-p}$. In particular, by setting $\varepsilon=\frac{\alpha_Q}{\alpha}-1$, the above inequality directly yields for all $\xi\in S$ that
\begin{equation}\label{eq4.3}
\begin{aligned}
    |u(\xi)|^{Q^\prime}&\leq \Bigg(|w(\xi)|+\bigg(1-\Big(\frac{\alpha}{\alpha_Q}\Big)^{Q-1}\bigg)^{\frac{1}{Q}}\Bigg)^{Q^\prime}\\
    &\leq (1+\varepsilon)|w(\xi)|^{Q^\prime}+C_\varepsilon \bigg(1-\Big(\frac{\alpha}{\alpha_Q}\Big)^{Q-1}\bigg)^{\frac{1}{Q-1}}  \\
    &=  \frac{\alpha_Q}{\alpha}|w(\xi)|^{Q^\prime}+\bigg(1-\Big(\frac{\alpha}{\alpha_Q}\Big)^{Q-1}\bigg)^{-\frac{1}{Q-1}}\bigg(1-\Big(\frac{\alpha}{\alpha_Q}\Big)^{Q-1}\bigg)^{\frac{1}{Q-1}}  \\
    &=\frac{\alpha_Q}{\alpha}|w(\xi)|^{Q^\prime}+1. 
\end{aligned}
\end{equation}
Now, applying the weighted Moser--Trudinger inequality on bounded subsets of the Heisenberg group, see, e.g., Lam--Lu--Tang \cite[Lemma 4.1]{Lam-Lu-Tang-2012}, together with \eqref{eq4.1}, \eqref{eq4.3}, and $\alpha\in(0,\alpha_Q)$, gives
\begin{equation}\label{eq4.4}
\begin{aligned}
 J & \leq \int_{S}\frac{\operatorname{\exp}\big(\alpha\big(1-\frac{\beta}{Q}\big) |u|^{Q^\prime}\big)}{r(\xi)^\beta}\,\mathrm{d}\xi\leq \operatorname{\exp}\Big(\alpha\Big(1-\frac{\beta}{Q}\Big)\Big)\int_{S}\frac{\operatorname{\exp}\big(\alpha_Q\big(1-\frac{\beta}{Q}\big) |w|^{Q^\prime}\big)}{r(\xi)^\beta}\,\mathrm{d}\xi  \\
 &\leq \operatorname{\exp}\Big(\alpha_Q\Big(1-\frac{\beta}{Q}\Big)\Big)C(Q,\beta)|S|^{1-\frac{\beta}{Q}} \leq C(Q,\beta)\bigg(1-\Big(\frac{\alpha}{\alpha_Q}\Big)^{Q-1}\bigg)^{-(1-\frac{\beta}{Q})}, 
\end{aligned}
\end{equation}
where $C(Q,\beta)>0$ is a constant, which may vary from line to line and depends only on $Q$ and $\beta$. It then follows from \eqref{eq4.2} and \eqref{eq4.4} that
\begin{align}\label{eq4.5}
  \textnormal{SCMT}(\alpha,\beta)\leq \bigg(1-\Big(\frac{\alpha}{\alpha_Q}\Big)^{Q-1}\bigg)^{-(1-\frac{\beta}{Q})}C(Q,\beta).
\end{align}
In what follows, we proceed to establish the sharpness of  $\alpha_Q$, that is, to verify that $\text{SCMT}(\alpha_Q,\beta)=+\infty$. To do this, we fix $s>0$ and define the Moser’s function sequence $\{w_k\}_k$ by
\begin{align}\label{eq4.6}
  w_k(\xi):=\frac{1}{c^{\frac{1}{Q}}_Q}\begin{cases}
    (\ln k)^{\frac{Q-1}{Q}}\quad\quad\quad\quad\quad\text{if}\quad r(\xi)\in[0,sk^{-1}],\\
    (\ln k)^{-\frac{1}{Q}}\ln\big(\frac{s}{r(\xi)}\big)\quad\quad\text{if}\quad r(\xi)\in[sk^{-1},s],\\
    0 \quad\quad\quad\quad\quad\quad\quad\quad\quad\text{if}\quad r(\xi)\in[s,+\infty).
\end{cases}
\end{align}
for all integers $k\geq2$. By arguing similarly as in Theorem \ref{thm1.1}, we can prove that $\|D_H w_k\|_Q=1$. Notice that for all $p\geq Q$, by using integration by parts, we also obtain
\begin{equation}\label{eq4.7}
\begin{aligned}
    \| w_k\|^p_p&=s^Qc^{-\frac{p}{Q}}_Q\bigg(\int_\Sigma\,\mathrm{d}\mu(\xi^\ast)\bigg)\bigg(\frac{(\ln k)^{(\frac{Q-1}{Q})p}}{Q}(k^{-1})^Q +\frac{1}{(\ln k)^{\frac{p}{Q}}}\int_{0}^{\ln k}t^p \operatorname{exp}(-Qt)\,\mathrm{d}t\bigg)  \\
    &=\frac{ps^Qc^{-\frac{p}{Q}}_Q}{Q}\bigg(\int_\Sigma\,\mathrm{d}\mu(\xi^\ast)\bigg)\frac{1}{(\ln k)^{\frac{p}{Q}}}\int_{0}^{\ln k}t^{p-1} \operatorname{exp}(-Qt)\,\mathrm{d}t \\
    &\leq \frac{ps^Qc^{-\frac{p}{Q}}_Q}{Q}\bigg(\int_\Sigma\,\mathrm{d}\mu(\xi^\ast)\bigg)\frac{1}{(\ln k)^{\frac{p}{Q}}}\int_{0}^{\infty}t^{p-1} \operatorname{exp}(-Qt)\,\mathrm{d}t  \\
    &=\frac{ps^Qc^{-\frac{p}{Q}}_Q \Gamma(p)}{Q^{p+1}}\bigg(\int_\Sigma\,\mathrm{d}\mu(\xi^\ast)\bigg)\frac{1}{(\ln k)^{\frac{p}{Q}}}. 
\end{aligned}
\end{equation}
As a result, one deduces that
\begin{align}\label{eq4.8}
 \|D_H w_k\|_Q=1\qquad \text{and}\qquad 0\leq \|w_k\|_p^p\leq O\bigg(\frac{1}{(\ln k)^{\frac{p}{Q}}}\bigg) \quad\text{for all}\quad p\geq Q.
\end{align}
The estimates obtained above will play an important role in the subsequent sections. In addition, one can deduce that
\begin{align*}
  \int_{\mathbb{H}^n}&\frac{\Phi\big(\alpha_Q \big(1-\frac{\beta}{Q}\big)|w_k|^{Q^\prime}\big)}{r(\xi)^\beta}\,\mathrm{d}\xi\geq \int_{\big\{\xi\in \mathbb{H}^n\colon~ r(\xi)<sk^{-1}\}}\frac{\Phi\big(\alpha_Q \big(1-\frac{\beta}{Q}\big)|w_k|^{Q^\prime}\big)}{r(\xi)^\beta}\,\mathrm{d}\xi\\
&=\int_\Sigma\int_{0}^{sk^{-1}}\Phi\bigg(\alpha_Q \bigg(1-\frac{\beta}{Q}\bigg)\bigg|c^{-\frac{1}{Q}}_Q(\ln k)^{\frac{Q-1}{Q}}\bigg|^{Q^\prime}\bigg)\eta^{Q-\beta-1}\,\mathrm{d}\eta\,\mathrm{d}\mu(\xi^\ast)\\
  &= \bigg(\int_\Sigma \,\mathrm{d}\mu(\xi^\ast)\bigg)\Phi\big((Q-\beta)\ln k\big)\int_{0}^{sk^{-1}}\eta^{Q-\beta-1}\,\mathrm{d}\eta\\
  &\gtrsim \bigg(\int_\Sigma \,\mathrm{d}\mu(\xi^\ast)\bigg)\,\ \frac{\operatorname{\exp}\big((Q-\beta)\ln k\big)(sk^{-1})^{Q-\beta}}{Q-\beta}\gtrsim 1.
\end{align*}
This together with \eqref{eq4.8} implies that
$$ \|w_k\|_Q^{-(Q-\beta)} \int_{\mathbb{H}^n}\frac{\Phi\big(\alpha_Q \big(1-\frac{\beta}{Q}\big)|w_k|^{Q^\prime}\big)}{r(\xi)^\beta}\,\mathrm{d}\xi\gtrsim (\ln k)^{\frac{Q-\beta}{Q}}\to+\infty\quad\text{as}\quad k\to\infty. $$
Thus, we have shown that $\text{SCMT}(\alpha_Q,\beta)=+\infty$, and therefore the sharpness of $\alpha_Q$ is established. By virtue of \eqref{eq4.7}, we may further write
\begin{align*}
 \| w_k\|^Q_Q&=s^Qc^{-1}_Q\bigg(\int_\Sigma\,\mathrm{d}\mu(\xi^\ast)\bigg)\bigg(\frac{(\ln k)^{Q-1}}{Q}(k^{-1})^Q+\frac{1}{\ln k}\int_{0}^{\ln k}t^Q \operatorname{exp}(-Qt)\,\mathrm{d}t\bigg)\\
    &\leq   s^Qc^{-1}_Q\bigg(\int_\Sigma\,\mathrm{d}\mu(\xi^\ast)\bigg)\bigg(\frac{(\ln k)^{Q-1}}{Q}(k^{-1})^Q+\frac{1}{\ln k}\int_{0}^{\infty}t^Q \operatorname{exp}(-Qt)\,\mathrm{d}t\bigg) \\
    &=s^Qc^{-1}_Q\bigg(\int_\Sigma\,\mathrm{d}\mu(\xi^\ast)\bigg)\bigg(\frac{(\ln k)^{Q-1}}{Q}(k^{-1})^Q+\frac{1}{\ln k}\frac{\Gamma(Q+1)}{Q^{Q+1}}\bigg).
\end{align*}
The estimates in \eqref{eq4.8} and above imply that there exist constants $C_1,C_2>0$, independent of $k$, such that
\begin{align*}
 \|D_H w_k\|_Q=1\qquad \text{and}\qquad 0\leq \|w_k\|_Q^Q\leq \frac{1}{\ln k}\bigg(C_1+C_2\bigg(\frac{\ln k}{k}\bigg)^Q\bigg).
\end{align*}
The information obtained above consequently yields
\begin{align*}
   \text{SCMT}(\alpha,\beta)&\geq \cfrac{1}{\| w_k\|_Q^{Q-\beta}}\int_{\{\xi\in\mathbb{H}^n\colon~r(\xi)<sk^{-1}\}}\frac{\Phi\big(\alpha \big(1-\frac{\beta}{Q}\big)|w_k|^{Q^\prime}\big)}{r(\xi)^\beta}\,\mathrm{d}\xi\\
   &\geq \cfrac{(\ln k)^{\frac{Q-\beta}{Q}}}{\Big(C_1+C_2\big(\frac{\ln k}{k}\big)^Q\Big)^{\frac{Q-\beta}{Q}}}\int_\Sigma\int_{0}^{sk^{-1}}\Phi\bigg(\alpha \bigg(1-\frac{\beta}{Q}\bigg)\bigg|c^{-\frac{1}{Q}}_Q(\ln k)^{\frac{Q-1}{Q}}\bigg|^{Q^\prime}\bigg)\\
   &\qquad\times\eta^{Q-\beta-1}\,\mathrm{d}\eta\,\mathrm{d}\mu(\xi^\ast)\\
   &\gtrsim \cfrac{(\ln k)^{\frac{Q-\beta}{Q}} \bigg(\displaystyle\int_\Sigma\,\mathrm{d}\mu(\xi^\ast)\bigg)(sk^{-1})^{Q-\beta}}{\Big(C_1+C_2\big(\frac{\ln k}{k}\big)^Q\Big)^{\frac{Q-\beta}{Q}}(Q-\beta)}~\operatorname{\exp}\Bigg(\alpha\Big(1-\frac{\beta}{Q}\Big)\bigg|c^{-\frac{1}{Q}}_Q(\ln k)^{\frac{Q-1}{Q}}\bigg|^{Q^\prime}\Bigg)\\
   &\gtrsim \cfrac{(\ln k)^{\frac{Q-\beta}{Q}}}{\Big(C_1+C_2\big(\frac{\ln k}{k}\big)^Q\Big)^{\frac{Q-\beta}{Q}}}~\operatorname{\exp}\Big((Q-\beta)\ln k\Big(\frac{\alpha}{\alpha_Q}-1\Big)\Big).
\end{align*}
For $k$ sufficiently large, it follows that $\frac{1}{(C_1+C_2(\frac{\ln k}{k})^Q)^{\frac{Q-\beta}{Q}}}\approx \frac{1}{C_1^{\frac{Q-\beta}{Q}}}$. The above inequality thereby implies that
\begin{align*}
   \text{SCMT}(\alpha,\beta) \gtrsim (\ln k)^{\frac{Q-\beta}{Q}}\operatorname{\exp}\Big((Q-\beta)\ln k\Big(\frac{\alpha}{\alpha_Q}-1\Big)\Big).
\end{align*}
Taking into account that $\alpha\underset{\approx}{<}\alpha_Q$, we may fix $k$ sufficiently large such that
$$ \ln k\approx\frac{1}{1-\frac{\alpha}{\alpha_Q}}\quad\text{and}\quad \operatorname{\exp}\Big((Q-\beta)\ln k\Big(\frac{\alpha}{\alpha_Q}-1\Big)\Big)\approx\operatorname{\exp}\big(-(Q-\beta)\big).$$
Thus, whenever $\alpha\underset{\approx}{<}\alpha_Q$, we obtain
$$\textnormal{SCMT}(\alpha,\beta) \gtrsim\frac{\operatorname{\exp}\big(-(Q-\beta)\big)}{\big(1-\frac{\alpha}{\alpha_Q}\big)^{\frac{Q-\beta}{Q}}}\approx\frac{(Q-1)^{\frac{Q-\beta}{Q}}~\operatorname{\exp}\big(-(Q-\beta)\big)}{\big(1-\big(\frac{\alpha}{\alpha_Q}\big)^{Q-1}\big)^{(1-\frac{\beta}{Q})}}, $$
where we have used $\frac{1-(\frac{\alpha}{\alpha_Q})^{Q-1}}{1-\frac{\alpha}{\alpha_Q}}\approx Q-1$ for $\alpha\underset{\approx}{<}\alpha_Q$. Hence, we deduce that
\begin{align}\label{eq4.9}
    \textnormal{SCMT}(\alpha,\beta)\geq \bigg(1-\Big(\frac{\alpha}{\alpha_Q}\Big)^{Q-1}\bigg)^{-(1-\frac{\beta}{Q})}c(Q,\beta).
\end{align}
where $c(Q,\beta)>0$ is a constant depending only on $Q$ and $\beta$. Now, by combining \eqref{eq4.5} and \eqref{eq4.9}, we conclude the proof. \qed

The following result, stated in Mahanta et al. \cite[Lemma 2.2]{DTAL-2026}, will be useful in our forthcoming discussions.
\begin{knownlem}\label{lem2.1}
 The embedding $HW^{1,Q}(\mathbb{H}^n)\hookrightarrow L^p (\mathbb{H}^n)$ holds for all $p\geq Q$. Also, the embedding $HW^{1,Q}(\mathbb{H}^n)\hookrightarrow L^p(B_R)$ is compact for all $R>0$ and $p\geq1$.
\end{knownlem}
Next, we derive a weighted embedding that will serve as a key ingredient in the maximization problem stated in Theorem \ref{thm1.4} and in the subsequent analysis.
\begin{lemma}\label{lem4.2}
Let $0\leq\beta<Q$ and $q\geq Q$ be real numbers. Then the weighted embedding $HW^{1,Q}(\mathbb{H}^n)\hookrightarrow L^q(\mathbb{H}^n,r(\xi)^{-\beta}\mathrm{d}\xi)$ is continuous. In fact, this embedding is compact for all $0<\beta<Q$.
\end{lemma}
\begin{proof}
The case $\beta=0$ follows directly from Lemma \ref{lem2.1}. Next, we take $\beta>0$ and $u\in HW^{1,Q}(\mathbb{H}^n)$. Let $t,t'>1$ be such that $\frac{1}{t}+\frac{1}{t'}=1$ and $\beta t'<Q$. Applying H\"{o}lder's inequality and Lemma \ref{lem2.1}, we obtain
\begin{align*}
   \int_{\mathbb{H}^n}\frac{|u|^q}{r(\xi)^\beta}\,\mathrm{d}\xi &\leq \bigg(\int_{B_1}|u|^{qt}\,\mathrm{d}\xi\bigg)^{\frac{1}{t}}\bigg(\int_{B_1}\frac{1}{r(\xi)^{\beta t^\prime}}\,\mathrm{d}\xi\bigg)^{\frac{1}{t^\prime}}+\int_{\mathbb{H}^n\setminus B_1}|u|^q\,\mathrm{d}\xi\\
   &\lesssim\|u\|^q\bigg(\int_{B_1}\frac{1}{r(\xi)^{\beta t^\prime}}\,\mathrm{d}\xi\bigg)^{\frac{1}{t^\prime}}+\|u\|^q\lesssim\|u\|^q<+\infty.
\end{align*}
This shows that $HW^{1,Q}(\mathbb{H}^n)\hookrightarrow L^q(\mathbb{H}^n,r(\xi)^{-\beta}\mathrm{d}\xi)$ holds for all $q\geq Q$ and $\beta\in (0,Q)$. Fixing $\beta\in(0,Q)$, it remains only to prove the compactness of the established embedding. Let $\{u_k\}_k$ be a bounded sequence in $HW^{1,Q}(\mathbb{H}^n)$. By Lemma \ref{lem2.1}, up to a subsequence, there exists $u\in HW^{1,Q}(\mathbb{H}^n)$ such that, for every $R>0$ and $p\ge1$, we have
\begin{align*}
 u_{k_j}\rightharpoonup u\quad\text{in}\quad HW^{1,Q}(\mathbb{H}^n),\quad
 u_{k_j}\to u\quad\text{in}\quad L^p(B_R)\quad\text{and}\quad
  u_{k_j}\to u\quad\text{a.e. in}\quad\mathbb{H}^n
\end{align*}
as $j\to\infty$. By Egoroff's theorem (see, e.g., Royden-Fitzpatrick \cite[{\S}18.1, p. 364]{Royden-Fitzpatrick-1988}), we can notice that for any open ball $B_R$ and $\delta>0$, there exists a measurable set $\mathcal{E}_\delta\subset B_R$ satisfying $|\mathcal{E}_\delta|<\delta$ so that $u_{k_j}\to u$ uniformly in $B_R\setminus \mathcal{E}_\delta$ as $j\to\infty$. The integral is thus decomposed as follows
\begin{align*}
    \lim_{j\to\infty}\int_{\mathbb{H}^n}\frac{|u_{k_j}- u|^q}{r(\xi)^\beta}\,\mathrm{d}\xi=\lim_{R\to\infty}\lim_{\delta\to 0}\lim_{j\to\infty}\bigg(\int_{ \mathcal{E}_\delta}+\int_{B_R\setminus \mathcal{E}_\delta}+\int_{\mathbb{H}^n\setminus B_R}\bigg)\frac{|u_{k_j}- u|^q}{r(\xi)^\beta}\,\mathrm{d}\xi.
\end{align*}
By using H\"{o}lder's inequality and  $HW^{1,Q}(\mathbb{H}^n)\hookrightarrow L^{q t^\prime}(\mathbb{H}^n,r(\xi)^{-\beta t^\prime}\,\mathrm{d}\xi)$, we get
\begin{align*}
 \int_{ \mathcal{E}_\delta}\frac{|u_{k_j}- u|^q}{r(\xi)^\beta}\,\mathrm{d}\xi&\leq   |\mathcal{E}_\delta|^{\frac{1}{t}}\bigg(\int_{ \mathbb{H}^n}\frac{|u_{k_j}- u|^{q t^\prime}}{r(\xi)^{\beta t^\prime }}\,\mathrm{d}\xi \bigg)^{\frac{1}{t^\prime}}\lesssim |\mathcal{E}_\delta|^{\frac{1}{t}}\sup_{j\in\mathbb{N}}\|u_{k_j}- u\|^q\lesssim |\mathcal{E}_\delta|^{\frac{1}{t}}.
\end{align*}
This yields that
$$ \lim_{R\to\infty}\lim_{\delta\to 0}\lim_{j\to\infty}\int_{ \mathcal{E}_\delta}\frac{|u_{k_j}- u|^q}{r(\xi)^\beta}\,\mathrm{d}\xi=0. $$
Moreover, we also have
\begin{align*}
 \int_{B_R\setminus \mathcal{E}_\delta}\frac{|u_{k_j}- u|^q}{r(\xi)^\beta}\,\mathrm{d}\xi\leq \bigg(\sup_{B_R\setminus \mathcal{E}_\delta}|u_{k_j}- u|\bigg)^q\int_{B_R} \frac{\mathrm{d}\xi}{r(\xi)^{\beta}}\lesssim \bigg(\sup_{B_R\setminus \mathcal{E}_\delta}|u_{k_j}- u|\bigg)^q,
\end{align*}
which directly gives
$$ \lim_{R\to\infty}\lim_{\delta\to 0}\lim_{j\to\infty}\int_{B_R\setminus \mathcal{E}_\delta}\frac{|u_{k_j}- u|^q}{r(\xi)^\beta}\,\mathrm{d}\xi=0.$$
It follows from the estimate
\begin{align*}
  \int_{\mathbb{H}^n\setminus B_R}\frac{|u_{k_j}- u|^q}{r(\xi)^\beta}\,\mathrm{d}\xi\leq \frac{1}{R^\beta}  \int_{\mathbb{H}^n}|u_{k_j}- u|^q\,\mathrm{d}\xi\lesssim \frac{1}{R^\beta}\sup_{j\in\mathbb{N}}\|u_{k_j}- u\|^q\lesssim \frac{1}{R^\beta}
\end{align*}
that
$$ \lim_{R\to\infty}\lim_{\delta\to 0}\lim_{j\to\infty} \int_{\mathbb{H}^n\setminus B_R}\frac{|u_{k_j}- u|^q}{r(\xi)^\beta}\,\mathrm{d}\xi=0.$$
Hence, summing up the preceding estimates, we obtain $ u_{k_j}\to u$ in $L^{q }(\mathbb{H}^n,r(\xi)^{-\beta }\,\mathrm{d}\xi)$ as $j\to\infty$, which completes the proof.
\end{proof}
The following result is a direct consequence of the above lemma.
\begin{corollary}\label{cor4.3}
Let $0<\beta<Q$ and $q\geq Q$ be real numbers. Then the best constant $S_q$ in the weighted embedding $HW^{1,Q}(\mathbb{H}^n)\hookrightarrow L^q(\mathbb{H}^n,r(\xi)^{-\beta}\mathrm{d}\xi)$ can be achieved, where
$$ S_q:=\inf_{u\in HW^{1,Q}(\mathbb{H}^n)\setminus\{0\}}\frac{\|u\|}{~\|u\|_{q,\beta}}. $$
\end{corollary}
The result below follows directly from Yang \cite[Lemma 2.1 and Lemma 2.2]{Yang-2012} and will be used repeatedly throughout the sequel.
\begin{knownlem}\label{lem2.2}
The function $\Phi(t)=\exp(t)-\sum_{j=0}^{Q-2}\frac{t^j}{j!}$ is increasing and convex on $[0,+\infty)$. In addition, for all $p\in[1,+\infty)$, there holds
\begin{align*}
\bigg(\exp(t)-\sum_{j=0}^{Q-2}\frac{t^j}{j!}\bigg)^p\leq \exp(p t)-\sum_{j=0}^{Q-2}\frac{(p t)^j}{j!}\quad\text{for all}\quad t\geq0.
\end{align*}
Moreover, the following properties of the function $\Phi$ hold:
\begin{itemize}
      \item [\textnormal{(a)}] For any $\zeta\geq 1$, there holds $\Phi(\zeta t)\geq \zeta^{Q-1}\Phi(t)$ for all $t\geq0$.
      \item [\textnormal{(b)}] For any $0\leq\zeta\leq1$, there holds $\Phi(\zeta t)\leq \zeta^{Q-1}\Phi(t)$ for all $t\geq0$.
  \end{itemize}
\end{knownlem}
Now, we establish a compactness result for the weighted integrability of the Young function in the borderline Sobolev space on the Heisenberg group, as follows.
\begin{proposition}\label{prop4.4}
Assume that $0<\alpha<\alpha_Q$ and $0<\beta<Q$. Let $\{u_k\}_k\subset HW^{1,Q}(\mathbb{H}^n)$ be a sequence such that $u_k\rightharpoonup u$  in $HW^{1,Q}(\mathbb{H}^n)$ as $k\to\infty$ and $\|D_H u_k\|_Q=1$. Then there holds
\begin{align*}
\lim_{k\to\infty}\int_{\mathbb{H}^n}\frac{\Phi\big(\alpha \big(1-\frac{\beta}{Q}\big)|u_k|^{Q^\prime}\big)}{{r(\xi)^\beta}}\,\mathrm{d}\xi=\int_{\mathbb{H}^n}\frac{\Phi\big(\alpha \big(1-\frac{\beta}{Q}\big)|u|^{Q^\prime}\big)}{{r(\xi)^\beta}}\,\mathrm{d}\xi.
\end{align*}
\end{proposition}
\begin{proof}
 Let the given hypotheses be satisfied. Define the map
$\Pi\colon [0,+\infty)\to\mathbb{R}$ by
 \begin{align*}
     \Pi(t):=\Phi\Big(\alpha \Big(1-\frac{\beta}{Q}\Big)t^{Q^\prime}\Big)-\frac{\big(\alpha \big(1-\frac{\beta}{Q}\big)\big)^{Q-1}}{(Q-1)!}~t^Q\quad\text{for all}\quad t\geq0.
 \end{align*}
 It follows from a direct computation that
 \begin{align*}
   \Pi^\prime(t)=\alpha\Big(1-\frac{\beta}{Q}\Big)Q^\prime\,\Phi\Big(\alpha \Big(1-\frac{\beta}{Q}\Big)t^{Q^\prime}\Big)t^{\frac{1}{Q-1}}\quad\text{for all}\quad t\geq0.
 \end{align*}
Keeping the above observation in mind, the mean value theorem, combined with the convexity of $t\mapsto |t|^p$ for $p\in[1,+\infty)$ and Lemma \ref{lem2.2}, yields some $\zeta\in[0,1]$ satisfying
\begin{align*}
    \big|&\Pi(|u_k|)-\Pi(|u|)\big|\leq \Pi^\prime\big(\zeta|u_k|+(1-\zeta)|u|\big)\big||u_k|-|u|\big|\\
    &\lesssim  \Phi\Big(\alpha \Big(1-\frac{\beta}{Q}\Big)\big(\zeta|u_k|+(1-\zeta)|u|\big)^{Q^\prime}\Big)\big(\zeta|u_k|+(1-\zeta)|u|\big)^{\frac{1}{Q-1}}|u_k-u|\\
    &\lesssim  \Phi\Big(\alpha \Big(1-\frac{\beta}{Q}\Big)\Big(\zeta|u_k|^{Q^\prime}+(1-\zeta)|u|^{Q^\prime}\Big)\Big)\big(|u_k|+|u|\big)^{\frac{1}{Q-1}}|u_k-u|\\
    &\lesssim  \bigg[\zeta\Phi\Big(\alpha \Big(1-\frac{\beta}{Q}\Big)|u_k|^{Q^\prime}\Big)+(1-\zeta)\Phi\Big(\alpha \Big(1-\frac{\beta}{Q}\Big)|u|^{Q^\prime}\Big)\bigg]\big(|u_k|+|u|\big)^{\frac{1}{Q-1}}|u_k-u|\\
    & \lesssim  \bigg[\Phi\Big(\alpha \Big(1-\frac{\beta}{Q}\Big)|u_k|^{Q^\prime}\Big)+\Phi\Big(\alpha \Big(1-\frac{\beta}{Q}\Big)|u|^{Q^\prime}\Big)\bigg]\big(|u_k|+|u|\big)^{\frac{1}{Q-1}}|u_k-u|.
\end{align*}
Let $\vartheta_j>1$ for each $j=1,2,3$, be specified later so that $\sum_{j=1}^{3}\frac{1}{\vartheta_j}=1$. Then the above inequality, together with H\"{o}lder's inequality, yields
\begin{align*}
\left|\int_{\mathbb{H}^n}\frac{\Pi(|u_k|)-\Pi(|u|)}{r(\xi)^\beta}\,\mathrm{d}\xi\right|&\lesssim  \left\|\Phi\Big(\alpha \Big(1-\frac{\beta}{Q}\Big)|u_k|^{Q^\prime}\Big)+\Phi\Big(\alpha \Big(1-\frac{\beta}{Q}\Big)|u|^{Q^\prime}\Big)\right\|_{\vartheta_1,\beta}\\
& \qquad\times\left\|\big(|u_k|+|u|\big)^{\frac{1}{Q-1}}\right\|_{\vartheta_2,\beta}\|u_k-u\|_{\vartheta_3,\beta}.
\end{align*}
Note that there also holds $D_H u_k\rightharpoonup D_H u$  in $L^{Q}(\mathbb{H}^n,\mathbb{R}^{2n})$ as $k\to\infty$ and thus, we have $\|D_H u\|_Q\leq1$. Take $\vartheta_1>1$ close to $1$ such that $\alpha\vartheta_1<\alpha_Q$, which is possible because of $\alpha<
\alpha_Q$. Moreover, by utilizing Theorem \ref{thm1.3}, Lemmas \ref{lem2.1} and \ref{lem2.2}, we obtain
\begin{align*}
&\left\|\Phi\Big(\alpha \Big(1-\frac{\beta}{Q}\Big)|u_k|^{Q^\prime}\Big)+\Phi\Big(\alpha \Big(1-\frac{\beta}{Q}\Big)|u|^{Q^\prime}\Big)\right\|_{\vartheta_1,\beta}\\&\leq \left\|\Phi\Big(\alpha \Big(1-\frac{\beta}{Q}\Big)|u_k|^{Q^\prime}\Big)\right\|_{\vartheta_1,\beta}+\left\|\Phi\Big(\alpha \Big(1-\frac{\beta}{Q}\Big)|u|^{Q^\prime}\Big)\right\|_{\vartheta_1,\beta}\\
&\leq\bigg(\int_{\mathbb{H}^n}\frac{\Phi\big(\alpha\vartheta_1 \big(1-\frac{\beta}{Q}\big)|u_k|^{Q^\prime}\big)}{r(\xi)^\beta}\,\mathrm{d}\xi\bigg)^{\frac{1}{\vartheta_1}}+\bigg(\int_{\mathbb{H}^n}\frac{\Phi\big(\alpha\vartheta_1 \big(1-\frac{\beta}{Q}\big)|u|^{Q^\prime}\big)}{r(\xi)^\beta}\,\mathrm{d}\xi\bigg)^{\frac{1}{\vartheta_1}}\\
&\lesssim \|u_k\|_Q^{\frac{Q-\beta}{\vartheta_1}}+\|u\|_Q^{\frac{Q-\beta}{\vartheta_1}}\lesssim 1,
\end{align*}
where we have invoked the fact that $\{u_k\}_k\subset HW^{1,Q}(\mathbb{H}^n)$ is uniformly bounded. Next, for the choice $\vartheta_2$ satisfying $\frac{\vartheta_2}{Q-1}\geq Q$, Lemma \ref{lem4.2} leads to
\begin{align*}
    \left\|\big(|u_k|+|u|\big)^{\frac{1}{Q-1}}\right\|_{\vartheta_2,\beta}&=\||u_k|+|u|\|^{\frac{1}{Q-1}}_{\frac{\vartheta_2}{Q-1},\beta}\leq\left(\|u_k\|_{\frac{\vartheta_2}{Q-1},\beta}+\|u\|_{\frac{\vartheta_2}{Q-1},\beta}\right)^{\frac{1}{Q-1}}\\
    &\leq\|u_k\|^{\frac{1}{Q-1}}_{\frac{\vartheta_2}{Q-1},\beta}+\|u\|^{\frac{1}{Q-1}}_{\frac{\vartheta_2}{Q-1},\beta}\lesssim 1.
\end{align*}
The above estimate is again relying on the uniform boundedness of $\{u_k\}_k\subset HW^{1,Q}(\mathbb{H}^n)$. Finally, taking $\vartheta_3\geq Q$, it follows from Lemma \ref{lem4.2} that $\|u_k-u\|_{\vartheta_3,\beta}\to 0$ as $k\to\infty$.
Hence, from this information, we conclude that
\begin{align}\label{eq3.111}
    \lim_{k\to\infty}\int_{\mathbb{H}^n}
    \frac{\Pi(|u_k|)}{r(\xi)^\beta}\,\mathrm{d}\xi=\int_{\mathbb{H}^n}
    \frac{\Pi(|u|)}{r(\xi)^\beta}\,\mathrm{d}\xi.
\end{align}
By using Lemma \ref{lem4.2}, we also have
\begin{align}\label{eq3.112}
    \lim_{k\to\infty}\int_{\mathbb{H}^n}\frac{|u_{k}|^Q}{r(\xi)^\beta}\,\mathrm{d}\xi= \int_{\mathbb{H}^n}\frac{|u|^Q}{r(\xi)^\beta}\,\mathrm{d}\xi.
\end{align}
In view of \eqref{eq3.111} and \eqref{eq3.112}, the desired result holds, which completes the proof.
\end{proof}
The following elementary result, established next, provides a key ingredient in the proof of Theorem \ref{thm1.4}.
\begin{lemma}\label{lem4.5}
For any $0<\alpha<\alpha_Q$ and $0\leq\beta<Q$, there holds
$$\textnormal{SCMT}(\alpha,\beta)=\sup_{u\in HW^{1,Q}(\mathbb{H}^n),\,\|D_H u\|_{Q}=1,\,\|u\|_Q= 1} \int_{\mathbb{H}^n}\frac{\Phi\big(\alpha \big(1-\frac{\beta}{Q}\big)|u|^{Q^\prime}\big)}{r(\xi)^\beta}\,\mathrm{d}\xi.$$
\end{lemma}
\begin{proof}
We only need to prove that
$$\textnormal{SCMT}(\alpha,\beta)\leq\sup_{u\in HW^{1,Q}(\mathbb{H}^n),\,\|D_H u\|_{Q}=1,\,\|u\|_Q= 1} \int_{\mathbb{H}^n}\frac{\Phi\big(\alpha \big(1-\frac{\beta}{Q}\big)|u|^{Q^\prime}\big)}{r(\xi)^\beta}\,\mathrm{d}\xi,$$
as the reverse inequality is immediate.
Take $u\in HW^{1,Q}(\mathbb{H}^n)$ such that $\|D_H u\|_{Q}\leq 1$. Define the function $w\colon \mathbb{H}^n\to \mathbb{R}$ by $$w(\xi)=\frac{u(\delta_\theta(\xi))}{\|D_H u\|_{Q}}\quad\text{with}\quad \theta=\frac{\|u\|_Q}{\|D_H u\|_{Q}}.$$
This shows that
 \begin{align*}
    \|w\|_Q=\theta^{-1}\bigg(\frac{\| u\|_Q}{\|D_H u\|_{Q}}\bigg) =1\quad\text{and}\quad\|D_H w\|_Q=1.
 \end{align*}
 Define \begin{align*}
    \quad {J}_{\alpha,\beta}(v)=\int_{\mathbb{H}^n}\frac{\Phi\big(\alpha \big(1-\frac{\beta}{Q}\big)|v|^{Q^\prime}\big)}{r(\xi)^\beta}\mathrm{d}\xi\quad\text{for all}\quad v\in HW^{1,Q}(\mathbb{H}^n).
\end{align*}
 Moreover, by applying Lemma \ref{lem2.2}, it follows that
 \begin{align*}
     {J}_{\alpha,\beta}(w)&=\frac{1}{\theta^{Q-\beta}}{J}_{\alpha,\beta}\bigg(\frac{u}{\|D_H u\|_Q}\bigg)\geq \frac{1}{\theta^{Q-\beta}\|D_H u\|_Q^Q}{J}_{\alpha,\beta}(u)\geq \frac{1}{\|u\|_Q^{Q-\beta}}{J}_{\alpha,\beta}(u).
 \end{align*}
Hence, the desired inequality follows, and the proof is complete.
\end{proof}
Finally, we establish the existence of a maximizer for the subcritical weighted Moser--Trudinger inequality, thus concluding the section.
\subsection{{Proof of Theorem \ref{thm1.4}}}
Let $\{u_k\}_k\subset HW^{1,Q}(\mathbb{H}^n)\setminus\{0\}$ be a maximizing sequence for $\textnormal{SCMT}(\alpha,\beta)$. By using Lemma \ref{lem4.5}, we can assume that $ \| u_k\|_Q=1$, $\| D_H u_k\|_Q=1$ and \begin{align*}
    \lim_{k\to\infty} \int_{\mathbb{H}^n}\frac{\Phi\big(\alpha \big(1-\frac{\beta}{Q}\big)|u_k|^{Q^\prime}\big)}{r(\xi)^\beta}\,\mathrm{d}\xi=\textnormal{SCMT}(\alpha,\beta).
\end{align*}
This shows that $\{u_k\}_k$ is a bounded sequence in $ HW^{1,Q}(\mathbb{H}^n)$. Also, by employing Lemma \ref{lem4.2}, up to a subsequence still denoted by itself, there exists $u\in HW^{1,Q}(\mathbb{H}^n)$ such that for all $q\geq Q$, we have
\begin{align}\label{eq4.10}
u_{k}\rightharpoonup u~~\text{in}~~ HW^{1,Q}(\mathbb{H}^n),~ u_{k}\to u~~\text{in}~~ L^q(\mathbb{H}^n,r(\xi)^{-\beta}\mathrm{d}\xi)~~\text{and}~~ u_{k}\to u~~\text{a.e. in}\ \mathbb{H}^n
\end{align}
as $k\to\infty$. Now, by using \eqref{eq4.10}, we obtain that  $ u_k\rightharpoonup u$ in $L^{Q}(\mathbb{H}^n)$ as $k\to\infty$ and $D_H u_k\rightharpoonup D_H u$  in $L^{Q}(\mathbb{H}^n,\mathbb{R}^{2n})$ as $k\to\infty$. Moreover, by the weak lower semicontinuity of the norm, we have
\begin{align}\label{eq4.11}
  \| u\|_Q\leq1\quad\text{and}\quad \|D_H u\|_Q\leq1.
\end{align}
In view of \eqref{eq4.10}, Proposition \ref{prop4.4} further implies that
\begin{align}\label{eq4.12}
\textnormal{SCMT}(\alpha,\beta)=\int_{\mathbb{H}^n}\frac{\Phi\big(\alpha \big(1-\frac{\beta}{Q}\big)|u|^{Q^\prime}\big)}{r(\xi)^\beta}\,\mathrm{d}\xi.
\end{align}
The positivity of $\textnormal{SCMT}(\alpha,\beta)$ ensures that $u\neq 0$. Indeed, if $u=0$, then \eqref{eq4.12} would imply $\textnormal{SCMT}(\alpha,\beta)=0$, which is contrary to $\textnormal{SCMT}(\alpha,\beta)>0$. Next, by recalling the definition of $\textnormal{SCMT}(\alpha,\beta)$, we obtain from \eqref{eq4.11} and \eqref{eq4.12} that
\begin{align*}
 \textnormal{SCMT}(\alpha,\beta)\geq \frac{1}{\|u\|_Q^{Q-\beta}} \int_{\mathbb{H}^n}\frac{\Phi\big(\alpha \big(1-\frac{\beta}{Q}\big)|u|^{Q^\prime}\big)}{r(\xi)^\beta}\,\mathrm{d}\xi= \frac{\textnormal{SCMT}(\alpha,\beta)}{\|u\|_Q^{Q-\beta}}\geq  \textnormal{SCMT}(\alpha,\beta).
\end{align*}
Hence, $u$ attains the supremum
$\textnormal{SCMT}(\alpha,\beta)$, and $\|u\|_Q=1$. Thus, $u$ is a maximizer, and the desired result now follows. This completes the proof. \qed
\section{The Caffarelli--Kohn--Nirenberg Interpolation-Type Inequality}\label{sec5}
In this section, we establish Theorem \ref{thm1.5} that concerns Caffarelli--Kohn--Nirenberg interpolation-type inequalities on the Heisenberg group. The proof covers an extended range of parameters in the borderline Sobolev setting. We then analyze the corresponding sharp constants in the asymptotic regime and relate them to those of the subcritical weighted Moser--Trudinger inequalities.
\subsection{ {Proof of Theorem \ref{thm1.5}}}
Let $\alpha_0$ be the supremum of all $\alpha>0$ such that for all $u\in HW^{1,Q}(\mathbb{H}^n)$ satisfying $\|D_H u\|_Q\leq 1$, one has
\begin{align*}
\int_{\mathbb{H}^n}\frac{\Phi\big(\alpha |u|^{Q^\prime}\big)}{r(\xi)^\beta}\,\mathrm{d}\xi
\leq C(Q,\beta)\bigg(\int_{\mathbb{H}^n}|u|^Q \,\mathrm{d}\xi\bigg)^{1-\frac{\beta}{Q}},
\end{align*}
where $C(Q,\beta)>0$ is a constant depending only on $Q$ and $\beta$. By the subcritical weighted Moser--Trudinger inequality established in Theorem \ref{thm1.3}, also proved by  Lam--Lu--Tang \cite[Theorem 1.1]{Lam-Lu-Tang-2014}, we obtain
$\alpha_0=\alpha_Q(1-\frac{\beta}{Q})$. Thus, for every $\alpha\in(0,\alpha_0)$, there exists a constant $C(Q,\beta)>0$ such that for all $u\in HW^{1,Q}(\mathbb{H}^n)\setminus\{0\}$ and all integers $k\geq Q-1$, there holds
$$ C(Q,\beta)\bigg(\frac{\|u\|_Q}{\|D_H u\|_Q}\bigg)^{Q-\beta} \geq \int_{\mathbb{H}^n}\Phi\bigg(\alpha\bigg|\frac{u}{\|D_H u\|_Q}\bigg|^{Q^\prime}\bigg)\,\frac{\mathrm{d}\xi}{r(\xi)^\beta}\geq\frac{\alpha^k} {k!}\bigg(\frac{\|u\|_{Q^\prime k,\beta}}{\|D_H u\|_Q}\bigg)^{Q^\prime k}.$$
As a result, we deduce that
\begin{align}\label{eq5.1}
  \|u\|_{Q^\prime k,\beta}\leq \bigg(\frac{k! C(Q,\beta)}{\alpha^k}\bigg)^{\frac{1}{Q^\prime k}} \|u\|_Q^{\frac{Q-\beta}{Q^\prime k}} \|D_H u\|_Q^{1-\frac{Q-\beta}{Q^\prime k}}
\end{align}
holds for all $u\in HW^{1,Q}(\mathbb{H}^n)$ and all integers $k\geq Q-1$. This shows that for all $\vartheta\in[Q,+\infty)$ there always exists an integer $k\geq Q-1$ satisfying $Q^\prime k\leq \vartheta<Q^\prime(k+1)$. Hence, we may take $\zeta\in(0,1]$ so that
$\frac{1}{\vartheta}=\frac{\zeta}{Q^\prime k}+\frac{1-\zeta}{Q^\prime(k+1)}$. It then follows from H\"{o}lder's inequality that
\begin{align}\label{eq5.2}
  \|u\|_{\vartheta,\beta}\leq \|u\|_{Q^\prime k,\beta}^\zeta \|u\|_{Q^\prime(k+1),\beta}^{1-\zeta}.
\end{align}
The estimates \eqref{eq5.1} and \eqref{eq5.2} together imply that, for all $u\in HW^{1,Q}(\mathbb{H}^n)$ and $\vartheta\in[Q,+\infty)$, there holds
\begin{align}\label{eq5.3}
 \|u\|_{\vartheta,\beta}\leq ((k+1)!)^{\frac{1}{\vartheta}}\alpha^{-\frac{1}{Q^\prime}} C(Q,\beta)^{\frac{1}{\vartheta}}\|u\|_Q^{\frac{Q-\beta}{\vartheta}} \|D_H u\|_Q^{1-\frac{Q-\beta}{\vartheta}}.
\end{align}
The assumption $\vartheta\geq Q^\prime k$ implies that
$((k+1)!)^{\frac{1}{\vartheta}}=\Gamma(k+2)^{\frac{1}{\vartheta}}\leq \big(\Gamma(\frac{\vartheta}{Q^\prime}+2)\big)^{\frac{1}{\vartheta}}.$
In view of \eqref{eq5.3}, the desired inequality \eqref{eq1.3} follows upon using the estimate given below $$C(Q,\beta,\vartheta)\approx \bigg(\Gamma\Big(\frac{\vartheta}{Q^\prime}+2\Big)\bigg)^{\frac{1}{\vartheta}} \alpha^{-\frac{1}{Q^\prime}} C(Q,\beta)^{\frac{1}{\vartheta}}.$$
The next step is to show that, for some $\lambda>0$, there holds
$C(Q,\beta,\vartheta)\simeq\lambda\vartheta^{\frac{1}{Q^\prime}}$ as $\vartheta\to\infty$. More precisely, there exists $\vartheta_0\geq Q$ such that for all $\vartheta\geq\vartheta_0$, we have
\begin{align}\label{eq5.4}
    \|u\|_{\vartheta,\beta}\leq \lambda \vartheta^{\frac{1}{Q^\prime}}\|u\|_Q^{\frac{Q-\beta}{\vartheta}} \|D_H u\|_Q^{1-\frac{Q-\beta}{\vartheta}},
\end{align}
for all $u\in HW^{1,Q}(\mathbb{H}^n)$. Due to Stirling’s asymptotic formula, one has
$$ \bigg(\Gamma\Big(\frac{\vartheta}{Q^\prime}+2\Big)\bigg)^{\frac{1}{\vartheta}}=(1+o_\vartheta(1))\bigg(\frac{\vartheta}{Q^\prime\operatorname{\exp}(1)}\bigg)^{\frac{1}{Q^\prime}}\quad\text{as}\quad \vartheta\to\infty.$$
By combining the above information, all together leads to the following estimate
$$\|u\|_{\vartheta,\beta}\leq (1+o_\vartheta(1))\bigg(\frac{\vartheta}{\alpha Q^\prime\operatorname{\exp}(1)}\bigg)^{\frac{1}{Q^\prime}}\|u\|_Q^{\frac{Q-\beta}{\vartheta}} \|D_H u\|_Q^{1-\frac{Q-\beta}{\vartheta}}\quad\text{as}\quad \vartheta\to\infty,$$
for all $u\in HW^{1,Q}(\mathbb{H}^n)$. This gives the validity of \eqref{eq5.4}. In addition, let $\boldsymbol{\Lambda}_{Q,\beta}$ be the infimum of all $\lambda>0$ such that \eqref{eq5.4} holds. Evidently, we obtain from the preceding estimate that $\boldsymbol{\Lambda}_{Q,\beta} \leq (\alpha Q^\prime\operatorname{\exp}(1))^{-\frac{1}{Q^\prime}}$ for all $\alpha\in(0,\alpha_0)$. Now, taking the limit $\alpha\to\alpha_0^{-}$ and using the value of $\alpha_0$ yields
\begin{align}\label{eq5.5}
 \boldsymbol{\Lambda}_{Q,\beta} \leq \bigg(\frac{1}{\alpha_Q\big(1-\frac{\beta}{Q}\big) Q^\prime\operatorname{\exp}(1)}\bigg)^{\frac{1}{Q^\prime}}.
\end{align}
To finish the proof, it suffices to show the reverse inequality of \eqref{eq5.5}. Let $\lambda>\boldsymbol{\Lambda}_{Q,\beta}$ be given. It follows from the definition of $\boldsymbol{\Lambda}_{Q,\beta}$ that there exists $\vartheta_0\geq Q$ such that \eqref{eq5.4} holds for all $\vartheta\geq\vartheta_0$. Now, we take $\alpha>0$, with its precise value to be determined later. Thus, for all $u\in HW^{1,Q}(\mathbb{H}^n)$ satisfying $\|D_Hu\|_Q\leq1$, we can write
\begin{align*}
    \int_{\mathbb{H}^n}\frac{\Phi\big(\alpha |u|^{Q^\prime}\big)}{r(\xi)^\beta}\,\mathrm{d}\xi=:I+J,
\end{align*}
where
\begin{align*}
    I=\int_{\mathbb{H}^n}\sum_{\substack{k\in\mathbb{N}\\ Q\leq Q'k<\vartheta_0}}\frac{\alpha^k}{k!}\frac{|u|^{Q'k}}{r(\xi)^\beta}\,\mathrm{d}\xi\qquad\text{and}\qquad J=\int_{\mathbb{H}^n}\sum_{\substack{k\in\mathbb{N}\\  Q'k\geq\vartheta_0}}\frac{\alpha^k}{k!}\frac{|u|^{Q'k}}{r(\xi)^\beta}\,\mathrm{d}\xi.
\end{align*}
To evaluate $I$, we begin by noting that, for each integer $k$ with $Q\leq Q'k<\vartheta_0$, one can find $\zeta\in(0,1]$ such that
$\frac{1}{Q^\prime k}=\frac{\zeta}{Q}+\frac{1-\zeta}{\vartheta_0}$. Then, the H\"older's inequality yields
\begin{align*}
 \|u\|_{Q^\prime k,\beta} \leq \|u\|_{Q,\beta}^\zeta \|u\|_{\vartheta_0,\beta}^{1-\zeta},
\end{align*}
for all $u\in HW^{1,Q}(\mathbb{H}^n)$. Recall that \eqref{eq1.3} and \eqref{eq5.4} remain valid when $\vartheta$ replaced by $Q$ and $\vartheta_0$. By applying \eqref{eq1.3} and \eqref{eq5.4} to the preceding estimate, we obtain for all $u\in HW^{1,Q}(\mathbb{H}^n)$ with $\|D_Hu\|_Q\leq1$ and each integer $k$ with $Q\leq Q'k<\vartheta_0$ that
\begin{align*}
 \|u\|_{Q^\prime k,\beta} &\leq\Big(C(Q,\beta) \|u\|_Q^{\frac{Q-\beta}{Q}}\|D_Hu\|_Q^{1-\frac{Q-\beta}{Q}}\Big)^\zeta \Big(\lambda \vartheta_0^{\frac{1}{Q^\prime}}\|u\|_Q^{\frac{Q-\beta}{\vartheta_0}} \|D_H u\|_Q^{1-\frac{Q-\beta}{\vartheta_0}}\Big)^{1-\zeta}\\
 &\leq C(Q,\beta)^\zeta \lambda^{1-\zeta} \vartheta_0^{\frac{{1-\zeta}}{Q^\prime}}\|u\|_Q^{\zeta(\frac{Q-\beta}{Q})}\|u\|_Q^{(1-\zeta)(\frac{Q-\beta}{\vartheta_0})}=C(Q,\beta)^\zeta \lambda^{1-\zeta} \vartheta_0^{\frac{{1-\zeta}}{Q^\prime}}\|u\|_Q^{\frac{Q-\beta}{Q^\prime k}}.
\end{align*}
Denoting by $D$ the maximum of the finitely many constants $C(Q,\beta)^{Q^\prime\zeta} \lambda^{Q^\prime(1-\zeta)} \vartheta_0^{1-\zeta}$ (note that $\zeta$ depends on $k$), we get
\begin{align}\label{eq5.6}
    I=\sum_{\substack{k\in\mathbb{N}\\ Q\leq Q'k<\vartheta_0}}\frac{\alpha^k}{k!}\|u\|_{Q^\prime k,\beta}^{Q^\prime k}\leq \sum_{\substack{k\in\mathbb{N}\\ Q\leq Q'k<\vartheta_0}}\frac{(\alpha D)^k}{k!}\|u\|_Q^{Q-\beta}\leq \operatorname{\exp}(\alpha D)\|u\|_Q^{Q-\beta}.
\end{align}
To this end, we compute $J$. By applying the monotone convergence theorem as stated in Royden--Fitzpatrick \cite[{\S}4.3, p. 83]{Royden-Fitzpatrick-1988} and \eqref{eq5.4}, we obtain for all $u\in HW^{1,Q}(\mathbb{H}^n)$ with $\|D_Hu\|_Q\leq1$ that
\begin{equation}\label{eq5.7}
\begin{aligned}
  J&=\sum_{\substack{k\in\mathbb{N}\\  Q'k\geq\vartheta_0}}\frac{\alpha^k}{k!}\|u\|_{Q^\prime k,\beta}^{Q^\prime k}\leq  \sum_{\substack{k\in\mathbb{N}\\  Q'k\geq\vartheta_0}}\frac{\alpha^k}{k!}\Big(\lambda \big(Q^\prime k\big)^{\frac{1}{Q^\prime}}\|u\|_Q^{\frac{Q-\beta}{Q^\prime k}} \|D_H u\|_Q^{1-\frac{Q-\beta}{Q^\prime k}}\Big)^{Q^\prime k}  \\
  &\leq \sum_{\substack{k\in\mathbb{N}\\  Q'k\geq\vartheta_0}}\frac{\alpha^k}{k!}\Big(\lambda \big(Q^\prime k\big)^{\frac{1}{Q^\prime}}\|u\|_Q^{\frac{Q-\beta}{Q^\prime k}}\Big)^{Q^\prime k} =\sum_{\substack{k\in\mathbb{N}\\  Q^\prime k\geq\vartheta_0}}\frac{k^k}{k!}\big(\alpha\lambda^{Q^\prime }Q^\prime \big)^k\|u\|_Q^{Q-\beta}.
\end{aligned}
\end{equation}
By d'Alembert's ratio test, the series appearing in \eqref{eq5.7} converges provided that $\alpha\lambda^{Q^\prime }Q^\prime<\operatorname{exp}(-1)$. In particular, there holds $\alpha\in\big(0,(\lambda^{Q^\prime }Q^\prime)^{-1}\operatorname{exp}(-1)\big)$. It follows from \eqref{eq5.6}, \eqref{eq5.7} and the definition of $\alpha_0$ that $\alpha_0\geq (\lambda^{Q^\prime }Q^\prime\operatorname{exp}(1))^{-1}$ for all $\lambda>\boldsymbol{\Lambda}_{Q,\beta}$. Passing to the limit $\lambda\to\boldsymbol{\Lambda}_{Q,\beta}^{+}$ and using the value of $\alpha_0$, we obtain
\begin{align}\label{eq5.8}
 \boldsymbol{\Lambda}_{Q,\beta} \geq \bigg(\frac{1}{\alpha_Q\big(1-\frac{\beta}{Q}\big) Q^\prime\operatorname{\exp}(1)}\bigg)^{\frac{1}{Q^\prime}}.
\end{align}
Hence, the desired conclusion follows from \eqref{eq5.5} and \eqref{eq5.8}. This completes the proof. \qed

\section{The Weighted Moser--Trudinger Inequality under Inhomogeneous Norm Constraints and Its Maximizer}\label{sec6}
In this section, inspired by the rearrangement-free approach of Lam--Lu--Zhang \cite{Lam-Lu-Zhang-2017}, a weighted Moser--Trudinger inequality under inhomogeneous norm constraints on the Heisenberg group is derived using Theorem \ref{thm1.3}. Its equivalence with the subcritical weighted Moser--Trudinger inequality is subsequently established. In light of the corresponding Euclidean analogs in \cite{Lam-2017, Lam-Lu-Zhang-2019}, this equivalence yields the existence of maximizers in the subcritical regime by avoiding radial reduction of Sobolev functions. Now, we begin with the following elementary result.
 \begin{lemma}\label{lem6.1}
 Let $\textnormal{CMT}_{a,b}(\alpha_Q,\beta)<+\infty$ be satisfied. Then, for every $\alpha\in(0,\alpha_Q)$, one obtains $\textnormal{SCMT}(\alpha,\beta)<+\infty$, together with the estimate
$$\textnormal{SCMT}(\alpha,\beta)\leq\Bigg(\cfrac{\big(\frac{\alpha}{\alpha_Q}\big)^{(\frac{Q-1}{Q})b}}{1-\big(\frac{\alpha}{\alpha_Q}\big)^{(\frac{Q-1}{Q})a}}\Bigg)^{\frac{Q}{b}(1-\frac{\beta}{Q})}\textnormal{CMT}_{a,b}(\alpha_Q,\beta). $$
In particular, it holds
$$\textnormal{SCMT}(\alpha,\beta)\leq\Bigg(\cfrac{\big(\frac{\alpha}{\alpha_Q}\big)^{Q-1}}{1-\big(\frac{\alpha}{\alpha_Q}\big)^{Q-1}}\Bigg)^{1-\frac{\beta}{Q}}\textnormal{CMT}_{Q,Q}(\alpha_Q,\beta).$$
\end{lemma}
\begin{proof}
In view of Lemma \ref{lem4.1} and the definition of $\textnormal{SCMT}(\alpha,\beta)$, under the given hypothesis, we may assume $u\in HW^{1,Q}(\mathbb{H}^n)$ with $\|D_H u\|_Q\leq 1$ and $\| u\|_Q= 1$. Now, we define $w\colon\mathbb{H}^n\to\mathbb{R}$ by
\begin{align*}
  w(\xi)=\Big(\frac{\alpha}{\alpha_Q}\Big)^{\frac{Q-1}{Q}}u(\delta_\theta(\xi)) \qquad\text{with}\qquad \theta=\Bigg(\cfrac{\big(\frac{\alpha}{\alpha_Q}\big)^{(\frac{Q-1}{Q})b}}{1-\big(\frac{\alpha}{\alpha_Q}\big)^{(\frac{Q-1}{Q})a}}\Bigg)^{\frac{1}{b}}
\end{align*}
and
\begin{align*}
     {J}_{\alpha,\beta}(v)=\int_{\mathbb{H}^n}\frac{\Phi\big(\alpha \big(1-\frac{\beta}{Q}\big)|v|^{Q^\prime}\big)}{r(\xi)^\beta}\,\mathrm{d}\xi\quad\text{for all}\quad v\in HW^{1,Q}(\mathbb{H}^n).
\end{align*}
By direct calculations, we obtain
$$ \|D_H w\|_Q^a=\Big(\frac{\alpha}{\alpha_Q}\Big)^{(\frac{Q-1}{Q})a}\|D_H u\|^a_Q\leq \Big(\frac{\alpha}{\alpha_Q}\Big)^{(\frac{Q-1}{Q})a}~~\text{and}~~\|w\|_Q^b=1- \Big(\frac{\alpha}{\alpha_Q}\Big)^{(\frac{Q-1}{Q})a}. $$
It follows that $\|D_H w\|_Q^a+\| w\|_Q^b\leq 1$. Hence, by combining the above estimates, we obtain from the definition of $\textnormal{CMT}_{a,b}(\alpha_Q,\beta)$ that
\begin{align*}
 {J}_{\alpha,\beta}(u)=\theta^{Q-\beta}J_{{\alpha_Q},\beta}(w)\leq\Bigg(\cfrac{\big(\frac{\alpha}{\alpha_Q}\big)^{(\frac{Q-1}{Q})b}}{1-\big(\frac{\alpha}{\alpha_Q}\big)^{(\frac{Q-1}{Q})a}}\Bigg)^{\frac{Q}{b}(1-\frac{\beta}{Q})}\textnormal{CMT}_{a,b}(\alpha_Q,\beta).
\end{align*}
The above inequality ensures the desired result, which completes the proof.
\end{proof}

The preceding results enable us to complete the proof of Theorem \ref{thm1.6} as follows.
\subsection{{Proof of Theorem \ref{thm1.6}}} It is enough to consider $\alpha=\alpha_Q$, as the case $\alpha\in(0,\alpha_Q)$ follows by the same argument. To proceed, we assume that $b\leq Q$ and $u\in HW^{1,Q}(\mathbb{H}^n)\setminus\{0\}$ with $\|D_H u\|^a_{Q}+\| u\|^b_Q\leq 1$. Indeed, if possible, let us choose $\|D_H u\|_Q=\vartheta\in(0,1)$ and $\| u\|_Q\leq(1-\vartheta^a)^{\frac{1}{b}}$. Take $\frac{1}{2}<\vartheta<1$. Define $w\colon\mathbb{H}^n\to\mathbb{R}$ by
$$w(\xi)=\frac{u(\delta_\theta(\xi))}{\vartheta}\qquad\text{with}\qquad\theta=\frac{(1-\vartheta^a)^{\frac{1}{b}}}{\vartheta}. $$
By straightforward calculations, we obtain
$$ \|D_H w\|_Q=\vartheta^{-1}\|D_H u\|_Q=1 \qquad\text{and}\qquad \| w\|_Q=(1-\vartheta^a)^{-\frac{1}{b}}\|u\|_Q\leq1.  $$
Next, applying Theorem \ref{thm1.3} and using $b\leq Q$, one sees that
\begin{align*}
  \int_{\mathbb{H}^n}&\frac{\Phi\big(\alpha_Q \big(1-\frac{\beta}{Q}\big)|u|^{Q^\prime}\big)}{r(\xi)^\beta}\,\mathrm{d}\xi  =\int_{\mathbb{H}^n}\frac{\Phi\big(\alpha_Q \big(1-\frac{\beta}{Q}\big)|u(\delta_\theta(\xi)|^{Q^\prime}\big)}{r(\delta_\theta(\xi))^\beta}\,\mathrm{d}(\delta_\theta(\xi))\\
  &=\theta ^{Q-\beta}\int_{\mathbb{H}^n}\frac{\Phi\big(\vartheta^{Q^\prime}\alpha_Q \big(1-\frac{\beta}{Q}\big)|w|^{Q^\prime}\big)}{r(\xi)^\beta}\,\mathrm{d}\xi\\
 &\leq \theta ^{Q-\beta} \textnormal{SCMT}\big(\vartheta^{Q^\prime}\alpha_Q,\beta\big)\leq \bigg(\frac{(1-\vartheta^a)^{\frac{1}{b}}}{\vartheta}\bigg)^{Q-\beta}\frac{C(Q,\beta)}{(1-\vartheta^Q)^{1-\frac{\beta}{Q}}}\\
 &\leq \bigg(\frac{(1-\vartheta^a)^{\frac{Q}{b}}}{1-\vartheta^Q}\bigg)^{1-\frac{\beta}{Q}}C(Q,\beta)\leq C(Q,a,b,\beta),
\end{align*}
where $C(Q,a,b,\beta)>0$ is a constant depending only on $Q,\,a,\,b$ and $\beta$. Now, we take $0<\vartheta\leq\frac{1}{2}$ and define $w\colon\mathbb{H}^n\to\mathbb{R}$ by $w(\xi)=2u(\delta_2(\xi))$ for all $\xi\in \mathbb{H}^n$.
It follows that $\|D_H w\|_Q=2\|D_H u\|_Q\leq 1$ and $\|w\|_Q=\| u\|_Q\leq 1$. Then we invoke Theorem \ref{thm1.3} to obtain
\begin{align*}
   \int_{\mathbb{H}^n}&\frac{\Phi\big(\alpha_Q \big(1-\frac{\beta}{Q}\big)|u|^{Q^\prime}\big)}{r(\xi)^\beta}\,\mathrm{d}\xi =\int_{\mathbb{H}^n}\frac{\Phi\big(\alpha_Q \big(1-\frac{\beta}{Q}\big)|u(\delta_2(\xi)|^{Q^\prime}\big)}{r(\delta_2(\xi))^\beta}\,\mathrm{d}(\delta_2(\xi))\\
   &= 2 ^{Q-\beta}\int_{\mathbb{H}^n}\frac{\Phi\big(2^{-Q^\prime}\alpha_Q \big(1-\frac{\beta}{Q}\big)|w|^{Q^\prime}\big)}{r(\xi)^\beta}\,\mathrm{d}\xi\leq 2 ^{Q-\beta} \text{SCMT}\big(2^{-Q^\prime}\alpha_Q,\beta\big)\\
   &\leq \frac{2 ^{Q-\beta}}{\big(1-{2^{-Q}}\big)^{1-\frac{\beta}{Q}}}\,C(Q,\beta)\leq C(Q,\beta),
\end{align*}
where $C(Q,\beta)>0$ is a constant depending only on   $Q$ and $\beta$. Thus, upon combining the above estimates, we infer that $\textnormal{CMT}_{a,b}(\alpha_Q,\beta)<+\infty$ for all $b\leq Q$.

Next, we need to prove the sharpness of $\alpha_Q$, that is, $\textnormal{CMT}_{a,b}(\alpha,\beta)=+\infty$ for all $\alpha>\alpha_Q$. For this purpose, we take the same Moser's function sequence $\{w_k\}_k$, which is defined in \eqref{eq4.6}. It follows from \eqref{eq4.8} that
$$\|D_H w_k\|_Q=1\qquad\text{and}\qquad \| w_k\|_Q=o_k(1) \quad\text{as}\quad k\to\infty.$$
Moreover, let $\{\zeta_k\}_k\subset (0,1)$ be a real sequence satisfying $\zeta_k^a+\zeta_k^b\| w_k\|^b_Q=1$. Thus, by its definition, $\zeta_k\to 1$ as $k\to\infty$. For each $k\in\mathbb{N}$, set
$ u_k(\xi)=\zeta_k w_k(\xi)$ for all $\xi\in\mathbb{H}^n$. By the construction, it is immediate that $\|D_H u_k\|^a_Q+\| u_k\|^b_Q=1$. Taking $\alpha>\alpha_Q$, we obtain
\begin{align*}
   \int_{\mathbb{H}^n}&\frac{\Phi\big(\alpha \big(1-\frac{\beta}{Q}\big)|u_k|^{Q^\prime}\big)}{r(\xi)^\beta}\,\mathrm{d}\xi\geq \int_{\{\xi\in\mathbb{H}^n\colon~r(\xi)<sk^{-1}\}}\frac{\Phi\big(\alpha \big(1-\frac{\beta}{Q}\big)|u_k|^{Q^\prime}\big)}{r(\xi)^\beta}\,\mathrm{d}\xi\\
   &= \int_\Sigma\int_{0}^{sk^{-1}}\Phi\left(\alpha \Big(1-\frac{\beta}{Q}\Big)\Big|\zeta_k c^{-\frac{1}{Q}}_Q(\ln k)^{\frac{Q-1}{Q}}\Big|^{Q^\prime}\right)
   \eta^{Q-\beta-1}\,\mathrm{d}\eta\,\mathrm{d}\mu(\xi^\ast)\\
   &\gtrsim \bigg(\displaystyle\int_\Sigma\,\mathrm{d}\mu(\xi^\ast)\bigg)\cfrac{(sk^{-1})^{Q-\beta}}{Q-\beta}\,\operatorname{\exp}\Bigg(\alpha\Big(1-\frac{\beta}{Q}\Big)\bigg|\zeta_k c^{-\frac{1}{Q}}_Q(\ln k)^{\frac{Q-1}{Q}}\bigg|^{Q^\prime}\Bigg)\\
   &\gtrsim \operatorname{\exp}\bigg((Q-\beta)\ln k\bigg(\frac{\alpha \zeta_k^{Q^\prime}}{\alpha_Q}-1\bigg)\bigg)\to+\infty\quad\text{as}\quad k\to\infty.
\end{align*}
This concludes that $\textnormal{CMT}_{a,b}(\alpha,\beta)=+\infty$ for all $\alpha>\alpha_Q$, and consequently establishes the sharpness of $\alpha_Q$.
In addition, by applying Lemma \ref{lem6.1}, there holds
$$ \sup_{\alpha\in(0,\alpha_Q)}\Bigg(\cfrac{1-\big(\frac{\alpha}{\alpha_Q}\big)^{(\frac{Q-1}{Q})a}}{\big(\frac{\alpha}{\alpha_Q}\big)^{(\frac{Q-1}{Q})b}}\Bigg)^{\frac{Q}{b}(1-\frac{\beta}{Q})}\text{SCMT}(\alpha,\beta)\leq \text{CMT}_{a,b}(\alpha_Q,\beta).$$
Let $\{v_k\}_k\subset HW^{1,Q}(\mathbb{H}^n)\setminus\{0\}$ be a maximizing sequence for $\text{CMT}_{a,b}(\alpha_Q,\beta)$. Evidently, we deduce that
$\|D_H v_k\|^a_Q+\| v_k\|^b_Q\leq1$
and \begin{align*}
    \lim_{k\to\infty} \int_{\mathbb{H}^n}\frac{\Phi\big(\alpha_Q \big(1-\frac{\beta}{Q}\big)|v_k|^{Q^\prime}\big)}{r(\xi)^\beta}\,\mathrm{d}\xi=\text{CMT}_{a,b}(\alpha_Q,\beta).
\end{align*}
For $k\in\mathbb{N}$, we define
$$ \ell_k(\xi)=\frac{v_k(\delta_{\theta_k} (\xi))}{\|D_H v_k\|_Q}\quad\text{with}\quad\theta_k=\bigg(\frac{1-\|D_H v_k\|^a_{Q}}{\|D_H v_k\|^b_Q}\bigg)^{\frac{1}{b}} .$$
This yields that
\begin{align*}
 \|D_H \ell_k\|_Q=1 \quad\text{and} \quad\|\ell_k\|_Q=\theta_k^{-1}\,\frac{\|v_k\|_Q}{\|D_H v_k\|_Q}\leq \theta_k^{-1}\,\frac{\big(1-\|D_Hv_k\|^a_Q\big)^{\frac{1}{b}}}{\|D_H v_k\|_Q}=1.
\end{align*}
Recalling the definition of $\textnormal{SCMT}(\alpha,\beta)$, we obtain
\begin{align*}
\int_{\mathbb{H}^n}&\frac{\Phi\big(\alpha_Q \big(1-\frac{\beta}{Q}\big)|v_k|^{Q^\prime}\big)}{r(\xi)^\beta}\,\mathrm{d}\xi =\int_{\mathbb{H}^n}\frac{\Phi\big(\alpha_Q \big(1-\frac{\beta}{Q}\big)|v_k(\delta_{\theta_k}(\xi)|^{Q^\prime}\big)}{r(\delta_{\theta_k}(\xi))^\beta}\,\mathrm{d}(\delta_{\theta_k}(\xi))\\
  &=\theta_k^{Q-\beta}\int_{\mathbb{H}^n}\frac{\Phi\big(\alpha_Q \|D_H v_k\|_Q^{Q^\prime}\big(1-\frac{\beta}{Q}\big)|\ell_k|^{Q^\prime}\big)}{r(\xi)^\beta}\,\mathrm{d}\xi\leq \theta_k ^{Q-\beta} \textnormal{SCMT}\big(\alpha_Q \|D_H v_k\|_Q^{Q^\prime},\beta\big)\\
 &=\bigg(\frac{1-\|D_H v_k\|^a_{Q}}{\|D_H v_k\|^b_Q}\bigg)^{\frac{Q}{b}(1-\frac{\beta}{Q})}\textnormal{SCMT}\big(\alpha_Q \|D_H v_k\|_Q^{Q^\prime},\beta\big)\\
 &\leq \sup_{\alpha\in(0,\alpha_Q)}\Bigg(\cfrac{1-\big(\frac{\alpha}{\alpha_Q}\big)^{(\frac{Q-1}{Q})a}}{\big(\frac{\alpha}{\alpha_Q}\big)^{(\frac{Q-1}{Q})b}}\Bigg)^{\frac{Q}{b}(1-\frac{\beta}{Q})}\text{SCMT}(\alpha,\beta).
\end{align*}
Passing to the limit $k\to\infty$, we get
\begin{align*}
  \text{CMT}_{a,b}(\alpha_Q,\beta)\leq\sup_{\alpha\in(0,\alpha_Q)}\Bigg(\cfrac{1-\big(\frac{\alpha}{\alpha_Q}\big)^{(\frac{Q-1}{Q})a}}{\big(\frac{\alpha}{\alpha_Q}\big)^{(\frac{Q-1}{Q})b}}\Bigg)^{\frac{Q}{b}(1-\frac{\beta}{Q})}\text{SCMT}(\alpha,\beta).
\end{align*}
The validity of \eqref{eq1.4} is thus immediate from the preceding estimates. To complete the proof, it is enough to verify that $\textnormal{CMT}_{a,b}(\alpha_Q,\beta)<+\infty$ implies $b\leq Q$. Indeed, if not, let there exists some $b$ with $b>Q$ such that $\textnormal{CMT}_{a,b}(\alpha_Q,\beta)<+\infty$. It follows from \eqref{eq1.4} that
\begin{align*}
    \limsup_{\alpha\to\alpha_Q^-}\Bigg(\cfrac{1-\big(\frac{\alpha}{\alpha_Q}\big)^{(\frac{Q-1}{Q})a}}{\big(\frac{\alpha}{\alpha_Q}\big)^{(\frac{Q-1}{Q})b}}\Bigg)^{\frac{Q}{b}(1-\frac{\beta}{Q})}\textnormal{SCMT}(\alpha,\beta)<+\infty.
\end{align*}
In view of Theorem \ref{thm1.3}, we also have
\begin{align*}
    \limsup_{\alpha\to\alpha_Q^-}\Bigg(\cfrac{1-\big(\frac{\alpha}{\alpha_Q}\big)^{(\frac{Q-1}{Q})a}}{\big(\frac{\alpha}{\alpha_Q}\big)^{(\frac{Q-1}{Q})b}}\Bigg)^{1-\frac{\beta}{Q}}\textnormal{SCMT}(\alpha,\beta)>0.
\end{align*}
This shows that
$$ \liminf_{\alpha\to\alpha_Q^-}\Bigg(\cfrac{1-\big(\frac{\alpha}{\alpha_Q}\big)^{(\frac{Q-1}{Q})a}}{\big(\frac{\alpha}{\alpha_Q}\big)^{(\frac{Q-1}{Q})b}}\Bigg)^{\frac{Q}{b}(1-\frac{\beta}{Q})}\Bigg(\cfrac{\big(\frac{\alpha}{\alpha_Q}\big)^{(\frac{Q-1}{Q})b}}{1-\big(\frac{\alpha}{\alpha_Q}\big)^{(\frac{Q-1}{Q})a}}\Bigg)^{1-\frac{\beta}{Q}}<+\infty, $$
which is contrary to the assumption $b>Q$. This finishes the proof. \qed

Next, we establish the continuity of $\textnormal{SCMT}(\cdot,\beta)$, which will be essential for the results that follow.
\begin{lemma}\label{lem6.2}
For all $0\leq\beta<Q$, the mapping $\alpha\mapsto \textnormal{SCMT}(\alpha,\beta)$ is continuous on $(0,\alpha_Q)$.
\end{lemma}
\begin{proof}
Let $0\leq\beta<Q$ and $\alpha\in(0,\alpha_Q)$ be fixed. For each $u\in HW^{1,Q}(\mathbb{H}^n)$, we define the map $J_u:(0,\alpha_Q)\to \mathbb{R}$ by
\begin{align*}
J_u(\alpha):=
\int_{\mathbb{H}^n}
\frac{\Phi\big(\alpha\big(1-\frac{\beta}{Q}\big)
|u|^{Q^\prime}\big)}{r(\xi)^\beta}\,\mathrm{d}\xi\quad\text{for all}\quad \alpha\in (0,\alpha_Q). \end{align*}
It follows from Lemma \ref{lem2.2} that for each
$u\in HW^{1,Q}(\mathbb{H}^n)$, the mapping $\alpha\mapsto J_u(\alpha)$ is convex on
$(0,\alpha_Q)$, thanks to the fact that the argument of the convex function $\Phi(\cdot)$ is linear in $\alpha$.
In view of Lemma \ref{lem4.1}, one sees that
\begin{align*}
\textnormal{SCMT}(\alpha,\beta)
=\sup_{u\in HW^{1,Q}(\mathbb{H}^n),\,\|D_H u\|_{Q}\leq 1,\,\|u\|_Q= 1} J_u(\alpha).
\end{align*}
Now, by applying Rockafellar \cite[{\S}5, Theorem 5.5, p. 35]{Rockafellar-1970}, we deduce that $\textnormal{SCMT}(\cdot,\beta)$ is convex on
$(0,\alpha_Q)$, being the pointwise supremum of convex functions. In addition, Theorem \ref{thm1.3} ensures that
$\textnormal{SCMT}(\alpha,\beta)<+\infty$ for every
$\alpha\in(0,\alpha_Q)$. Hence,
$\textnormal{SCMT}(\cdot,\beta)$ is a finite-valued convex
function on $(0,\alpha_Q)$. Thus, the continuity of the map $\alpha\mapsto \textnormal{SCMT}(\alpha,\beta)$ on $(0,\alpha_Q)$ follows directly from Rockafellar \cite[{\S}10, Theorem 10.1, p. 82]{Rockafellar-1970}. This finishes the proof.
\end{proof}
The following limit characterizes the behavior of $\textnormal{SCMT}(\cdot,\beta)$ near zero.
\begin{lemma}\label{lem6.3}
If $0<\beta<Q$, $a>0$, $b>0$ and $\alpha\in(0,\alpha_Q]$, then there holds
\begin{align*}
    \lim_{\zeta\to 0^+}\Bigg(\cfrac{1-\big(\frac{\zeta}{\alpha}\big)^{(\frac{Q-1}{Q})a}}{\big(\frac{\zeta}{\alpha}\big)^{(\frac{Q-1}{Q})b}}\Bigg)^{\frac{Q}{b}(1-\frac{\beta}{Q})}\textnormal{SCMT}(\zeta,\beta) =0.
\end{align*}
\end{lemma}
\begin{proof}
Let $\{\zeta_k\}_k$ be a real sequence such that $\zeta_k\to 0$ as $k\to\infty$. For each $k\in\mathbb{N}$, let $u_k\in HW^{1,Q}(\mathbb{H}^n)\setminus\{0\}$ be a maximizer for $\textnormal{SCMT}(\zeta_k,\beta)$. In view of Lemma \ref{lem4.1}, it follows that $\|u_k\|_Q=1$, $\|D_H u_k\|_Q\leq1$ and
\begin{align}\label{eq6.1}
    \int_{\mathbb{H}^n}\frac{\Phi\big(\zeta_k \big(1-\frac{\beta}{Q}\big)|u_k|^{Q^\prime}\big)}{r(\xi)^\beta}\,\mathrm{d}\xi=\textnormal{SCMT}(\zeta_k,\beta).
\end{align}
Thus, the sequence $\{u_k\}_k$ is bounded in $HW^{1,Q}(\mathbb{H}^n)$. Applying H\"{o}lder's inequality, Lemma \ref{lem2.2}, Theorem \ref{thm1.3}, and Lemma \ref{lem4.2}, together with the boundedness of $\{u_k\}_k$ in $HW^{1,Q}(\mathbb{H}^n)$, we infer, for all sufficiently large $k$, that
\begin{equation}\label{eq6.2}
\begin{aligned}
\int_{\mathbb{H}^n}&\frac{\operatorname{\exp}\big(\zeta_k \big(1-\frac{\beta}{Q}\big)|u_k|^{Q^\prime}\big)|u_k|^Q}{r(\xi)^\beta}\,\mathrm{d}\xi \\&=\int_{\mathbb{H}^n}\frac{\Phi\big(\zeta_k \big(1-\frac{\beta}{Q}\big)|u_k|^{Q^\prime}\big)|u_k|^Q}{r(\xi)^\beta}\,\mathrm{d}\xi +\sum_{j=0}^{Q-2}\frac{\big(\zeta_k \big(1-\frac{\beta}{Q}\big)\big)^{j}}{j!}\int_{\mathbb{H}^n}\frac{|u_k|^{Q+Q^\prime j}}{r(\xi)^\beta}\,\mathrm{d}\xi  \\
&\lesssim\bigg(\int_{\mathbb{H}^n}\frac{\Phi\big(2\zeta_k \big(1-\frac{\beta}{Q}\big)|u_k|^{Q^\prime}\big)}{r(\xi)^\beta}\,\mathrm{d}\xi\bigg)^{\frac{1}{2}}\bigg(\int_{\mathbb{H}^n}\frac{|u_k|^{2Q}}{r(\xi)^\beta}\,\mathrm{d}\xi\bigg)^{\frac{1}{2}} +\sum_{j=0}^{Q-2}\int_{\mathbb{H}^n}\frac{|u_k|^{Q+Q^\prime j}}{r(\xi)^\beta}\,\mathrm{d}\xi \\
&\lesssim 1. 
\end{aligned}
\end{equation}
The estimates \eqref{eq6.1} and \eqref{eq6.2} further yield that
\begin{align*}
\Bigg(&\cfrac{1-\big(\frac{\zeta_k}{\alpha}\big)^{(\frac{Q-1}{Q})a}}{\big(\frac{\zeta_k}{\alpha}\big)^{(\frac{Q-1}{Q})b}}\Bigg)^{\frac{Q}{b}(1-\frac{\beta}{Q})}\textnormal{SCMT}(\zeta_k,\beta)\\&= \Bigg(\cfrac{1-\big(\frac{\zeta_k}{\alpha}\big)^{(\frac{Q-1}{Q})a}}{\big(\frac{\zeta_k}{\alpha}\big)^{(\frac{Q-1}{Q})b}}\Bigg)^{\frac{Q}{b}(1-\frac{\beta}{Q})} \int_{\mathbb{H}^n}\frac{\Phi\big(\zeta_k \big(1-\frac{\beta}{Q}\big)|u_k|^{Q^\prime}\big)}{r(\xi)^\beta}\,\mathrm{d}\xi\\
&\leq \Bigg(\cfrac{1-\big(\frac{\zeta_k}{\alpha}\big)^{(\frac{Q-1}{Q})a}}{\big(\frac{\zeta_k}{\alpha}\big)^{(\frac{Q-1}{Q})b}}\Bigg)^{\frac{Q}{b}(1-\frac{\beta}{Q})} \frac{\big(\zeta_k \big(1-\frac{\beta}{Q}\big)\big)^{Q-1}}{(Q-1)!}\int_{\mathbb{H}^n}\frac{\operatorname{\exp}\big(\zeta_k \big(1-\frac{\beta}{Q}\big)|u_k|^{Q^\prime}\big)|u_k|^Q}{r(\xi)^\beta}\,\mathrm{d}\xi\\
&\lesssim \Bigg(\cfrac{1-\big(\frac{\zeta_k}{\alpha}\big)^{(\frac{Q-1}{Q})a}}{\big(\frac{\zeta_k}{\alpha}\big)^{(\frac{Q-1}{Q})b}}\Bigg)^{\frac{Q}{b}(1-\frac{\beta}{Q})}\zeta_k^{Q-1}\to 0\quad\text{as}\qquad k\to\infty.
\end{align*}
Hence, the desired assertion follows, and the proof is complete.
\end{proof}
To complement the preceding estimate, the following result describes the endpoint behavior of $\textnormal{SCMT}(\cdot,\beta)$ for $\alpha$ in the admissible range.
\begin{lemma}\label{lem6.4}
For $0<\beta<Q$ and $a>0$, if $(0<\alpha< \alpha_Q,\,b> 0)$ or $(\alpha= \alpha_Q,\, 0<b<Q)$, then there holds
\begin{align*}
    \lim_{\zeta\to \alpha^-}\Bigg(\cfrac{1-\big(\frac{\zeta}{\alpha}\big)^{(\frac{Q-1}{Q})a}}{\big(\frac{\zeta}{\alpha}\big)^{(\frac{Q-1}{Q})b}}\Bigg)^{\frac{Q}{b}(1-\frac{\beta}{Q})}\textnormal{SCMT}(\zeta,\beta) =0.
\end{align*}
\end{lemma}
\begin{proof}
Assume that the stated hypotheses are satisfied. If $\alpha<\alpha_Q$, then $\textnormal{SCMT}(\zeta,\beta)\leq\textnormal{SCMT}(\alpha,\beta)<+\infty$ for all $\zeta\in(0,\alpha)$, and the conclusion is immediate. Thus, it suffices to prove that
 \begin{align*}
    \lim_{\zeta\to \alpha^-_Q}\Bigg(\cfrac{1-\big(\frac{\zeta}{\alpha_Q}\big)^{(\frac{Q-1}{Q})a}}{\big(\frac{\zeta}{\alpha_Q}\big)^{(\frac{Q-1}{Q})b}}\Bigg)^{\frac{Q}{b}(1-\frac{\beta}{Q})}\textnormal{SCMT}(\zeta,\beta) =0.
\end{align*}
The upper bound in Theorem \ref{thm1.3} gives
\begin{align*}
     \limsup_{\zeta\to \alpha^-_Q}\Bigg(\cfrac{1-\big(\frac{\zeta}{\alpha_Q}\big)^{(\frac{Q-1}{Q})a}}{\big(\frac{\zeta}{\alpha_Q}\big)^{(\frac{Q-1}{Q})Q}}\Bigg)^{1-\frac{\beta}{Q}}\textnormal{SCMT}(\zeta,\beta)<+\infty.
\end{align*}
Now, by using the above estimate, it follows that
\begin{align*}
     &\limsup_{\zeta\to \alpha^-_Q}\Bigg(\cfrac{1-\big(\frac{\zeta}{\alpha_Q}\big)^{(\frac{Q-1}{Q})a}}{\big(\frac{\zeta}{\alpha_Q}\big)^{(\frac{Q-1}{Q})b}}\Bigg)^{\frac{Q}{b}(1-\frac{\beta}{Q})}\textnormal{SCMT}(\zeta,\beta)\\
     &=\limsup_{\zeta\to \alpha^-_Q}\Bigg[\Bigg(\cfrac{1-\big(\frac{\zeta}{\alpha_Q}\big)^{(\frac{Q-1}{Q})a}}{\big(\frac{\zeta}{\alpha_Q}\big)^{(\frac{Q-1}{Q})b}}\Bigg)^{\frac{Q}{b}(1-\frac{\beta}{Q})}\Bigg(\cfrac{\big(\frac{\zeta}{\alpha_Q}\big)^{(\frac{Q-1}{Q})Q}}{1-\big(\frac{\zeta}{\alpha_Q}\big)^{(\frac{Q-1}{Q})a}}\Bigg)^{1-\frac{\beta}{Q}}\\
     &\qquad\qquad\quad\times\Bigg(\cfrac{1-\big(\frac{\zeta}{\alpha_Q}\big)^{(\frac{Q-1}{Q})a}}{\big(\frac{\zeta}{\alpha_Q}\big)^{(\frac{Q-1}{Q})Q}}\Bigg)^{1-\frac{\beta}{Q}}\textnormal{SCMT}(\zeta,\beta)\Bigg]\\
     &\leq\lim_{\zeta\to \alpha^-_Q}\bigg(1-\Big(\frac{\zeta}{\alpha_Q}\Big)^{(\frac{Q-1}{Q})a}\bigg)^{(\frac{Q}{b}-1)(1-\frac{\beta}{Q})} \limsup_{\zeta\to \alpha^-_Q}\Bigg(\cfrac{1-\big(\frac{\zeta}{\alpha_Q}\big)^{(\frac{Q-1}{Q})a}}{\big(\frac{\zeta}{\alpha_Q}\big)^{(\frac{Q-1}{Q})Q}}\Bigg)^{1-\frac{\beta}{Q}}\textnormal{SCMT}(\zeta,\beta)\\
     &=0.
\end{align*}
Since the quantity in the statement is nonnegative, this finishes the proof.
\end{proof}
Finally, this section concludes with a proof of the existence of maximizers for the weighted Moser--Trudinger inequality under inhomogeneous norm constraints.
\subsection{{Proof of Theorem \ref{thm1.7}}}
Let the given hypotheses be satisfied. Due to Lemmas \ref{lem6.3} and \ref{lem6.4}, we deduce that
\begin{align*}
    \lim_{\zeta\to 0^+}\Bigg(\cfrac{1-\big(\frac{\zeta}{\alpha}\big)^{(\frac{Q-1}{Q})a}}{\big(\frac{\zeta}{\alpha}\big)^{(\frac{Q-1}{Q})b}}\Bigg)^{\frac{Q}{b}(1-\frac{\beta}{Q})}\textnormal{SCMT}(\zeta,\beta)<\textnormal{CMT}_{a,b}(\alpha,\beta)
\end{align*}
and
\begin{align*}
    \lim_{\zeta\to \alpha^-}\Bigg(\cfrac{1-\big(\frac{\zeta}{\alpha}\big)^{(\frac{Q-1}{Q})a}}{\big(\frac{\zeta}{\alpha}\big)^{(\frac{Q-1}{Q})b}}\Bigg)^{\frac{Q}{b}(1-\frac{\beta}{Q})}\textnormal{SCMT}(\zeta,\beta) <\textnormal{CMT}_{a,b}(\alpha,\beta).
\end{align*}
Moreover, since \eqref{eq1.4} holds, Lemma \ref{lem6.2} implies that there exists some $\zeta\in(0,\alpha)$ such that
\begin{align*}
   \Bigg(\cfrac{1-\big(\frac{\zeta}{\alpha}\big)^{(\frac{Q-1}{Q})a}}{\big(\frac{\zeta}{\alpha}\big)^{(\frac{Q-1}{Q})b}}\Bigg)^{\frac{Q}{b}(1-\frac{\beta}{Q})}\textnormal{SCMT}(\zeta,\beta)=\textnormal{CMT}_{a,b}(\alpha,\beta).
\end{align*}
By Theorem \ref{thm1.4}, there exists a maximizer $u\in HW^{1,Q}(\mathbb{H}^n)\setminus\{0\}$ for $\textnormal{SCMT}(\zeta,\beta)$. It follows from Lemma \ref{lem4.1} that $\|u\|_Q=1$, $\|D_H u\|_Q\leq1$ and \begin{align*}
    \int_{\mathbb{H}^n}\frac{\Phi\big(\zeta \big(1-\frac{\beta}{Q}\big)|u|^{Q^\prime}\big)}{r(\xi)^\beta}\,\mathrm{d}\xi=\textnormal{SCMT}(\zeta,\beta).
\end{align*}
Define a function $w\colon \mathbb{H}^n\to\mathbb{R}$ by
\begin{align*}
    w(\xi)=\Big(\frac{\zeta}{\alpha}\Big)^{\frac{Q-1}{Q}} u(\delta_\theta(\xi))\qquad\text{with}\qquad \theta= \Bigg(\cfrac{\big(\frac{\zeta}{\alpha}\big)^{(\frac{Q-1}{Q})b}}{1-\big(\frac{\zeta}{\alpha}\big)^{(\frac{Q-1}{Q})a}}\Bigg)^{\frac{1}{b}}.
\end{align*}
By straightforward calculations, we obtain
\begin{align*}
    \|D_H w\|_Q^a=\Big(\frac{\zeta}{\alpha}\Big)^{(\frac{Q-1}{Q})a}\|D_H u\|_Q^a\leq \Big(\frac{\zeta}{\alpha}\Big)^{(\frac{Q-1}{Q})a}\quad\text{and}\quad \|w\|_Q^b=1- \bigg(\frac{\zeta}{\alpha}\bigg)^{(\frac{Q-1}{Q})a}.
\end{align*}
This implies that $\|D_H w\|_Q^a+\|w\|_Q^b\leq 1$. Now, by gathering all of the above information, we get
\begin{align*}
 \textnormal{CMT}_{a,b}(\alpha,\beta)&=\Bigg(\cfrac{1-\big(\frac{\zeta}{\alpha}\big)^{(\frac{Q-1}{Q})a}}{\big(\frac{\zeta}{\alpha}\big)^{(\frac{Q-1}{Q})b}}\Bigg)^{\frac{Q}{b}(1-\frac{\beta}{Q})}\textnormal{SCMT}(\zeta,\beta)\\
 &=\Bigg(\cfrac{1-\big(\frac{\zeta}{\alpha}\big)^{(\frac{Q-1}{Q})a}}{\big(\frac{\zeta}{\alpha}\big)^{(\frac{Q-1}{Q})b}}\Bigg)^{\frac{Q}{b}(1-\frac{\beta}{Q})}\int_{\mathbb{H}^n}\frac{\Phi\big(\zeta \big(1-\frac{\beta}{Q}\big)|u|^{Q^\prime}\big)}{r(\xi)^\beta}\,\mathrm{d}\xi\\
 &=\Bigg(\cfrac{1-\big(\frac{\zeta}{\alpha}\big)^{(\frac{Q-1}{Q})a}}{\big(\frac{\zeta}{\alpha}\big)^{(\frac{Q-1}{Q})b}}\Bigg)^{\frac{Q}{b}(1-\frac{\beta}{Q})}\int_{\mathbb{H}^n}\frac{\Phi\big(\zeta \big(1-\frac{\beta}{Q}\big)|u(\delta_\theta(\xi)|^{Q^\prime}\big)}{r(\delta_\theta(\xi))^\beta}\,\mathrm{d}(\delta_\theta(\xi))\\
  &=\Bigg(\cfrac{1-\big(\frac{\zeta}{\alpha}\big)^{(\frac{Q-1}{Q})a}}{\big(\frac{\zeta}{\alpha}\big)^{(\frac{Q-1}{Q})b}}\Bigg)^{\frac{Q}{b}(1-\frac{\beta}{Q})}\theta ^{Q-\beta}\int_{\mathbb{H}^n}\frac{\Phi\big(\alpha \big(1-\frac{\beta}{Q}\big)|w|^{Q^\prime}\big)}{r(\xi)^\beta}\,\mathrm{d}\xi\\
 &=\int_{\mathbb{H}^n}\frac{\Phi\big(\alpha \big(1-\frac{\beta}{Q}\big)|w|^{Q^\prime}\big)}{r(\xi)^\beta}\,\mathrm{d}\xi.
\end{align*}
Hence, we conclude that $w\in HW^{1,Q}(\mathbb{H}^n)\setminus\{0\}$ is actually a maximizer for $\textnormal{CMT}_{a,b}(\alpha,\beta)$. This completes the proof. \qed

\section{The Improved Subcritical Weighted Moser--Trudinger Inequality and Its Maximizer}\label{sec7}
This section establishes a strengthened subcritical weighted Moser--Trudinger inequality on the Heisenberg group in Theorem \ref{thm1.8}, together with the existence of corresponding maximizers as stated in Theorem \ref{thm1.9}, extending the result of Ishiwata et al. \cite{Ishiwata-Nakamura-Wadade-2014} for the subcritical weighted Moser--Trudinger inequality to the present setting. Now, we begin with the following preliminary result, which plays a key role in proving Theorem \ref{thm1.8}.
\begin{lemma}\label{lem7.1}
For all $0\leq\beta<Q$, $p\geq Q$ and $u\in HW^{1,Q}(\mathbb{H}^n)$ satisfying $\|D_H u\|_Q\leq1$, there holds
 \begin{align*}
     \int_{\mathbb{H}^n}\frac{|u|^p}{r(\xi)^\beta}\,\mathrm{d}\xi\leq C\bigg(\int_{\mathbb{H}^n}|u|^Q\,\mathrm{d}\xi\bigg)^{1-\frac{\beta}{Q}},
 \end{align*}
where $C=C(Q,p,\beta)>0$ is a constant depending only on $Q$, $p$ and $\beta$.
\end{lemma}
\begin{proof}
Thanks to the density of $C_0^\infty(\mathbb{H}^n)$ in $HW^{1,Q}(\mathbb{H}^n)$, the desired result need only be verified for   $u\in C_0^\infty(\mathbb{H}^n)\setminus\{0\}$, $u\geq0$ and $\|D_Hu\|_Q\leq1$. Define $S=\{\xi\in \mathbb{H}^n\colon \,u(\xi)>1\}$. Notice that $S\subset \mathbb{H}^n $ is a bounded  set and there holds
\begin{align*}
        |S|=\int_S\,\mathrm{d}\xi\leq \int_S|u|^Q\,\mathrm{d}\xi\leq\int_{\mathbb{H}^n}|u|^Q\,\mathrm{d}\xi.
\end{align*}
Divide the whole integral as follows
$$  \int_{\mathbb{H}^n}\frac{|u|^p}{r(\xi)^\beta}\,\mathrm{d}\xi= \int_{\mathbb{H}^n\setminus S}\frac{|u|^p}{r(\xi)^\beta}\,\mathrm{d}\xi+ \int_{S}\frac{|u|^p}{r(\xi)^\beta}\,\mathrm{d}\xi=:I+J.$$
It is standard to see that $u(\xi)\leq 1$ for all $\xi\in\mathbb{H}^n\setminus S$ and thus, we deduce that
\begin{align*}
    I&\leq \int_{\{\xi\in \mathbb{H}^n,\, u(\xi)\leq 1,\,r(\xi)\geq\|u\|_Q\}}\frac{|u|^Q}{r(\xi)^\beta}\,\mathrm{d}\xi+\int_{\{\xi\in \mathbb{H}^n,\, u(\xi)\leq 1,\,r(\xi)\leq\|u\|_Q\}}\frac{|u|^Q}{r(\xi)^\beta}\,\mathrm{d}\xi\\
    &\leq\frac{1}{\|u\|_Q^\beta}\int_{\mathbb{H}^n}|u|^Q\,\mathrm{d}\xi+\int_\Sigma\int_{0}^{\|u\|_Q}\eta^{Q-\beta-1}\,\mathrm{d}\eta\,\mathrm{d}\mu(\xi^\ast)\\
    &=\|u\|_Q^{Q-\beta}+\frac{\int_\Sigma\,\mathrm{d}\mu(\xi^\ast)}{Q-\beta} \|u\|_Q^{Q-\beta}\leq C(Q,\beta)\|u\|_Q^{Q-\beta},
\end{align*}
where $C(Q,\beta)>0$ is a constant depending only on $Q$ and $\beta$. On the other hand, let $v\colon\mathbb{H}^n\to\mathbb{R}$ be such that $v(\xi)=u(\xi)-1$ for all $\xi\in S$. It follows that $\|D_H v\|_{L^Q(S)}\leq1$ and $v\in HW^{1,Q}_0(S)$. In addition,  we obtain for all $p\geq Q$ that
\begin{align*}
  |u|^p\leq 2^{p-1}\big(1+C(Q,p,\beta)\big)\operatorname{\exp}\Big(\alpha_{Q}\Big(1-\frac{\beta}{Q}\Big)|v|^{Q^\prime}\Big)\quad\text{for all}\quad \xi\in S,
\end{align*}
where $C(Q,p,\beta)>0$ is a constant depending on $Q$, $p$ and $\beta$. Due to the weighted Moser--Trudinger inequality established by Lam--Lu--Tang \cite[Lemma 4.1]{Lam-Lu-Tang-2012} on bounded subsets of the Heisenberg group and using the previous inequality, we get
\begin{align*}
 J &  \leq 2^{p-1}\big(1+C(Q,p,\beta)\big)\int_{S}\frac{\operatorname{\exp}\big(\alpha_Q\big(1-\frac{\beta}{Q}\big) |v|^{Q^\prime}\big)}{r(\xi)^\beta}\,\mathrm{d}\xi\\
 &\leq 2^{p-1}\big(1+C(Q,p,\beta)\big)C(Q,\beta)|S|^{1-\frac{\beta}{Q}}\leq C(Q,p,\beta)\|u\|_Q^{Q-\beta},
\end{align*}
where $C=C(Q,p,\beta)>0$ is a constant depending only on $Q$, $p$ and $\beta$. By gathering all the above information, we conclude that the desired result holds. This completes the proof.
\end{proof}
Now, we are in a position to prove Theorem \ref{thm1.8} as follows.
\subsection{{Proof of Theorem \ref{thm1.8}}}
Let $0<\alpha<\alpha_{Q}$ and $0\leq \beta<Q$ be given. In addition, assume that $u\in HW^{1,Q}(\mathbb{H}^n)$ with $\|D_H u\|_Q\leq1$. Now, divide the whole integral in the following manner
\begin{align*}
\int_{\mathbb{H}^n}\frac{\operatorname{\exp}\big(\alpha \big(1-\frac{\beta}{Q}\big)|u|^{Q^\prime}\big)|u|^Q}{r(\xi)^\beta}\,\mathrm{d}\xi=:I+II,
\end{align*}
where
\begin{align*}
    I=\int_{\mathbb{H}^n}\frac{\Phi\big(\alpha \big(1-\frac{\beta}{Q}\big)|u|^{Q^\prime}\big)|u|^Q}{r(\xi)^\beta}\mathrm{d}\xi\quad\text{and}\quad II=\int_{\mathbb{H}^n}\sum_{j=0}^{Q-2}\frac{\big(\alpha \big(1-\frac{\beta}{Q}\big)\big)^j}{j!}\bigg(\frac{|u|^{Q^\prime(Q-1+j)}}{r(\xi)^\beta}\bigg)\mathrm{d}\xi.
\end{align*}
Take $t>1$ such that $\frac{1}{t}+\frac{1}{t^\prime}=1$. Let $t>1$ be close to $1$ with $t\alpha<\alpha_Q$, which is possible since $\alpha<\alpha_Q$. An application of H\"{o}lder's inequality, together with Theorem \ref{thm1.3} and Lemmas \ref{lem2.2} and \ref{lem7.1}, yields
\begin{equation}\label{eq7.1}
\begin{aligned}
I&\leq\bigg(\int_{\mathbb{H}^n}\frac{\Phi\big(t\alpha \big(1-\frac{\beta}{Q}\big)|u|^{Q^\prime}\big)}{r(\xi)^\beta}\,\mathrm{d}\xi\bigg)^{\frac{1}{t}}\bigg(\int_{\mathbb{H}^n}\frac{|u|^{Qt^\prime}}{r(\xi)^\beta}\,\mathrm{d}\xi\bigg)^{\frac{1}{t^\prime}}  \\
&\lesssim \bigg(\int_{\mathbb{H}^n}|u|^Q\,\mathrm{d}\xi\bigg)^{\frac{1}{t}(1-\frac{\beta}{Q})}\bigg(\int_{\mathbb{H}^n}|u|^Q\,\mathrm{d}\xi\bigg)^{\frac{1}{t^\prime}(1-\frac{\beta}{Q})} =\bigg(\int_{\mathbb{H}^n}|u|^Q\,\mathrm{d}\xi\bigg)^{1-\frac{\beta}{Q}}. 
\end{aligned}
\end{equation}
Invoking Lemma \ref{lem7.1} once more, we obtain
\begin{equation}\label{eq7.2}
\begin{aligned}
        II&=\sum_{j=0}^{Q-2}\frac{\big(\alpha \big(1-\frac{\beta}{Q}\big)\big)^j}{j!}\int_{\mathbb{H}^n}\frac{|u|^{Q^\prime(Q-1+j)}}{r(\xi)^\beta}\,\mathrm{d}\xi  \\
        &\lesssim\sum_{j=0}^{Q-2}\frac{\big(\alpha \big(1-\frac{\beta}{Q}\big)\big)^j}{j!}\bigg(\int_{\mathbb{H}^n}|u|^Q\,\mathrm{d}\xi\bigg)^{1-\frac{\beta}{Q}} \lesssim \bigg(\int_{\mathbb{H}^n}|u|^Q\,\mathrm{d}\xi\bigg)^{1-\frac{\beta}{Q}}. 
    \end{aligned}
\end{equation}
The desired result now follows directly from \eqref{eq7.1} and \eqref{eq7.2}. To finish the proof, we only need to show that $\textnormal{SCMT}_{\operatorname{e}}(\alpha_Q,\beta)=+\infty$. Recall the Moser's function sequence $\{w_k\}_k$ that is defined in \eqref{eq4.6}. Note that \eqref{eq4.8} remains valid. This yields the following estimate
\begin{align*}
\int_{\mathbb{H}^n}&\frac{\operatorname{\exp}\big(\alpha_Q \big(1-\frac{\beta}{Q}\big)|w_k|^{Q^\prime}\big)|w_k|^Q}{r(\xi)^\beta}\,\mathrm{d}\xi\\&\geq\int_{\{\xi\in\mathbb{H}^n\colon\, r(\xi)< sk^{-1}\}}\frac{\operatorname{\exp}\big(\alpha_Q \big(1-\frac{\beta}{Q}\big)|w_k|^{Q^\prime}\big)|w_k|^Q}{r(\xi)^\beta}\,\mathrm{d}\xi\\
&=\int_\Sigma\int_{0}^{sk^{-1}}\operatorname{\exp}\bigg(\alpha_Q \bigg(1-\frac{\beta}{Q}\bigg)\bigg|c^{-\frac{1}{Q}}_Q(\ln k)^{\frac{Q-1}{Q}}\bigg|^{Q^\prime}\bigg)\bigg|c^{-\frac{1}{Q}}_Q(\ln k)^{\frac{Q-1}{Q}}\bigg|^Q\eta^{Q-\beta-1}\,\mathrm{d}\eta\,\mathrm{d}\mu(\xi^\ast)\\
  &= \bigg(\int_\Sigma \,\mathrm{d}\mu(\xi^\ast)\bigg)\operatorname{\exp}\big((Q-\beta)\ln k\big)c^{-1}_Q(\ln k)^{Q-1}\int_{0}^{sk^{-1}}\eta^{Q-\beta-1}\,\mathrm{d}\eta\\
  &= \bigg(\int_\Sigma \,\mathrm{d}\mu(\xi^\ast)\bigg)\,\ \frac{\operatorname{\exp}\big((Q-\beta)\ln k\big)c^{-1}_Q(\ln k)^{Q-1}(sk^{-1})^{Q-\beta}}{Q-\beta}\gtrsim (\ln k)^{Q-1}.
\end{align*}
Hence, we obtain
\begin{align*}
    \textnormal{SCMT}_{\operatorname{e}}(\alpha_Q,\beta)&\geq\cfrac{1}{\|w_k\|_Q^{Q(1-\frac{\beta}{Q})}} \int_{\mathbb{H}^n}\frac{\operatorname{\exp}\big(\alpha_Q \big(1-\frac{\beta}{Q}\big)|w_k|^{Q^\prime}\big)|w_k|^Q}{r(\xi)^\beta}\,\mathrm{d}\xi\gtrsim (\ln k)^{\frac{Q^2-\beta}{Q}}\\
    &\to+\infty\quad\text{as}\quad k\to\infty.
\end{align*}
This establishes the sharpness of $\alpha_{Q}$, which completes the proof. \qed

The proof of the following result is analogous to those of Lemmas \ref{lem4.1} and \ref{lem4.5} and is thus omitted.
\begin{lemma}\label{lem7.2}
 For any $0<\alpha<\alpha_Q$ and $0\leq\beta<Q$, the following result holds
 \begin{align*}
  \textnormal{SCMT}_{\operatorname{e}}(\alpha,\beta)&=\sup_{u\in HW^{1,Q}(\mathbb{H}^n),\,\|D_H u\|_{Q}\leq 1,\,\|u\|_Q= 1} \int_{\mathbb{H}^n}\frac{\operatorname{\exp}\big(\alpha \big(1-\frac{\beta}{Q}\big)|u|^{Q^\prime}\big)|u|^Q}{r(\xi)^\beta}\,\mathrm{d}\xi\\
  &=\sup_{u\in HW^{1,Q}(\mathbb{H}^n),\,\|D_H u\|_{Q}= 1,\,\|u\|_Q= 1} \int_{\mathbb{H}^n}\frac{\operatorname{\exp}\big(\alpha \big(1-\frac{\beta}{Q}\big)|u|^{Q^\prime}\big)|u|^Q}{r(\xi)^\beta}\,\mathrm{d}\xi.
 \end{align*}
\end{lemma}
The following proposition establishes a stronger form of the subcritical weighted Moser--Trudinger inequality needed for the subsequent analysis.
\begin{proposition}\label{prop7.3}
 Let $0<\alpha<\alpha_Q$, $0\leq\beta<Q$ and $p\geq Q$ be real numbers. Then there holds
 \begin{align*}
    \Lambda_p(\alpha,\beta):&=\sup_{u\in HW^{1,Q}(\mathbb{H}^n),~\|D_H u\|_{Q}\leq 1}\cfrac{1}{\|u\|_p^{p(1-\frac{\beta}{Q})}} \int_{\mathbb{H}^n}\frac{\operatorname{\exp}\big(\alpha \big(1-\frac{\beta}{Q}\big)|u|^{Q^\prime}\big)|u|^p}{r(\xi)^\beta}\,\mathrm{d}\xi<+\infty.
 \end{align*}
Moreover, the constant $\alpha_Q$ is sharp, that is, the above supremum is infinite for all $\alpha\geq \alpha_Q$.
\end{proposition}
\begin{proof}
By using the density argument, we can assume $u\in C_0^\infty(\mathbb{H}^n)\setminus\{0\}$, $u\geq0$ and $\|D_Hu\|_Q\leq1$. Define the set $\Omega=\{\xi\in \mathbb{H}^n\colon \,u(\xi)>1\}$. It follows that $\Omega$ is a bounded set in $\mathbb{H}^n$, and there holds
\begin{align}\label{eq7.3}
        |\Omega|=\int_\Omega\,\mathrm{d}\xi\leq \int_\Omega|u|^p\,\mathrm{d}\xi\leq\int_{\mathbb{H}^n}|u|^p\,\mathrm{d}\xi.
\end{align}
Now, we split the whole integral as follows
\begin{align*}
 \int_{\mathbb{H}^n}\frac{\operatorname{\exp}\big(\alpha \big(1-\frac{\beta}{Q}\big)|u|^{Q^\prime}\big)|u|^p}{r(\xi)^\beta}\,\mathrm{d}\xi&=\left(\int_{\mathbb{H}^n\setminus\Omega}+\int_\Omega\right)\frac{\operatorname{\exp}\big(\alpha \big(1-\frac{\beta}{Q}\big)|u|^{Q^\prime}\big)|u|^p}{r(\xi)^\beta}\,\mathrm{d}\xi=:I+J.
\end{align*}
Since $u(\xi)\leq 1$ for all $\xi\in\mathbb{H}^n\setminus \Omega$ and $\alpha\in(0,\alpha_Q)$, therefore we have
\begin{align*}
    I&\leq \operatorname{\exp}\Big(\alpha \Big(1-\frac{\beta}{Q}\Big)\Big)\int_{\{\xi\in \mathbb{H}^n,\, u(\xi)\leq 1\}}\frac{|u|^p}{r(\xi)^\beta}\,\mathrm{d}\xi\\
    &\leq \operatorname{\exp}\Big(\alpha \Big(1-\frac{\beta}{Q}\Big)\Big)\Bigg(\int_{\big\{\xi\in \mathbb{H}^n,\, u(\xi)\leq 1,\,r(\xi)\geq\|u\|_p^{\frac{p}{Q}}\big\}}+\int_{\big\{\xi\in \mathbb{H}^n,\, u(\xi)\leq 1,\,r(\xi)\leq\|u\|_p^{\frac{p}{Q}}\big\}}\Bigg)\frac{|u|^p}{r(\xi)^\beta}\,\mathrm{d}\xi\\
    &\leq\operatorname{\exp}\Big(\alpha \Big(1-\frac{\beta}{Q}\Big)\Big)\bigg({\|u\|_p^{-\frac{p\beta}{Q}}}\int_{\mathbb{H}^n}|u|^p\,\mathrm{d}\xi+\int_\Sigma\int_{0}^{\|u\|_p^{\frac{p}{Q}}}\eta^{Q-\beta-1}\,\mathrm{d}\eta\,\mathrm{d}\mu(\xi^\ast)\bigg)\\
    &\leq\operatorname{\exp}\Big(\alpha_Q \Big(1-\frac{\beta}{Q}\Big)\Big)\bigg(\|u\|_p^{p(1-\frac{\beta}{Q})}+\frac{\int_\Sigma\,\mathrm{d}\mu(\xi^\ast)}{Q-\beta} \|u\|_p^{p(1-\frac{\beta}{Q})}\bigg) \leq C(Q,\beta)\|u\|_p^{p(1-\frac{\beta}{Q})},
\end{align*}
where $C(Q,\beta)>0$ is a constant depending only on $Q$ and $\beta$.  Next, let $w\colon\Omega\to\mathbb{R}$ be a function such that $w(\xi)=u(\xi)-1$ for all $\xi\in \Omega$. Thus, we conclude that $\|D_H w\|_{L^Q(\Omega)}\leq1$ and $w\in HW^{1,Q}_0(\Omega)$. Observe that for any $\vartheta\geq Q$, it holds that
\begin{align*}
  |u|^\vartheta\leq 2^{\vartheta-1}\big(1+C(Q,\vartheta,\beta)\big)\operatorname{\exp}\Big(\alpha_{Q}\Big(1-\frac{\beta}{Q}\Big)|w|^{Q^\prime}\Big)\quad\text{for all}\quad \xi\in \Omega,
\end{align*}
where $C(Q,\vartheta,\beta)>0$ is a constant depending only on $Q$, $\vartheta$ and $\beta$. Let $\zeta,\,\zeta^\prime>1$ be such that $\frac{1}{\zeta}+\frac{1}{\zeta^\prime}=1$. By choosing $\zeta>1$ close to $1$ and $\varepsilon>0$ close to $0$, we can assume that $(1+\varepsilon)\zeta\alpha<\alpha_Q$, which is possible since $\alpha<\alpha_Q$. Thus, by applying H\"{o}lder's inequality, \eqref{eq4.999}, \eqref{eq7.3} and the weighted Moser--Trudinger inequality established by Lam--Lu--Tang \cite[Lemma 4.1]{Lam-Lu-Tang-2012} on bounded subsets of the Heisenberg group, we deduce from the above information that
\begin{align*}
    J&\leq\bigg( \int_\Omega\frac{\operatorname{\exp}\big(\zeta\alpha \big(1-\frac{\beta}{Q}\big)|u|^{Q^\prime}\big)}{r(\xi)^\beta}\,\mathrm{d}\xi\bigg)^{\frac{1}{\zeta}} \bigg( \int_\Omega\frac{|u|^{p \zeta^\prime}}{r(\xi)^\beta}\,\mathrm{d}\xi\bigg)^{\frac{1}{\zeta^\prime}}\\
    &\lesssim \bigg( \int_\Omega\frac{\operatorname{\exp}\big(\zeta\alpha \big(1-\frac{\beta}{Q}\big)\big((1+\varepsilon)|w|^{Q^\prime}+C_\varepsilon\big)\big)}{r(\xi)^\beta}\mathrm{d}\xi\bigg)^{\frac{1}{\zeta}}\bigg( \int_\Omega\frac{\operatorname{\exp}\big(\alpha_Q \big(1-\frac{\beta}{Q}\big)|w|^{Q^\prime}\big)}{r(\xi)^\beta}\mathrm{d}\xi\bigg)^{\frac{1}{\zeta^\prime}}\\
    &\lesssim \bigg( \int_\Omega\frac{\operatorname{\exp}\big((1+\varepsilon)\zeta\alpha \big(1-\frac{\beta}{Q}\big)|w|^{Q^\prime}\big)}{r(\xi)^\beta}\,\mathrm{d}\xi\bigg)^{\frac{1}{\zeta}} \bigg( \int_\Omega\frac{\operatorname{\exp}\big(\alpha_Q \big(1-\frac{\beta}{Q}\big)|w|^{Q^\prime}\big)}{r(\xi)^\beta}\,\mathrm{d}\xi\bigg)^{\frac{1}{\zeta^\prime}}\\
    &\lesssim \bigg(\int_{\mathbb{H}^n}|u|^p\,\mathrm{d}\xi\bigg)^{{\frac{1}{\zeta}(1-\frac{\beta}{Q})}} \bigg(\int_{\mathbb{H}^n}|u|^p\,\mathrm{d}\xi\bigg)^{{\frac{1}{\zeta^\prime}(1-\frac{\beta}{Q})}}=\bigg(\int_{\mathbb{H}^n}|u|^p\,\mathrm{d}\xi\bigg)^{{1-\frac{\beta}{Q}}}.
\end{align*}
The desired inequality now follows by gathering the above estimates. To conclude the proof, it remains to establish the sharpness of $\alpha_Q$. In fact, it suffices to prove that $\Lambda_p(\alpha_Q,\beta)=+\infty$. For this purpose, let us reconsider Moser's function sequence $\{w_k\}_k$ introduced in \eqref{eq4.6}. In view of the validity of \eqref{eq4.8}, the following estimate follows
\begin{align*}
\int_{\mathbb{H}^n}&\frac{\operatorname{\exp}\big(\alpha_Q \big(1-\frac{\beta}{Q}\big)|w_k|^{Q^\prime}\big)|w_k|^p}{r(\xi)^\beta}\,\mathrm{d}\xi\\&\geq\int_{\{\xi\in\mathbb{H}^n\colon\, r(\xi)< sk^{-1}\}}\frac{\operatorname{\exp}\big(\alpha_Q \big(1-\frac{\beta}{Q}\big)|w_k|^{Q^\prime}\big)|w_k|^p}{r(\xi)^\beta}\,\mathrm{d}\xi\\
&=\int_\Sigma\int_{0}^{sk^{-1}}\operatorname{\exp}\bigg(\alpha_Q \bigg(1-\frac{\beta}{Q}\bigg)\bigg|c^{-\frac{1}{Q}}_Q(\ln k)^{\frac{Q-1}{Q}}\bigg|^{Q^\prime}\bigg)\bigg|c^{-\frac{1}{Q}}_Q(\ln k)^{\frac{Q-1}{Q}}\bigg|^p\eta^{Q-\beta-1}\,\mathrm{d}\eta\,\mathrm{d}\mu(\xi^\ast)\\
  &= \bigg(\int_\Sigma \,\mathrm{d}\mu(\xi^\ast)\bigg)\operatorname{\exp}\big((Q-\beta)\ln k\big)c^{-\frac{p}{Q}}_Q(\ln k)^{p(\frac{Q-1}{Q})}\int_{0}^{sk^{-1}}\eta^{Q-\beta-1}\,\mathrm{d}\eta\\
  &= \bigg(\int_\Sigma \,\mathrm{d}\mu(\xi^\ast)\bigg)\,\ \frac{\operatorname{\exp}\big((Q-\beta)\ln k\big)c^{-\frac{p}{Q}}_Q(\ln k)^{p(\frac{Q-1}{Q})}(sk^{-1})^{Q-\beta}}{Q-\beta}\gtrsim (\ln k)^{p(\frac{Q-1}{Q})}.
\end{align*}
This yields that
\begin{align*}
    \Lambda_p(\alpha_Q,\beta)&\geq\cfrac{1}{\|w_k\|_p^{p(1-\frac{\beta}{Q})}} \int_{\mathbb{H}^n}\frac{\operatorname{\exp}\big(\alpha_Q \big(1-\frac{\beta}{Q}\big)|w_k|^{Q^\prime}\big)|w_k|^p}{r(\xi)^\beta}\mathrm{d}\xi\gtrsim (\ln k)^{p(\frac{Q-1}{Q})+\frac{p}{Q}(1-\frac{\beta}{Q})}\\
    &\to+\infty\quad\text{as}\quad k\to\infty.
\end{align*}
This establishes the sharpness of $\alpha_Q$. Hence, the proof is complete.
\end{proof}
In what follows, the compactness result plays a crucial role in establishing maximizers for the subcritical weighted Moser--Trudinger inequality stated in Theorem \ref{thm1.9}.
\begin{proposition}\label{prop7.4}
Let $0<\alpha<\alpha_Q$, $0<\beta<Q$ and $p\geq Q$ be real numbers. Then for any sequence $\{u_k\}_k\subset HW^{1,Q}(\mathbb{H}^n)$ satisfying $u_k\rightharpoonup u$ in $HW^{1,Q}(\mathbb{H}^n)$ as $k\to\infty$ and $\|D_H u_k\|_Q=1$, there holds
\begin{align*}
\lim_{k\to\infty}\int_{\mathbb{H}^n}\frac{\operatorname{\exp}\big(\alpha \big(1-\frac{\beta}{Q}\big)|u_k|^{Q^\prime}\big)|u_k|^p}{r(\xi)^\beta}\,\mathrm{d}\xi=\int_{\mathbb{H}^n}\frac{\operatorname{\exp}\big(\alpha \big(1-\frac{\beta}{Q}\big)|u|^{Q^\prime}\big)|u|^p}{r(\xi)^\beta}\,\mathrm{d}\xi.
\end{align*}
\end{proposition}
\begin{proof} The proof follows that of Proposition \ref{prop4.4}, the only difference being the use of Proposition \ref{prop7.3} instead of Theorem \ref{thm1.3}. The details are left to the reader.
\end{proof}
To end this section, we now prove Theorem \ref{thm1.9} as follows.
\subsection{Proof of Theorem \ref{thm1.9}}
Let $\{u_k\}_k\subset HW^{1,Q}(\mathbb{H}^n)\setminus\{0\}$ be a maximizing sequence for $\textnormal{SCMT}_{\operatorname{e}}(\alpha,\beta)$. In virtue of Lemma \ref{lem7.2}, we have $\| u_k\|_Q=1$, $\| D_H u_k\|_Q=1$ and
\begin{align*}
 \lim_{k\to\infty} \int_{\mathbb{H}^n}\frac{\operatorname{\exp}\big(\alpha \big(1-\frac{\beta}{Q}\big)|u_k|^{Q^\prime}\big)|u_k|^Q}{r(\xi)^\beta}\,\mathrm{d}\xi=\textnormal{SCMT}_{\operatorname{e}}(\alpha,\beta).
\end{align*}
It follows that $\{u_k\}_k$ is a bounded sequence in $ HW^{1,Q}(\mathbb{H}^n)$. Hence, passing to a subsequence if necessary, we may assume that $u_k\rightharpoonup u$ in $ HW^{1,Q}(\mathbb{H}^n)$ as $k\to\infty$. Observe that \eqref{eq4.10} and \eqref{eq4.11} are still true. Due to \eqref{eq4.10} and Proposition \ref{prop7.4}, we have
\begin{align}\label{eq7.4}
\textnormal{SCMT}_{\operatorname{e}}(\alpha,\beta) =\int_{\mathbb{H}^n}\frac{\operatorname{\exp}\big(\alpha \big(1-\frac{\beta}{Q}\big)|u|^{Q^\prime}\big)|u|^Q}{r(\xi)^\beta}\,\mathrm{d}\xi.
\end{align}
Next, we claim that $u\neq0$. In fact, if we assume $u=0$, then \eqref{eq7.4} gives $\textnormal{SCMT}_{\operatorname{e}}(\alpha,\beta)=0$, contradicting $\textnormal{SCMT}_{\operatorname{e}}(\alpha,\beta)>0$. It follows from the definition of $\textnormal{SCMT}_{\operatorname{e}}(\alpha,\beta)$,  \eqref{eq4.11} and \eqref{eq7.4} that
\begin{align*}
\textnormal{SCMT}_{\operatorname{e}}(\alpha,\beta) &\geq \frac{1}{\|u\|_Q^{Q-\beta}} \int_{\mathbb{H}^n}\frac{\operatorname{\exp}\big(\alpha \big(1-\frac{\beta}{Q}\big)|u|^{Q^\prime}\big)|u|^Q}{r(\xi)^\beta}\,\mathrm{d}\xi= \frac{\textnormal{SCMT}_{\operatorname{e}}(\alpha,\beta)}{\|u\|_Q^{Q-\beta}}\\
&\geq  \textnormal{SCMT}_{\operatorname{e}}(\alpha,\beta).
\end{align*}
This shows that $u$ attains the supremum
$\textnormal{SCMT}_{\operatorname{e}}(\alpha,\beta)$, and $\|u\|_Q=1$. Hence, $u$ is a maximizer, and the proof is complete.\qed

\section{Existence of Normalized Ground State Solution}\label{sec8}
In this section, we deal with the existence of normalized ground states of \eqref{main problem} as mentioned in Theorem \ref{thm1.11}. The proof proceeds by characterizing the mountain pass level variationally, establishing compactness for the corresponding Palais--Smale sequences, and identifying their nontrivial limit through the limiting Euler-Lagrange equation. The Poho\v{z}aev characterization then provides the ground-state property. In what follows, we assume, without further mention, that \eqref{g1}--\eqref{g4} are satisfied. The hypothesis \eqref{g5} is used only in Lemma \ref{lem8.13}\textnormal{(ii)}, Lemma \ref{lem8.14} and the last step of the proof of Theorem \ref{thm1.11}, and we state it there. Next, we begin with a preliminary result that will be used repeatedly below.
\begin{lemma}\label{lem8.1}
 Let $0<\beta<Q$, $a>0$, $0<m<1$, $\vartheta\geq1$ and $\alpha>0$ be such that $\vartheta\alpha m^{\frac{1}{Q-1}}<\alpha_Q$. Then for every $u\in HW^{1,Q}(\mathbb{H}^n)$ with $\|u\|_Q=a$ and $\|D_H u\|_Q^Q\leq m$, there holds
 \begin{align*}
  \int_{\mathbb{H}^n}\frac{\left[\Phi\big(\alpha \big(1-\frac{\beta}{Q}\big)|u|^{Q^\prime}\big)\right]^\vartheta}{r(\xi)^\beta}\,\mathrm{d}\xi\leq \frac{a^{Q-\beta}}{m^{1-\frac{\beta}{Q}}}\,\textnormal{SCMT}(\vartheta\alpha m^{\frac{1}{Q-1}},\beta).
 \end{align*}
 In particular, the above integral is uniformly bounded on $\{u\in S(a)\colon\,\|D_Hu\|_Q^Q\leq m\}$.
\end{lemma}
\begin{proof}
By the given hypotheses, $a>0$ and thus, we have $u\neq 0$. In addition, one has $t:=\|D_H u\|_Q\in (0, m^\frac{1}{Q}]$. By setting $v=\frac{u}{t}$, it follows that $\|D_H v\|_Q=1$ and $\|v\|_Q=\frac{a}{t}$. Define $\gamma:=\vartheta\alpha t^{Q^\prime}$ and $\gamma_0:=\vartheta\alpha m^{\frac{1}{Q-1}}$. Recalling that $t^Q\leq m$, we obtain $t^{Q^\prime}\leq m^{\frac{1}{Q-1}}$. This yields that $ \gamma\leq\gamma_0<\alpha_Q$ and $(\frac{\gamma}{\gamma_0})^{Q-1}=\frac{t^Q}{m}\leq 1$. Now, by applying Lemma \ref{lem2.2}, one sees that $\Phi(s)^\vartheta\leq \Phi(\vartheta s)$ and $\Phi(\zeta s)\leq \zeta^{Q-1}\Phi(s)$ for all $0\leq\zeta\leq 1$. Hence, we deduce from Theorem \ref{thm1.3} that
\begin{align*}
  \int_{\mathbb{H}^n}&\frac{\left[\Phi\big(\alpha \big(1-\frac{\beta}{Q}\big)|u|^{Q^\prime}\big)\right]^\vartheta}{r(\xi)^\beta}\,\mathrm{d}\xi\leq  \int_{\mathbb{H}^n}\frac{\Phi\big(\gamma \big(1-\frac{\beta}{Q}\big)|v|^{Q^\prime}\big)}{r(\xi)^\beta}\,\mathrm{d}\xi\\
  &\leq \bigg(\frac{\gamma}{\gamma_0}\bigg)^{Q-1}\int_{\mathbb{H}^n}\frac{\Phi\big(\gamma_0 \big(1-\frac{\beta}{Q}\big)|v|^{Q^\prime}\big)}{r(\xi)^\beta}\,\mathrm{d}\xi\leq \frac{t^Q}{m}\,\textnormal{SCMT}(\gamma_0,\beta)\bigg(\frac{a}{t}\bigg)^{Q-\beta}\\
  &=\frac{a^{Q-\beta} t^\beta}{m} \,\textnormal{SCMT}(\gamma_0,\beta)\leq \frac{a^{Q-\beta}}{m^{1-\frac{\beta}{Q}}}\,\textnormal{SCMT}(\vartheta\alpha m^{\frac{1}{Q-1}},\beta).
\end{align*}
This completes the proof.
\end{proof}
The following lemma provides a compactness result that will be employed in the proof of the subsequent result.
\begin{lemma}\label{lem8.2}
Let $0<\beta<Q$, $a>0$ and $\{u_k\}_k\subset S(a)$ be a sequence with $ \limsup_{k\to\infty}\|D_H u_k\|_Q^Q<1$.
If $u_k\rightharpoonup u$ in  $HW^{1,Q}(\mathbb{H}^n)$ as $k\to\infty$ and $u_k\to u$ a.e. in $\mathbb{H}^n$ as $k\to\infty$, then there exists $\alpha>\alpha_Q$ such that for all $q>Q$, there holds
 \begin{align*}
     \frac{|u_k|^q\Phi\big(\alpha \big(1-\frac{\beta}{Q}\big)|u_k|^{Q^\prime}\big)}{r(\xi)^\beta}\to \frac{|u|^q\Phi\big(\alpha \big(1-\frac{\beta}{Q}\big)|u|^{Q^\prime}\big)}{r(\xi)^\beta}\quad\text{in}\quad L^1(\mathbb{H}^n)\quad\text{as}\quad k\to\infty.
 \end{align*}
\end{lemma}
\begin{proof}
Assume that the given hypotheses are satisfied. Next, we fix $m\in(0,1)$ such that $\limsup_{k\to\infty}\|D_H u_k\|_Q^Q<m<1.$ It follows that $\|D_H u_k\|_Q^Q\leq m$ for $k$ large enough. Recalling that $m<1$, we may choose $\alpha>\alpha_Q$ close to $\alpha_Q$ and $\vartheta>1$ close to $1$ so that $\vartheta\alpha m^{\frac{1}{Q-1}}<\alpha_Q$. For each $k\in\mathbb{N}$, we define $g_k\colon\mathbb{H}^n\to\mathbb{R}$ and $g\colon\mathbb{H}^n\to\mathbb{R}$ as follows
\begin{align*}
    g_k(\xi):=\frac{\Phi\big(\alpha \big(1-\frac{\beta}{Q}\big)|u_k|^{Q^\prime}\big)}{r(\xi)^{\frac{\beta}{\vartheta}}}\quad\text{and}\quad g(\xi):=\frac{\Phi\big(\alpha \big(1-\frac{\beta}{Q}\big)|u|^{Q^\prime}\big)}{r(\xi)^{\frac{\beta}{\vartheta}}}
\end{align*}
for a.e. $\xi\in\mathbb{H}^n$. It follows that $g_k\to g$ a.e. in $\mathbb{H}^n$ as $k\to\infty$. Moreover, an argument analogous to that in Lemma \ref{lem8.1} together with Lemma \ref{lem3.1} shows that the sequence $\{g_k\}_k$ is bounded in $L^\vartheta(\mathbb{H}^n)$ and thus, we  obtain $g_k\rightharpoonup g$ in $L^\vartheta(\mathbb{H}^n)$ as $k\to\infty$. Take $\vartheta^\prime>1$ such that $\frac{1}{\vartheta}+\frac{1}{\vartheta^\prime}=1$. The hypotheses together with Lemma \ref{lem4.2} imply that $u_k\to u$ in  $L^{q\vartheta^\prime}(\mathbb{H}^n,r(\xi)^{-\beta}\mathrm{d}\xi)$ as $k\to\infty$ and hence, there holds $u_k\,r(\xi)^{-\frac{\beta}{q\vartheta^\prime}}\rightarrow u\, r(\xi)^{-\frac{\beta}{q\vartheta^\prime}}$ in $L^{q\vartheta^\prime}(\mathbb{H}^n)$ as $k\to\infty$. Moreover, the continuity of the map $f\mapsto |f|^q$ from $L^{q\vartheta^\prime}(\mathbb{H}^n)$ into $L^{\vartheta^\prime}(\mathbb{H}^n)$ yields $|u_k|^qr(\xi)^{-\frac{\beta}{\vartheta^\prime}}\rightarrow |u|^q r(\xi)^{-\frac{\beta}{\vartheta^\prime}} $ in $L^{\vartheta^\prime}(\mathbb{H}^n)$ as $k\to\infty$. Now, by considering all the preceding arguments, we obtain
\begin{align*}
   \lim_{k\to\infty} \int_{\mathbb{H}^n}\frac{|u_k|^q\Phi\big(\alpha \big(1-\frac{\beta}{Q}\big)|u_k|^{Q^\prime}\big)}{r(\xi)^\beta}\,\mathrm{d}\xi=\int_{\mathbb{H}^n}\frac{|u|^q\Phi\big(\alpha \big(1-\frac{\beta}{Q}\big)|u|^{Q^\prime}\big)}{r(\xi)^\beta}\,\mathrm{d}\xi.
\end{align*}
Hence, the desired convergence follows from  $u_k\to u$ a.e. in $\mathbb{H}^n$ as $k\to\infty$ and Scheff\'{e}'s lemma \cite[{\S}5.10(ii), p.~55]{Williams-1991}. This finishes the proof.
\end{proof}
The compactness established above yields the following convergence of the nonlinear terms.
\begin{lemma}\label{lem8.3}
Let $0<\beta<Q$, $a>0$ and $\{u_k\}_k\subset S(a)$ be a sequence with $ \limsup_{k\to\infty}\|D_H u_k\|_Q^Q<1$.
If $u_k\rightharpoonup u$ in  $HW^{1,Q}(\mathbb{H}^n)$ as $k\to\infty$ and $u_k\to u$ a.e. in $\mathbb{H}^n$ as $k\to\infty$, then there hold
 \begin{align}\label{eq9.9999}
     \frac{G(u_k)}{r(\xi)^\beta}\to  \frac{G(u)}{r(\xi)^\beta}\quad\text{and }\quad\frac{g(u_k)u_k}{r(\xi)^\beta}\to \frac{g(u)u}{r(\xi)^\beta}\quad\text{in}\quad L^1(\mathbb{H}^n)\quad\text{as}\quad k\to\infty.
 \end{align}
\end{lemma}
\begin{proof}
Assume that the stated hypotheses hold. It follows from Lemma \ref{lem4.2} that $u_k\to u$ in  $L^{\tau}(\mathbb{H}^n,r(\xi)^{-\beta}\mathrm{d}\xi)$ as $k\to\infty$ for all $\tau>2Q-\beta$. Now, arguing similarly as in Lemma \ref{lem8.2}, we obtain $|u_k|^\tau r(\xi)^{-\beta}\to |u|^\tau r(\xi)^{-\beta}$ in $L^{1}(\mathbb{H}^n)$ as $k\to\infty$. By applying Lemma \ref{lem8.2}, whenever $\alpha>\alpha_Q$ and $q>Q$, we have
\begin{align*}
     \frac{|u_k|^q\Phi\big(\alpha \big(1-\frac{\beta}{Q}\big)|u_k|^{Q^\prime}\big)}{r(\xi)^\beta} \to\frac{|u|^q\Phi\big(\alpha \big(1-\frac{\beta}{Q}\big)|u|^{Q^\prime}\big)}{r(\xi)^\beta}\quad\text{in}\quad L^1(\mathbb{H}^n)\quad\text{as}\quad k\to\infty.
 \end{align*}
Taking into account \eqref{eq1.6}, $u_k\to u$ a.e. in $\mathbb{H}^n$ as $k\to\infty$, and the preceding convergences, the generalized dominated convergence theorem of Lebesgue, as stated in Royden--Fitzpatrick \cite[{\S}4.4, Theorem 19, p. 89]{Royden-Fitzpatrick-1988}, ensures the validity of \eqref{eq9.9999}. This completes the proof.
\end{proof}
\subsection{The Minimax Strategy} To analyze the geometry of $\mathcal{J}$ on the sphere $S(a)$, we consider the scaling operator $\mathcal{F}\colon HW^{1,Q}(\mathbb{H}^n)\times \mathbb{R}\to HW^{1,Q}(\mathbb{H}^n)$ defined by
\begin{align*}
   \mathcal{F}(u,s)(\xi):=\theta u(\delta_\theta(\xi)) \quad\text{with}\quad\theta:=\operatorname{\exp}(s)\quad\text{for a.e.}\quad\xi\in \mathbb{H}^n,
\end{align*}
for all $(u,s)\in  HW^{1,Q}(\mathbb{H}^n)\times \mathbb{R}$. This scaling preserves the $L^Q$-norm, in the sense that $\|\mathcal{F}(u,s)\|_Q=\|u\|_Q$ for all $s\in\mathbb{R}$. The functional $\widetilde{\mathcal{J}}\colon  HW^{1,Q}(\mathbb{H}^n)\times \mathbb{R}\to \mathbb{R}$ is defined  through $ \widetilde{\mathcal{J}}(u,s)=\mathcal{J}(\mathcal{F}(u,s))$ such that
\begin{align*}
 \widetilde{\mathcal{J}}(u,s)=\frac{\operatorname{\exp}(Qs)}{Q}\|D_Hu\|_Q^Q- \frac{1}{\operatorname{\exp}((Q-\beta)s)}\int_{\mathbb{H}^n}\frac{G(\operatorname{\exp}(s)u)}{r(\xi)^\beta} \,\mathrm{d}\xi
\end{align*}
for all $(u,s)\in  HW^{1,Q}(\mathbb{H}^n)\times \mathbb{R}$. By differentiating $\widetilde{\mathcal{J}}$ with respect to $s$, we obtain $\partial_s \widetilde{\mathcal{J}}(u,s)=\mathcal{P}(\mathcal{F}(u,s))$ for all $(u,s)\in  HW^{1,Q}(\mathbb{H}^n)\times \mathbb{R}$, where the Poho\v{z}aev functional $\mathcal{P}$ defined in \eqref{eq9.999}.

Next, we establish the mountain pass type geometry of $\widetilde{\mathcal{J}}$ on $S(a)\times\mathbb{R}$. To this end, we prove the following two preliminary results.
\begin{lemma}\label{lem8.4}
Let $u\in S(a)$ be any arbitrarily fixed function in $HW^{1,Q}(\mathbb{H}^n)$. Then the following results hold:
\begin{itemize}
    \item[\textnormal{(i)}] $\|D_H \mathcal{F}(u,s)\|_Q\to 0$ and $\mathcal{J}(\mathcal{F}(u,s))\to 0$ as $s\to-\infty$.
    \item[\textnormal{(ii)}] $\|D_H \mathcal{F}(u,s)\|_Q\to \infty$ and $\mathcal{J}(\mathcal{F}(u,s))\to -\infty$ as $s\to\infty$.
\end{itemize}
\end{lemma}
\begin{proof}
By direct computations, we obtain
\begin{align*}
    \int_{\mathbb{H}^n}|\mathcal{F}(u,s)|^Q\,\mathrm{d}\xi=a^Q,~~  \int_{\mathbb{H}^n}|D_H \mathcal{F}(u,s)|_H^Q\,\mathrm{d}\xi=\operatorname{\exp}(Qs)\int_{\mathbb{H}^n}|D_H u|_H^Q\,\mathrm{d}\xi
\end{align*}
and
\begin{align*}
   \int_{\mathbb{H}^n}\frac{|\mathcal{F}(u,s)|^\zeta}{r(\xi)^\beta}\,\mathrm{d}\xi =\operatorname{\exp}((\zeta-Q+\beta)s) \int_{\mathbb{H}^n}\frac{|u|^\zeta}{r(\xi)^\beta}\,\mathrm{d}\xi\quad\text{for all}\quad\zeta\geq Q.
\end{align*}
For any fixed $\zeta>Q$, the above estimates imply that
\begin{align*}
\|D_H \mathcal{F}(u,s)\|_Q\to 0\qquad\text{and}\qquad \|\mathcal{F}(u,s)\|_{\zeta,\beta}\to 0  \quad\text{as}\quad s\to-\infty.
\end{align*}
Since $\|D_H \mathcal{F}(u,s)\|_Q\to 0$ as $s\to-\infty$, there exist $s_0<0$ and $m\in(0,1)$ such that $\|D_H \mathcal{F}(u,s)\|^Q_Q\leq m$ for all $s\leq s_0$. Now, we assume that $\vartheta,\,\vartheta^\prime>1$ satisfying $\frac{1}{\vartheta}+\frac{1}{\vartheta^\prime}=1$. Take $\alpha>\alpha_Q$ close to $\alpha_Q$, $\vartheta>1$ close to $1$ such that $\vartheta\alpha m^{\frac{1}{Q-1}}<\alpha_Q$. By applying Lemma \ref{lem8.1}, we can find a uniform constant $C=C(Q,\vartheta,a,m,\beta)>0$ such that
\begin{align*}
    \Bigg(\int_{\mathbb{H}^n}\frac{\left[\Phi\big(\alpha \big(1-\frac{\beta}{Q}\big)|\mathcal{F}(u,s)|^{Q^\prime}\big)\right]^\vartheta}{r(\xi)^\beta}\,\mathrm{d}\xi\Bigg)^{\frac{1}{\vartheta}}\leq C\quad\text{for all}\quad s\leq s_0.
\end{align*}
As a result, we obtain from \eqref{eq1.6} and H\"{o}lder's inequality that
\begin{align*}
    \int_{\mathbb{H}^n}\frac{G(\mathcal{F}(u,s))}{r(\xi)^\beta}\,\mathrm{d}\xi\leq \varepsilon\|\mathcal{F}(u,s)\|^\tau_{\tau,\beta}+C_\varepsilon C \|\mathcal{F}(u,s)\|^q_{q\vartheta^\prime,\beta}\quad\text{for all}\quad s\leq s_0.
\end{align*}
This allows us to conclude that
\begin{align*}
  \|D_H \mathcal{F}(u,s)\|_Q\to 0\qquad\text{and}\qquad \int_{\mathbb{H}^n}\frac{G(\mathcal{F}(u,s))}{r(\xi)^\beta}\,\mathrm{d}\xi\to 0  \quad\text{as}\quad s\to-\infty.
\end{align*}
Hence, we obtain $\mathcal{J}(\mathcal{F}(u,s))\to 0$ as $s\to-\infty$, as asserted in \textnormal{(i)}. Next, one can observe that $\|D_H \mathcal{F}(u,s)\|_Q\to \infty$ as $s\to\infty$. By using \eqref{g4}, $p>2Q-\beta$, and the previous estimates, we have
\begin{equation}\label{eq8.999}
\begin{aligned}
 \mathcal{J}(\mathcal{F}(u,s))&\leq \frac{1}{Q} \|D_H \mathcal{F}(u,s)\|_Q^Q-\frac{\mu}{p}  \int_{\mathbb{H}^n}\frac{|\mathcal{F}(u,s)|^p}{r(\xi)^\beta}\,\mathrm{d}\xi \\
 &=\frac{\operatorname{\exp}(Qs)}{Q}\|D_H u\|_Q^Q-\frac{\mu}{p}\operatorname{\exp}((p-Q+\beta)s)\|u\|^p_{p,\beta}
   \to-\infty  
\end{aligned}
\end{equation}
as $s\to\infty$. This establishes \textnormal{(ii)}, which completes the proof.
\end{proof}
Next, we establish the energy separation needed for the subsequent variational analysis.
\begin{lemma}\label{lem8.5}
There exists $K>0$ small enough such that $$ 0<\sup_{u\in\mathcal{A}}\mathcal{J}(u)<\inf_{u\in\mathcal{B}}\mathcal{J}(u),$$ where
\begin{align*}
 \mathcal{A}:=\left\{u\in S(a)\colon\, \|D_Hu\|^Q_Q\leq K\right\} \quad\text{and} \quad \mathcal{B}:=\left\{u\in S(a)\colon\, \|D_Hu\|^Q_Q=2K\right\}.
\end{align*}
\end{lemma}
\begin{proof}
The Caffarelli--Kohn--Nirenberg interpolation-type inequality that already established in Theorem \ref{thm1.5} yields
\begin{align}\label{eq8.1}
  \|u\|_{\vartheta,\beta}\leq C(Q,\beta,\vartheta) \|u\|_Q^{\frac{Q-\beta}{\vartheta}}\|D_Hu\|_Q^{1-\frac{Q-\beta}{\vartheta}}
\end{align}
for all $u\in HW^{1,Q}(\mathbb{H}^n)$, where $0<\beta<Q$, $\vartheta\geq Q$ and $C(Q,\beta,\vartheta)>0$ is a constant  depending only on $Q$, $\beta$ and $\vartheta$. Fix $m=\frac{1}{2}$ and let $u\in S(a)$ be such that $\|D_H u\|_Q^Q\leq m$. Take $t,\,t^\prime>1$ satisfying $\frac{1}{t}+\frac{1}{t^\prime}=1$. Next, we choose $\alpha>\alpha_Q$ close to $\alpha_Q$, $t>1$ close to $1$ so that $t\alpha m^{\frac{1}{Q-1}}<\alpha_Q$. Note that $m$, $\alpha$ and $t$ are fixed before $K$ is chosen. By using \eqref{eq1.6}, H\"{o}lder's inequality, and Lemma \ref{lem8.1}, we deduce that
\begin{align}\label{eq8.2}
    \int_{\mathbb{H}^n}\frac{G(u)}{r(\xi)^\beta}\,\mathrm{d}\xi\leq \varepsilon\|u\|^\tau_{\tau,\beta}+\Tilde{C}\|u\|^q_{qt^\prime,\beta},
\end{align}
where $\tau>2Q-\beta$, $q>Q$, and $\Tilde{C}>0$ is a suitable constant depending only on $Q$, $\beta$, $a$ and $g$ (through $\varepsilon$, $C_\varepsilon$, $\alpha$, $q$ and $t$), but not on $K$. In view of \eqref{eq8.1} and \eqref{eq8.2}, one sees that
\begin{equation}\label{eq8.3}
\begin{aligned}
\int_{\mathbb{H}^n}\frac{G(u)}{r(\xi)^\beta}\,\mathrm{d}\xi&\leq \varepsilon C(Q,\beta,\tau)^\tau \|u\|_Q^{\tau(\frac{Q-\beta}{\tau})}\|D_Hu\|_Q^{\tau(1-\frac{Q-\beta}{\tau})} \\
&\qquad+\Tilde{C} C(Q,\beta,qt^\prime)^q \|u\|_Q^{q(\frac{Q-\beta}{qt^\prime})}\|D_Hu\|_Q^{q(1-\frac{Q-\beta}{qt^\prime})}\\
&=C_1\big(\|D_Hu\|_Q^Q\big)^{\frac{\tau-Q+\beta}{Q}} a^{Q-\beta}+C_2\big(\|D_Hu\|_Q^Q\big)^{\frac{qt^\prime-Q+\beta}{Qt^\prime}} a^{\frac{Q-\beta}{t^\prime}}, 
\end{aligned}
\end{equation}
where $C_1=\varepsilon C(Q,\beta,\tau)^\tau$ and $C_2=\Tilde{C} C(Q,\beta,qt^\prime)^q$. It follows from \eqref{g3} that $G(u)\geq 0$ for all $u\in HW^{1,Q}(\mathbb{H}^n)$. In particular, for any $u,\,v\in HW^{1,Q}(\mathbb{H}^n)$, we obtain from \eqref{eq8.9990} that
\begin{equation}\label{eq8.4}
\mathcal{J}(v)-\mathcal{J}(u)\geq\frac{1}{Q}\|D_H v\|_Q^Q-\frac{1}{Q}\|D_H u\|_Q^Q -\int_{\mathbb{H}^n}\frac{G(v)}{r(\xi)^\beta}\,\mathrm{d}\xi.
\end{equation}
Recall that $K\leq\frac{1}{4}$ (see \eqref{eq8.5}). Let $u,\,v\in HW^{1,Q}(\mathbb{H}^n)$ be such that $\|D_H u\|_Q^Q\leq K$ and $\|D_H v\|_Q^Q= 2K$. Thus, the estimates \eqref{eq8.3} and \eqref{eq8.4} together imply that
\begin{align*}
 \mathcal{J}(v)-\mathcal{J}(u)&\geq\frac{K}{Q} -C_1\left(2K\right)^{\frac{\tau-Q+\beta}{Q}} a^{Q-\beta}-C_2\left(2K\right)^{\frac{qt^\prime-Q+\beta}{Qt^\prime}} a^{\frac{Q-\beta}{t^\prime}}\\
 &=K\Big(\frac{1}{Q}-C_12^{\frac{\tau-Q+\beta}{Q}} K^{\frac{\tau-Q+\beta}{Q}-1}a^{Q-\beta}-C_2 2^{\frac{qt^\prime-Q+\beta}{Qt^\prime}} K^{\frac{qt^\prime-Q+\beta}{Qt^\prime}-1} a^{\frac{Q-\beta}{t^\prime}}\Big).
\end{align*}
Notice that $0<\beta<Q$, $\tau-Q+\beta>Q$ and $t^\prime=\frac{t}{t-1}$. It follows that $(q-Q)t^\prime>Q-\beta$ for any $q>Q$ once $t>1$ close to $1$. Define the constant $K$ by
\begin{align}\label{eq8.5}
    K:=\min\Bigg\{\frac{1}{4},\Bigg(\frac{1}{4QC_1 2^{\frac{\tau-Q+\beta}{Q}}a^{Q-\beta}}\Bigg)^{\frac{Q}{\tau-2Q+\beta}},\Bigg(\frac{1}{4QC_2 2^{\frac{qt^\prime-Q+\beta}{Qt^\prime}} a^{\frac{Q-\beta}{t^\prime}}}\Bigg)^{\frac{Qt^\prime}{(q-Q)t^\prime-Q+\beta}}\Bigg\}.
\end{align}
This yields that $ \mathcal{J}(v)-\mathcal{J}(u)\geq\frac{K}{2Q}>0$ and hence, the desired result follows. This finishes the proof. 
\end{proof}

As a byproduct of Lemma \ref{lem8.5}, we have the following result.
\begin{corollary}\label{cor8.6}
 Let $K$ be given in \eqref{eq8.5}. Then there holds $\mathcal{J}(u)>0$ for all $u\in S(a)$ with $\|D_H u\|_Q^Q\leq K$. Furthermore, we also have
 \begin{align*}
   \mathcal{J}^\ast:=\inf\left\{\mathcal{J}(u)\colon\, u\in S(a)\quad\text{and}\quad \|D_H u\|_Q^Q=\frac{K}{2}\right\}>0.
 \end{align*}
\end{corollary}
\begin{proof}
Fix $u\in S(a)$ with $\|D_Hu\|_Q^Q\leq K$, where $K>0$ is given in \eqref{eq8.5}. By an argument analogous to that in Lemma \ref{lem8.5}, we deduce from \eqref{eq8.3} that
\begin{align*}
 \mathcal{J}(u)&\geq \frac{1}{Q}\|D_H u\|_Q^Q -C_1\big(\|D_Hu\|_Q^Q\big)^{\frac{\tau-Q+\beta}{Q}} a^{Q-\beta}-C_2\big(\|D_Hu\|_Q^Q\big)^{\frac{qt^\prime-Q+\beta}{Qt^\prime}} a^{\frac{Q-\beta}{t^\prime}}\\
 &\geq\|D_H u\|_Q^Q \bigg(\frac{1}{Q}-C_1 K^{\frac{\tau-Q+\beta}{Q}-1}a^{Q-\beta}-C_2 K^{\frac{qt^\prime-Q+\beta}{Qt^\prime}-1} a^{\frac{Q-\beta}{t^\prime}}\bigg)\\
 &\geq \frac{1}{Q}\|D_H u\|_Q^Q\Bigg(1-\frac{1}{4}\bigg(\frac{1}{ 2^{\frac{\tau-Q+\beta}{Q}}}+\frac{1}{2^{\frac{qt^\prime-Q+\beta}{Qt^\prime}}}\bigg)\Bigg)\geq \frac{1}{2Q}\|D_H u\|_Q^Q>0.
\end{align*}
This shows that $ \mathcal{J}^\ast\geq \frac{K}{4Q}>0$. Hence, the proof is complete.
\end{proof}
Next, we fix $u_0\in S(a)$. In view of Lemma \ref{lem8.4} and Corollary \ref{cor8.6}, we can find $s_1=s_1(u_0,a,\beta,\mu)<0$ and $s_2=s_2(u_0,a,\beta,\mu)>0$ such that the functions $u_{1,\mu}=\mathcal{F}(u_0,s_1)$ and $u_{2,\mu}=\mathcal{F}(u_0,s_2)$ satisfy $\|D_H u_{1,\mu}\|_Q^Q<\frac{K}{2}$, $\|D_H u_{2,\mu}\|_Q^Q> 2K$, $\mathcal{J}(u_{1,\mu})>0$ and $\mathcal{J}(u_{2,\mu})<0$, where $K$ is given in \eqref{eq8.5}. To proceed, we follow the approach of Jeanjean \cite{Jeanjean-1997} and define the following mountain pass level by
\begin{align*}
\gamma_\mu(a,\beta):=\inf_{\psi\in\Gamma}\max_{t\in[0,1]} \mathcal{J}(\psi(t)),
\end{align*}
where
\begin{align*}
    \Gamma:=\left\{\psi\in C\left([0,1],S(a)\right):\, \|D_H \psi(0)\|_Q^Q<\frac{K}{2}\quad\text{and}\quad \mathcal{J}(\psi(1))<0\right\}.
\end{align*}
Then, for each $\psi\in\Gamma$, it follows from Corollary \ref{cor8.6} that $\|D_H \psi(1)\|_Q^Q>K$. Hence, there exists some $t_0\in(0,1)$ such that $\|D_H \psi(t_0)\|_Q^Q=\frac{K}{2}$ and consequently, we obtain
\begin{align*}
    \max_{t\in[0,1]} \mathcal{J}(\psi(t))\geq \mathcal{J}(\psi(t_0))\geq \mathcal{J}^\ast>0.
\end{align*}
This directly implies that $\gamma_\mu(a,\beta)\geq  \mathcal{J}^\ast>0$.

Now, we describe the mountain pass level's asymptotic behavior.
\begin{lemma}\label{lem8.7}
 Let $\kappa=\frac{p-Q+\beta}{Q}>1$. For every $u_0\in S(a)$, setting $A_0=\|D_H u_0\|_Q^Q$ and $B_0=\|u_0\|^p_{p,\beta}$, there holds
 \begin{align}\label{eq8.gamma}
 \gamma_\mu(a,\beta)\leq\left(1-\frac{1}{\kappa}\right)\left(\frac{A_0}{Q}\right)^{\frac{\kappa}{\kappa-1}}\left(\frac{\kappa\mu B_0}{p}\right)^{-\frac{1}{\kappa-1}}.
 \end{align}
 In particular, $\gamma_\mu(a,\beta)\to 0$ as $\mu\to\infty$.
\end{lemma}
\begin{proof} Let $u_0\in S(a)$, and let $A_0$, $B_0$ and $\kappa$ be as above. Note that $\kappa>1$ since $p>2Q-\beta$. Define the path $\psi_0\colon[0,1]\to S(a)$ by $\psi_0(t)=\mathcal{F}(u_0,(1-t)s_1+ts_2)\in\Gamma$ for all $t\in[0,1]$, where $s_1$ and $s_2$ are chosen for this $u_0$ as above. By repeating the argument used in \eqref{eq8.999} and invoking \eqref{g4} with $p>2Q-\beta$, we deduce that
\begin{align*}
\gamma_\mu(a,\beta)&\leq\max_{t\in[0,1]} \mathcal{J}(\psi_0(t))\leq \max_{\zeta> 0}\left(\frac{A_0}{Q}\zeta-\frac{\mu B_0}{p}\,\zeta^{\kappa}\right)\\
&=\left(1-\frac{1}{\kappa}\right)\left(\frac{A_0}{Q}\right)^{\frac{\kappa}{\kappa-1}}\left(\frac{\kappa\mu B_0}{p}\right)^{-\frac{1}{\kappa-1}}
\to 0 \quad\text{as}\quad\mu\to\infty.
\end{align*}
This finishes the proof.
\end{proof}
By adapting the argument of Jeanjean \cite[Proposition 2.2 and Lemma 2.4]{Jeanjean-1997}, let $\{u_k\}_k\subset S(a)$ be the $\textnormal{(PS)}$-sequence for $\mathcal{J}$ at the level $\gamma_\mu(a,\beta)$, obtained by setting $u_k=\mathcal{F}(v_k,s_k)$, where $\{(v_k,s_k)\}_k\subset S(a)\times\mathbb{R}$ is the $\textnormal{(PS)}$-sequence for $\widetilde{\mathcal{J}}$ at the same level. In particular, we have $\mathcal{J} (u_k)=\gamma_\mu(a,\beta)+o_k(1)$ as $k\to\infty$ and
\begin{align*}
  \left\|\left.\mathcal{J}^\prime\right|_{S(a)}(u_k)\right\|_\ast=\sup_{v\in T_{u_k}S(a),\,\|v\|=1}\left|\langle{\mathcal{J}^\prime (u_k),v}\rangle\right|=o_k(1)\quad\text{as}\quad k\to\infty,
\end{align*}
 where $T_{u}S(a)=\left\{v\in HW^{1,Q}(\mathbb{H}^n)\colon\,\int_{\mathbb{H}^n}|u|^{Q-2}uv\,\mathrm{d}\xi=0\right\}$ is the tangent space of $S(a)$ at $u$. Next, we define the functional $\Psi\colon HW^{1,Q}(\mathbb{H}^n)\to\mathbb{R}$ by
\begin{align}\label{eq8.9999}
    \Psi(u)=\frac{1}{Q}\int_{\mathbb{H}^n}|u|^Q\,\mathrm{d}\xi\quad\text{for all}\quad u\in HW^{1,Q}(\mathbb{H}^n).
\end{align}
It is trivial to see that $S(a)=\Psi^{-1}(\{a^Q/Q\})$. Thus, there exists some sequence $\{\lambda_k\}_k\subset \mathbb{R}$ such that $ \|\mathcal{J}^\prime(u_k)-\lambda_k\Psi^\prime(u_k)\|_\ast=o_k(1)$ as $k\to\infty$. As a result, we obtain
\begin{align*}
    -\Delta_{H,Q} u_k-\lambda_k|u_k|^{Q-2}u_k
=\frac{g(u_k)}{r(\xi)^\beta}+o_k(1)
\quad\text{in}\quad\left(HW^{1,Q}(\mathbb{H}^n)\right)^\ast
\end{align*}
as $k\to\infty$. Evidently, there  holds $\partial_s \widetilde{\mathcal{J}}(v_k,s_k)=o_k(1)$ as $k\to\infty$. This  gives $\mathcal{P}(u_k)=o_k(1)$ as $k\to\infty$, that is, one has
\begin{align*}
 \|D_Hu_k\|_Q^Q+ \int_{\mathbb{H}^n}\frac{\left((Q-\beta)G(u_k)-g(u_k)u_k\right)}{r(\xi)^\beta}\,\mathrm{d}\xi=o_k(1)\quad\text{as}\quad k\to\infty.
\end{align*}
 It follows from Lemmas \ref{lem8.8} and \ref{lem8.9} that $\{u_k\}_k$ is a bounded sequence in $HW^{1,Q}(\mathbb{H}^n)$, while the sequence  $\{\lambda_k\}_k\subset \mathbb{R}$ satisfies
\begin{align*}
    \lambda_k=\frac{1}{a^Q}\left(\|D_Hu_k\|_Q^Q-\int_{\mathbb{H}^n}\frac{g(u_k)u_k}{r(\xi)^\beta}\,\mathrm{d}\xi\right)+o_k(1)\quad\text{as}\quad k\to\infty.
\end{align*}

The nonlinear term admits the following uniform estimate along the Palais--Smale sequence.
\begin{lemma}\label{lem8.8}
There holds
\begin{align*}
    \limsup_{k\to\infty}\int_{\mathbb{H}^n}\frac{G(u_k)}{r(\xi)^\beta}\,\mathrm{d}\xi\leq \left(\frac{Q}{\theta-2Q+\beta}\right) \gamma_\mu(a,\beta).
\end{align*}
\end{lemma}
\begin{proof}
It follows from $\mathcal{J} (u_k)=\gamma_\mu(a,\beta)+o_k(1)$ and $\mathcal{P}(u_k)=o_k(1)$ as $k\to\infty$ that
\begin{align*}
  Q\mathcal{J} (u_k)+ \mathcal{P}(u_k)=Q \gamma_\mu(a,\beta)+o_k(1)\quad\text{as}\quad k\to\infty.
\end{align*}
This yields at once that
\begin{align*}
  2\|D_Hu_k\|_Q^Q-  \int_{\mathbb{H}^n}\frac{\left(\beta G(u_k)+g(u_k)u_k\right)}{r(\xi)^\beta}\,\mathrm{d}\xi=Q\gamma_\mu(a,\beta)+o_k(1)
\end{align*}
as $k\to\infty$. Now, again by using $\mathcal{J} (u_k)=\gamma_\mu(a,\beta)+o_k(1)$ as $k\to\infty$, we obtain from the above estimate that
\begin{align*}
   2Q\gamma_\mu(a,\beta)+  \int_{\mathbb{H}^n}\frac{\left(\left(2Q-\beta\right) G(u_k)-g(u_k)u_k\right)}{r(\xi)^\beta}\,\mathrm{d}\xi=Q\gamma_\mu(a,\beta)+o_k(1)
\end{align*}
as $k\to\infty$. Thus, we deduce from \eqref{g3} that
\begin{align*}
 Q\gamma_\mu(a,\beta)+o_k(1)\geq (\theta-2Q+\beta) \int_{\mathbb{H}^n}\frac{G(u_k)}{r(\xi)^\beta}\,\mathrm{d}\xi \quad\text{as}\quad k\to\infty.
\end{align*}
Hence, the desired result follows from the above estimate, which completes the proof.
\end{proof}
The next estimate provides the strict gradient bound required to apply the subcritical weighted Moser--Trudinger inequality.
\begin{lemma}\label{lem8.9}
  There holds
  \begin{align*}
    \limsup_{k\to\infty}\|D_H u_k\|_Q^Q\leq Q\left(\frac{\theta-Q+\beta}{\theta-2Q+\beta}\right) \gamma_\mu(a,\beta).
\end{align*}
Furthermore, if $\mu>\mu^\ast$, where $\mu^\ast=\mu^\ast(Q,\beta,a,p,\theta)>0$ is given by \eqref{eq8.mustar} below, then 
\begin{align*}
    \limsup_{k\to\infty}\|D_H u_k\|_Q^Q<1.
\end{align*}
\end{lemma}
\begin{proof}
Recalling that $\mathcal{J}(u_k)=\gamma_\mu(a,\beta)+o_k(1)$ as $k\to\infty$, we deduce
\begin{align*}
  \|D_H u_k\|_Q^Q= Q\gamma_\mu(a,\beta)+Q\int_{\mathbb{H}^n}\frac{G(u_k)}{r(\xi)^\beta}\,\mathrm{d}\xi+o_k(1) \quad\text{as}\quad k\to\infty.
\end{align*}
This together with Lemma \ref{lem8.8} yields
\begin{align*}
    \limsup_{k\to\infty}\|D_H u_k\|_Q^Q\leq Q\left(\frac{\theta-Q+\beta}{\theta-2Q+\beta}\right) \gamma_\mu(a,\beta).
\end{align*}
Denote $L:=Q\left(\frac{\theta-Q+\beta}{\theta-2Q+\beta}\right)$ and $\kappa=\frac{p-Q+\beta}{Q}$. Now, we define $\mu^\ast=\mu^\ast(Q,\beta,a,p,\theta)$ explicitly as follows
\begin{align}\label{eq8.mustar}
 \mu^\ast:=\frac{p}{\kappa}\left[L\left(1-\frac{1}{\kappa}\right)\right]^{\kappa-1}\frac{1}{Q^{\kappa}\,\mathcal{C}_p^{\,p}\,a^{Q-\beta}}>0,
\end{align}
and $\mathcal{C}_p=C(Q,\beta,p)$ is the best constant in \eqref{eq1.3} with $\vartheta=p$, that is, we have
\begin{align*}
\mathcal{C}_p^{\,p}:=\sup_{u\in HW^{1,Q}(\mathbb{H}^n)\setminus\{0\}}\frac{\|u\|^p_{p,\beta}}{\|u\|_Q^{Q-\beta}\|D_H u\|_Q^{p-Q+\beta}}<+\infty,
\end{align*}
by Theorem \ref{thm1.5}, since $p>2Q-\beta>Q$. Indeed, for every $u_0\in S(a)$, we have $\frac{A_0^\kappa}{B_0}=\frac{\|D_H u_0\|_Q^{p-Q+\beta}}{\|u_0\|_{p,\beta}^p}$. If $v\in HW^{1,Q}(\mathbb{H}^n)\setminus\{0\}$ and $u_0=\frac{av}{\|v\|_Q}$, then $u_0\in S(a)$ and 
\begin{align*}
\frac{A_0^\kappa}{B_0}=\frac{1}{a^{Q-\beta}}\,\frac{\|v\|_Q^{Q-\beta}\|D_H v\|_Q^{p-Q+\beta}}{\|v\|^p_{p,\beta}}.
\end{align*}
Hence, we obtain $\inf_{u_0\in S(a)}\frac{A_0^\kappa}{B_0}=\frac{1}{\mathcal{C}_p^{\,p}\,a^{Q-\beta}}$. Next, let us take $\mu>\mu^\ast$. Then we can choose $u_0\in S(a)$ such that
\begin{align*}
\mu>\frac{p}{\kappa}\left[L\left(1-\frac{1}{\kappa}\right)\right]^{\kappa-1}\frac{1}{Q^{\kappa}}\,\frac{A_0^\kappa}{B_0}=\frac{p}{\kappa B_0}\left[L\left(1-\frac{1}{\kappa}\right)\right]^{\kappa-1}\left(\frac{A_0}{Q}\right)^{\kappa},
\end{align*}
and this inequality is equivalent to $$L\left(1-\frac{1}{\kappa}\right)\left(\frac{A_0}{Q}\right)^{\frac{\kappa}{\kappa-1}}\left(\frac{\kappa\mu B_0}{p}\right)^{-\frac{1}{\kappa-1}}<1.$$ Moreover, by using \eqref{eq8.gamma}, we get $L\,\gamma_\mu(a,\beta)<1$ for all $\mu>\mu^\ast$.
This shows that
\begin{align*}
    \limsup_{k\to\infty}\|D_H u_k\|_Q^Q<1\quad\text{for all}\quad\mu>\mu^\ast.
\end{align*}
Hence, the proof is complete.
\end{proof}
The boundedness and asymptotic behavior of the Lagrange multipliers are summarized below.
\begin{lemma}\label{lem8.10}
Assume that $\mu>\mu^\ast$, where $\mu^\ast$ is given in Lemma \ref{lem8.9}. Then the sequence $\{\lambda_k\}_k\subset\mathbb{R}$ is bounded. In addition, they hold
\begin{align*}
 \limsup_{k\to\infty}|\lambda_k|\leq \frac{2Q}{a^Q}\left(\frac{\theta-\frac{1}{2}\left(Q-\beta\right)}{\theta-2Q+\beta}\right) \gamma_\mu(a,\beta)
\end{align*}
and
\begin{align*}
 \limsup_{k\to\infty}\lambda_k= -\left(\frac{Q-\beta}{a^Q}\right)\liminf_{k\to\infty}\int_{\mathbb{H}^n}\frac{G(u_k)}{r(\xi)^\beta}\,\mathrm{d}\xi.
\end{align*}
\end{lemma}
\begin{proof}
It follows directly from $\mathcal{P}(u_k)=o_k(1)$ as $k\to\infty$ and Lemmas \ref{lem8.8} and \ref{lem8.9} that
\begin{align*}
 \limsup_{k\to\infty} \int_{\mathbb{H}^n}\frac{g(u_k)u_k}{r(\xi)^\beta}\,\mathrm{d}\xi \leq\left(\frac{\theta Q}{\theta-2Q+\beta}\right)\gamma_\mu(a,\beta).
\end{align*}
This together with  Lemma \ref{lem8.9} and the fact that
\begin{align}\label{eq8.7}
    \lambda_k=\frac{1}{a^Q}\left(\|D_H u_k\|_Q^Q-\int_{\mathbb{H}^n}\frac{g(u_k)u_k}{r(\xi)^\beta}\,\mathrm{d}\xi\right)+o_k(1)\quad\text{as}\quad k\to\infty,
\end{align}
we conclude that $\{\lambda_k\}_k\subset\mathbb{R}$ is a bounded sequence such that
\begin{align*}
 \limsup_{k\to\infty}|\lambda_k|\leq \frac{2Q}{a^Q}\left(\frac{\theta-\frac{1}{2}\left(Q-\beta\right)}{\theta-2Q+\beta}\right) \gamma_\mu(a,\beta).
\end{align*}
Moreover, again by utilizing $\mathcal{P}(u_k)=o_k(1)$ as $k\to\infty$, we obtain from \eqref{eq8.7} that
\begin{align*}
 \limsup_{k\to\infty}\lambda_k= -\bigg(\frac{Q-\beta}{a^Q}\bigg)\liminf_{k\to\infty}\int_{\mathbb{H}^n}\frac{G(u_k)}{r(\xi)^\beta}\,\mathrm{d}\xi.
\end{align*}
The proof is thus complete.
\end{proof}
The boundedness of the sequence $\{u_k\}_k$ in $HW^{1,Q}(\mathbb{H}^n)$ ensures that, up to a subsequence, there exists $u\in HW^{1,Q}(\mathbb{H}^n)$ such that
$u_k\rightharpoonup u \quad\text{in }HW^{1,Q}(\mathbb{H}^n)$ as $k\to\infty$. The result below guarantees that this weak limit is nonzero whenever $\mu>\mu^\ast$.

Now, we prove that the weak limit of the Palais--Smale sequence is nontrivial at the mountain-pass level as follows.
\begin{lemma}\label{lem8.11}
 Let $\{u_k\}_k\subset S(a)$ be the $\textnormal{(PS)}$-sequence at level $\gamma_\mu(a,\beta)$ for $\left.\mathcal{J}\right|_{S(a)}$ such that $u_k\rightharpoonup u$ in $HW^{1,Q}(\mathbb{H}^n)$ as $k\to\infty$. Then, it holds that $u\neq 0$ for any $\mu>\mu^\ast$.
\end{lemma}
\begin{proof}
By applying Lemmas \ref{lem8.3} and \ref{lem8.9} under the assumption $\mu>\mu^\ast$, we have
\begin{align*}
\int_{\mathbb{H}^n}\frac{G(u_k)}{r(\xi)^\beta}\,\mathrm{d}\xi\to\int_{\mathbb{H}^n}\frac{G(u)}{r(\xi)^\beta}\,\mathrm{d}\xi\quad\text{and}\quad \int_{\mathbb{H}^n}\frac{g(u_k)u_k}{r(\xi)^\beta}\,\mathrm{d}\xi\to\int_{\mathbb{H}^n}\frac{g(u)u}{r(\xi)^\beta}\,\mathrm{d}\xi
\end{align*}
as $k\to\infty$. To reach a contradiction, we assume that $u=0$. The above convergences, together with \eqref{eq8.7} and Lemma \ref{lem8.10}, yield
\begin{align*}
   \limsup_{k\to\infty} \lambda_k=0\quad\text{and}\quad \lambda_k a^Q=\|D_H u_k\|_Q^Q+o_k(1)\quad\text{as}\quad k\to\infty.
\end{align*}
Thus, we deduce that
\begin{align*}
    0\leq\liminf_{k\to\infty}\|D_H u_k\|_Q^Q\leq\limsup_{k\to\infty}\|D_H u_k\|_Q^Q=a^Q \limsup_{k\to\infty} \lambda_k=0,
\end{align*}
which gives that $\|D_H u_k\|_Q^Q=o_k(1)$ as $k\to\infty$. It follows that $\gamma_\mu(a,\beta)=0$, contrary to $\gamma_\mu(a,\beta) >0$. This shows that $u\neq0$, which completes the proof.
\end{proof}
The following proposition ensures that the nontrivial weak limit of the Palais--Smale sequence inherits the Euler-Lagrange equation, providing a key step toward the existence of a normalized solution and strong convergence.
\begin{proposition}\label{prop8.12}
For $\mu>\mu^\ast$, let $\{u_k\}_k\subset S(a)$ be a bounded $\textnormal{(PS)}$-sequence at level $\gamma_\mu(a,\beta)$ for $\left.\mathcal{J}\right|_{S(a)}$ such that $u_k\rightharpoonup u\neq 0$ in $HW^{1,Q}(\mathbb{H}^n)$ as $k\to\infty$. Furthermore, assume that $\{\lambda_k\}_k\subset\mathbb{R}$ is a sequence satisfying $\lambda_k\to\lambda\,(<0)$ as $k\to\infty$ and there holds
\begin{align}\label{eq9.110}
    -\Delta_{H,Q} u_k-\lambda_k|u_k|^{Q-2}u_k
=\frac{g(u_k)}{r(\xi)^\beta}+o_k(1)
\quad\text{in}\quad\left(HW^{1,Q}(\mathbb{H}^n)\right)^\ast
\end{align}
as $k\to\infty$. Then the weak limit $u$ solves the following Euler-Lagrange equation:
\begin{align}\label{eq9.111}
    -\Delta_{H,Q} u-\lambda|u|^{Q-2}u
=\frac{g(u)}{r(\xi)^\beta}
\quad\text{in}\quad\left(HW^{1,Q}(\mathbb{H}^n)\right)^\ast.
\end{align}
\end{proposition}
\begin{proof}
 Assume that the given hypotheses are satisfied. By applying Lemmas \ref{lem2.1} and \ref{lem4.2}, we get
 \begin{align}\label{eq8.8}
    \begin{cases}
        u_k\to u\quad\text{in}\quad L^p (B_R)\quad\text{for all}\quad R>0\quad\text{and}\quad p\in [1,+\infty),\\
         u_k\to u\quad\text{in}\quad L^q(\mathbb{H}^n,r(\xi)^{-\beta}\mathrm{d}\xi)\quad\text{for all}\quad q\in [Q,+\infty),\\
         u_k\to u\quad\text{a.e. in}\quad \mathbb{H}^n\quad\text{as}\quad k\to\infty.
    \end{cases}
\end{align}
It follows from $\lambda_k\to\lambda$ in $\mathbb{R}$ as $k\to\infty$ that \eqref{eq9.110} is equivalently rewritten as
\begin{align}\label{eq8.9}
    -\Delta_{H,Q} u_k-\lambda|u_k|^{Q-2}u_k
=\frac{g(u_k)}{r(\xi)^\beta}+o_k(1)
\quad\text{in}\quad\left(HW^{1,Q}(\mathbb{H}^n)\right)^\ast\quad\text{as}\quad k\to\infty.
\end{align}
Define the auxiliary Euler-Lagrange functional $\mathcal{L}\colon HW^{1,Q}(\mathbb{H}^n)\to\mathbb{R}$ by
\begin{align*}
    \mathcal{L}(v):=\frac{1}{Q}\big(\|D_H v\|_Q^Q-\lambda \|v\|^Q_Q\big)- \int_{\mathbb{H}^n}\frac{G(v)}{r(\xi)^\beta} \,\mathrm{d}\xi\quad\text{for all}\quad v\in HW^{1,Q}(\mathbb{H}^n).
\end{align*}
Notice that $\mathcal{L}$ is a $C^1$-functional on $HW^{1,Q}(\mathbb{H}^n)$ and its differential satisfies
\begin{align*}
    \langle{\mathcal{L}^\prime(v),w\rangle}=\langle{\mathcal{J}^\prime(v)-\lambda\Psi^\prime(v),w\rangle}\quad\text{for all}\quad v,w\in HW^{1,Q}(\mathbb{H}^n),
\end{align*}
where $\mathcal{J}$ and $\Psi$ are defined in \eqref{eq8.9990} and \eqref{eq8.9999}, respectively. In view of \eqref{eq8.9}, one deduces that $\mathcal{L}^\prime(u_k)\to 0$ in $\left(HW^{1,Q}(\mathbb{H}^n)\right)^\ast$ as $k\to\infty$. Next, we choose $R>0$ and $\psi\in C^\infty_0(\mathbb{H}^n,[0,1])$ such that $\psi\equiv 1$ in $B_R$ and $\psi\equiv 0$ in $\mathbb{H}^n\setminus B_{2R}$. Thus, we obtain from $u_k\rightharpoonup u$ in $HW^{1,Q}(\mathbb{H}^n)$ as $k\to\infty$ that
\begin{align}\label{eq8.10}
    \langle{\mathcal{L}^\prime(u_k)-\mathcal{L}^\prime(u),(u_k-u)\psi}\rangle=o_k(1)\quad\text{as}\quad k\to\infty.
\end{align}
The convexity of the map $z\mapsto\frac{1}{Q}|z|^Q$ guarantees that the map $z\mapsto |z|^{Q-2}z$ is monotone. In addition, for each $k\in\mathbb{N}$, we have
\begin{align*}
    \Theta_k(u_k,u):&=\big(|D_H u_k|_H^{Q-2} D_H u_k-|D_H u|_H^{Q-2} D_H u, D_H u_k-D_H u\big)_H\\
    &\quad-\lambda\big(|u_k|^{Q-2}u_k-|u|^{Q-2}u\big)(u_k-u)\geq 0\quad\text{a.e. in}\quad \mathbb{H}^n.
\end{align*}
Now, by using well-known Simon’s inequality \cite{Simon-1978}, there exists a constant $c_Q>0$ depending only on $Q$ such that
\begin{equation}\label{eq8.11}
\begin{aligned}
&c^{-1}_{Q} \int_{B_R} |D_H u_k-D_H u|_H^{Q}\,\mathrm{d}\xi\\
&\leq c^{-1}_{Q}\int_{B_R} \left(|D_H u_k-D_H u|_H^{Q}-\lambda|u_k-u|^Q\right)\,\mathrm{d}\xi \\ &\leq    \int_{B_R} \Theta_k(u_k,u) \,\mathrm{d}\xi\leq \int_{\mathbb{H}^n}\Theta_k(u_k,u) \psi\,\mathrm{d}\xi \\
&= - \int_{\mathbb{H}^n}\big(|D_H u_k|_H^{Q-2}D_H u_k-|D_H u|_H^{Q-2}D_H u,D_H\psi)_H ( u_k-u)\,\mathrm{d}\xi   \\
&\qquad +\int_{\mathbb{H}^n}\frac{(g(u_k)-g(u))(u_k-u)\psi}{r(\xi)^\beta}\,\mathrm{d}\xi+o_k(1)\quad\text{as }k\to\infty, 
\end{aligned}
\end{equation}
where we have used \eqref{eq8.10} to obtain the last estimate. Due to H\"{o}lder's inequality and \eqref{eq8.8}, we obtain
\begin{align*}
    \bigg|\int_{\mathbb{H}^n}&\big(|D_H u_k|_H^{Q-2}D_H u_k-|D_H u|_H^{Q-2}D_H u,D_H\psi)_H ( u_k-u)\,\mathrm{d}\xi  \bigg|
		\\ & \leq \|D_H \psi\|_\infty\big(\|D_H u_k\|^{Q-1}_Q+\|D_H u\|^{Q-1}_Q\big) \bigg(\int_{B_{2R}}|u_k-u|^Q\,\mathrm{d}\xi\bigg)^{\frac{1}{Q}}\\
        &\lesssim \bigg(\int_{B_{2R}}|u_k-u|^Q\,\mathrm{d}\xi\bigg)^{\frac{1}{Q}}\to 0\quad\text{as }k\to\infty.
\end{align*}
This shows that
\begin{align}\label{eq8.12}
  \lim_{k\to\infty} \int_{\mathbb{H}^n}\big(|D_H u_k|_H^{Q-2}D_H u_k-|D_H u|_H^{Q-2}D_H u,D_H\psi)_H ( u_k-u)\,\mathrm{d}\xi=0.
\end{align}
Set
\begin{align*}
    I:=\varepsilon\int_{B_{2R}}\frac{(|u_k|^{\tau-1}+|u|^{\tau-1})|u
    _k-u|}{r(\xi)^\beta}\,\mathrm{d}\xi
\end{align*}
and \begin{align*}
    II:=C_\varepsilon\int_{B_{2R}}\frac{\big(|u_k|^{q-1} \Phi\big(\alpha\big(1-\frac{\beta}{Q}\big)|u_k|^{Q^\prime}\big)+|u|^{q-1} \Phi\big(\alpha\big(1-\frac{\beta}{Q}\big)|u|^{Q^\prime}\big)\big)|u
    _k-u|}{r(\xi)^\beta}\,\mathrm{d}\xi.
\end{align*}
Take $1<t<\frac{Q}{\beta}$ and $\nu\geq Q$ with $\frac{1}{s}+\frac{1}{s^\prime}=1$ for $s\in\{t,\nu\}$. Then, by applying H\"{o}lder's inequality and \eqref{eq8.8}, we get
\begin{align*}
I& \lesssim\Bigg[\bigg(\int_{B_{2R}}\frac{|u_k|^{(\tau-1)\nu^\prime}}{r(\xi)^\beta}\,\mathrm{d}\xi\bigg)^{\frac{1}{\nu^\prime}}+\bigg(\int_{B_{2R}}\frac{|u|^{(\tau-1)\nu^\prime}}{r(\xi)^\beta}\,\mathrm{d}\xi\bigg)^{\frac{1}{\nu^\prime}}\Bigg]\bigg(\int_{B_{2R}}\frac{|u
_k-u|^\nu}{r(\xi)^\beta}\,\mathrm{d}\xi\bigg)^{\frac{1}{\nu}}\\
&\lesssim \Big(\|u_k\|^{\tau-1}_{L^{(\tau-1)\nu^\prime t^\prime}(B_{2R})}+\|u\|^{\tau-1}_{L^{(\tau-1)\nu^\prime t^\prime}(B_{2R})}\Big)\bigg(\int_{B_{2R}}\frac{1}{r(\xi)^{\beta t}}\,\mathrm{d}\xi\bigg)^{\frac{1}{\nu^\prime t}}\bigg(\int_{B_{2R}}\frac{|u
_k-u|^\nu}{r(\xi)^\beta}\,\mathrm{d}\xi\bigg)^{\frac{1}{\nu}}\\
& \lesssim \bigg(\int_{\mathbb{H}^n}\frac{|u
_k-u|^\nu}{r(\xi)^\beta}\,\mathrm{d}\xi\bigg)^{\frac{1}{\nu}}\to 0\quad\text{as}\quad k\to\infty.
\end{align*}
The condition $\mu>\mu^\ast$ and the fact that $\{u_k\}_k\subset S(a)$, together with Lemma \ref{lem8.9}, gives $ \limsup_{k\to\infty}\|D_H u_k\|_Q^Q<1$. Let $\vartheta,\vartheta^\prime>1$ be such that $\frac{1}{\vartheta}+\frac{1}{\vartheta^\prime}=1$. Fix $m\in(0,1)$ such that $\limsup_{k\to\infty}\|D_H u_k\|_Q^Q<m<1.$ It follows that $\|D_H u_k\|_Q^Q\leq m$ for $k$ large enough. Recall that $m<1$, thus we may choose $\alpha>\alpha_Q$ close to $\alpha_Q$ and $\vartheta>1$ close to $1$ so that $\vartheta\alpha m^{\frac{1}{Q-1}}<\alpha_Q$. Now, by employing H\"{o}lder's inequality, Lemmas \ref{lem3.1}, \ref{lem2.2}, \ref{lem8.1} and \eqref{eq8.8}, we obtain
\begin{align*}
II&\lesssim \Bigg[\bigg(\int_{\mathbb{H}^n}\frac{|u_k|^{(q-1)\vartheta^\prime}|u_k-u|^{\vartheta^\prime}}{r(\xi)^{\beta}}\,\mathrm{d}\xi\bigg)^{\frac{1}{\vartheta^\prime}}\bigg(\int_{\mathbb{H}^n}\frac{\left[\Phi\big(\alpha\big(1-\frac{\beta}{Q}\big)|u_k|^{Q^\prime}\big)\right]^\vartheta}{r(\xi)^{\beta}}\,\mathrm{d}\xi\bigg)^{\frac{1}{\vartheta}}\\
&\quad\quad+\bigg(\int_{\mathbb{H}^n}\frac{|u|^{(q-1)\vartheta^\prime}|u_k-u|^{\vartheta^\prime}}{r(\xi)^{\beta}}\,\mathrm{d}\xi\bigg)^{\frac{1}{\vartheta^\prime}}\bigg(\int_{\mathbb{H}^n}\frac{\Phi\big(\vartheta\alpha\big(1-\frac{\beta}{Q}\big)|u|^{Q^\prime}\big)}{r(\xi)^{\beta}}\,\mathrm{d}\xi\bigg)^{\frac{1}{\vartheta}}\Bigg]\\
&\lesssim \Bigg[\bigg(\int_{\mathbb{H}^n}\frac{|u_k|^{q \vartheta^\prime}}{r(\xi)^{\beta}}\,\mathrm{d}\xi\bigg)^{\frac{q-1}{q \vartheta^\prime}}\bigg(\int_{\mathbb{H}^n}\frac{\left[\Phi\big(\alpha\big(1-\frac{\beta}{Q}\big)|u_k|^{Q^\prime}\big)\right]^\vartheta}{r(\xi)^{\beta}}\,\mathrm{d}\xi\bigg)^{\frac{1}{\vartheta}} \\
&\quad\quad+\bigg(\int_{\mathbb{H}^n}\frac{|u|^{q \vartheta^\prime}}{r(\xi)^{\beta}}\,\mathrm{d}\xi\bigg)^{\frac{q-1}{q \vartheta^\prime}}  \bigg(\int_{\mathbb{H}^n}\frac{\Phi\big(\vartheta\alpha\big(1-\frac{\beta}{Q}\big)|u|^{Q^\prime}\big)}{r(\xi)^{\beta}}\,\mathrm{d}\xi\bigg)^{\frac{1}{\vartheta}}\Bigg]\bigg(\int_{\mathbb{H}^n}\frac{|u_k-u|^{q \vartheta^\prime}}{r(\xi)^{\beta}}\,\mathrm{d}\xi\bigg)^{\frac{1}{q \vartheta^\prime}}\\
&\lesssim \bigg(\int_{\mathbb{H}^n}\frac{|u_k-u|^{q \vartheta^\prime}}{r(\xi)^{\beta}}\,\mathrm{d}\xi\bigg)^{\frac{1}{q \vartheta^\prime}}\to 0\quad\text{as}\quad k\to\infty.
\end{align*}
Taking into account the preceding estimates and \eqref{eq1.5}, we have
\begin{align*}
 \bigg|\int_{\mathbb{H}^n}\frac{(g(u_k)-g(u))(u_k-u)\psi}{r(\xi)^\beta}\,\mathrm{d}\xi\bigg|&\leq \int_{B_{2R}}\cfrac{(|g(u_k)|+|g(u)|)|u_k-u|}{r(\xi)^\beta}\, \mathrm{d}\xi\leq I+II\to 0
\end{align*}
as $k\to\infty$. This yields that
\begin{equation}\label{eq8.13}
  \lim_{k\to\infty} \int_{\mathbb{H}^n}\frac{(g(u_k)-g(u))(u_k-u)\psi}{r(\xi)^\beta}\,\mathrm{d}\xi=0.
\end{equation}
Letting $k\to\infty$ in \eqref{eq8.11}, we obtain from \eqref{eq8.12} and \eqref{eq8.13} that $D_H u_k\to D_H u$ in $L^{Q}(B_{R})$ as $k\to\infty$ for all $R>0$. Hence, up to a subsequence still denoted by itself, we infer that $D_H u_k\to D_H u$ a.e. in $\mathbb{H}^n$ as $k\to\infty$. It follows that $|D_H u_k|_H^{Q-2}D_H u_k \to |D_H u|_H^{Q-2}D_H u$ a.e. in $\mathbb{H}^n$ as $k\to\infty$. This together with the boundedness of $\{|D_H u_k|_H^{Q-2}D_H u_k\}_k$ in $L^{\frac{Q}{Q-1}}(\mathbb{H}^n,\mathbb{R}^{2n})$, yields
\begin{align}\label{eq8.14}
    |D_H u_k|_H^{Q-2}D_H u_k \rightharpoonup |D_H u|_H^{Q-2}D_H u\quad\text{in}\quad L^{\frac{Q}{Q-1}}(\mathbb{H}^n,\mathbb{R}^{2n})\quad\text{as}\quad k\to\infty.
\end{align}
Moreover, \eqref{eq8.8} and the boundedness of $\{|u_k|^{Q-2}u_k\}_k$ in $L^{\frac{Q}{Q-1}}(\mathbb{H}^n)$ imply
\begin{align}\label{eq8.15}
|u_k|^{Q-2}u_k\rightharpoonup |u|^{Q-2}u
\quad\text{in}\quad L^{\frac{Q}{Q-1}}(\mathbb{H}^n)\quad\text{as}\quad k\to\infty.
\end{align}
Fix $v\in C^\infty_0(\mathbb{H}^n)$. Now, by using the density of $C^\infty_0(\mathbb{H}^n)$ in $HW^{1,Q}(\mathbb{H}^n)$, \eqref{eq8.14} and \eqref{eq8.15} give
\begin{equation*}
 \lim_{k\to\infty} \int_{\mathbb{H}^n}|D_H u_k|_H^{Q-2}(D_H u_k,D_H v)_H\,\mathrm{d}\xi=\int_{\mathbb{H}^n}|D_H u|_H^{Q-2}(D_H u,D_H v)_H\,\mathrm{d}\xi
\end{equation*}
and
\begin{equation*}
  \lim_{k\to\infty}\int_{\mathbb{H}^n}|u_k|^{Q-2} u_k v\,\mathrm{d}\xi=\int_{\mathbb{H}^n}|u|^{Q-2} uv\,\mathrm{d}\xi.
\end{equation*}
Moreover, by using Lemma \ref{lem8.3}, we obtain
\begin{align*}
 \sup_{k\in\mathbb{N}}\int_{\mathbb{H}^n}\frac{g(u_k)u_k}{r(\xi)^\beta}\,\mathrm{d}\xi<+\infty.
\end{align*}
Note that $g(u_k) r(\xi)^{-\beta}\in L^1_{\textit{loc}} (\mathbb{H}^n)$. In view of \eqref{eq8.8}, we also have $u_k\to u$ in $L^1 _{\textit{loc}}(\mathbb{H}^n)$ as $k\to\infty$. Hence, by argument analogous to that of  de Souza--do \'{O} \cite[Lemma 3.1]{de Souza-do O-2011}, it follows that $g(u_k) r(\xi)^{-\beta}\to g(u) r(\xi)^{-\beta}$ in $L^1_{\textit{loc}}(\mathbb{H}^n)$ as $k\to\infty$. Recalling $\textnormal{supp}(v)$ is a compact set in $\mathbb{H}^n$, we infer that
\begin{align*}
    \lim_{k\to\infty}\int_{\mathbb{H}^n}\frac{g(u_k)v}{r(\xi)^\beta}\,\mathrm{d}\xi=\int_{\mathbb{H}^n}\frac{g(u)v}{r(\xi)^\beta}\,\mathrm{d}\xi\quad\text{for all}\quad v\in C^\infty_0(\mathbb{H}^n).
\end{align*}
By gathering the above convergences and passing to the limit in \eqref{eq8.9}, we obtain that
\begin{align*}
    \int_{\mathbb{H}^n}|D_H u|_H^{Q-2}(D_H u,D_H v)_H\,\mathrm{d}\xi-\lambda \int_{\mathbb{H}^n}|u|^{Q-2} uv\,\mathrm{d}\xi=\int_{\mathbb{H}^n}\frac{g(u)v}{r(\xi)^\beta}\,\mathrm{d}\xi
\end{align*}
for all $v\in C^\infty_0(\mathbb{H}^n)$. Due to the denseness of $C^\infty_0(\mathbb{H}^n)$ in $HW^{1,Q}(\mathbb{H}^n)$, we infer from the above identity that \eqref{eq9.111} holds. This finishes the proof.
\end{proof}
The following lemma characterizes the unique maximum of $\widetilde{\mathcal{J}}$ along the scaling path and its relation to the Poho\v{z}aev manifold.
\begin{lemma}\label{lem8.13}
The following results hold:
\begin{itemize}
    \item [\textnormal{(i)}] For every $u\in S(a)$, the map $s\mapsto\widetilde{\mathcal{J}}(u,s)$ attains its maximum at some  $s_u\in\mathbb{R}$ and $\mathcal{F}(u,s_u)\in \mathcal{P}(a)$.
    \item [\textnormal{(ii)}] If \eqref{g5} holds, then for every $u\in \mathcal{P}(a)$, one has $\mathcal{J}(u)=\max_{s\in\mathbb R}\mathcal{J}(\mathcal{F}(u,s))$.
\end{itemize}
\end{lemma}
\begin{proof}
 By applying Lemma \ref{lem8.4} and Corollary \ref{cor8.6}, we have $\widetilde{\mathcal{J}}(u,s)\to 0$ as $s\to-\infty$, $\widetilde{\mathcal{J}}(u,s)>0$ for $s$ negative enough and $\widetilde{\mathcal{J}}(u,s)\to -\infty$ as $s\to\infty$. This shows that the map $s\mapsto\widetilde{\mathcal{J}}(u,s)$ attains its global maximum at some $s_u\in\mathbb{R}$, with  $\widetilde{\mathcal{J}}(u,s_u)=\max_{s\in\mathbb{R}} \widetilde{\mathcal{J}}(u,s)$, $\partial_s\widetilde{\mathcal{J}}(u,s_u)=0$ and thus, we obtain $\mathcal{F}(u,s_u)\in \mathcal{P}(a)$. This proves $\textnormal{(i)}$. To prove $\textnormal{(ii)}$, we denote
 \begin{align*}
    \mathcal{M}(\sigma):=\frac{\widetilde{G}(\sigma)}{|\sigma|^{2Q-\beta-1}\sigma}\quad\text{for all}\quad\sigma\in\mathbb{R}\setminus\{0\}.
 \end{align*}
 In virtue of \eqref{g5}, the map $\sigma\mapsto \mathcal{M}(\sigma)$ is strictly increasing on $(-\infty,0)$ and $(0,+\infty)$. If $u(\xi)>0$, then the map $s\mapsto\mathcal{M}(\operatorname{\exp}(s)u(\xi))|u(\xi)|^{2Q-\beta-1}u(\xi)$ is strictly increasing. If $u(\xi)<0$, then $\operatorname{\exp}(s)u(\xi)$ decreases in $s$, so the map $s\mapsto\mathcal{M}(\operatorname{\exp}(s)u(\xi))$ is strictly decreasing, and since $|u(\xi)|^{2Q-\beta-1}u(\xi)<0$, the product is again strictly increasing. To proceed further, a direct calculation gives $$ \partial_s\widetilde{\mathcal{J}}(u,s)=\operatorname{\exp}(Qs)\big(\|D_H u\|_Q^Q-\Pi_u(s)\big),$$ where 
\begin{align*}
\Pi_u(s)= \int_{\mathbb{H}^n}\frac{\mathcal{M}(\operatorname{\exp}(s)u) |u|^{2Q-\beta-1}u}{r(\xi)^\beta}\,\mathrm{d}\xi
\end{align*}
and
\begin{align*}
    \mathcal{M}(\operatorname{\exp}(s)u) |u|^{2Q-\beta-1}u=\frac{\widetilde{G}(\operatorname{\exp}(s)u)}{\operatorname{\exp}((2Q-\beta)s)}.
\end{align*}
If $u\in\mathcal{P}(a)$, then we have $\|D_H u\|_Q^Q=\Pi_u(0)$ and  $\partial_s\widetilde{\mathcal{J}}(u,s)=\operatorname{\exp}(Qs)\left(\Pi_u(0)-\Pi_u(s) \right)$. Notice that the map $s\mapsto \Pi_u(s)$ is strictly increasing on $\mathbb{R}$, thanks to $u\in S(a)$, $a>0$ and the set $\{\xi\in\mathbb{H}^n\colon\, u(\xi)\neq 0\}$ has positive measure. Thus, we infer that $\partial_s\widetilde{\mathcal{J}}(u,s)$ changes its sign from positive to negative as $s$ increases through zero. It follows that $s=0$ is the unique global maximum of the map $s\mapsto\widetilde{\mathcal{J}}(u,s)$. This shows  $\textnormal{(ii)}$, which completes the proof.
\end{proof}
The following result plays a key role in establishing the relation between the mountain pass and Poho\v{z}aev levels.
\begin{lemma}\label{lem8.999}
  Let $u\in S(a)$ and let $K$ be as defined in \eqref{eq8.5}. Then the following statements hold:
  \begin{itemize}
      \item [\textnormal{(i)}] If $\mathcal{J}(u)<0$, then $\mathcal{P}(u)<0$.
      \item [\textnormal{(ii)}] If $\|D_H u\|_Q^Q\leq K$, then $\mathcal{P}(u)\geq \big(1-\frac{1}{4Q}\big)\|D_H u\|_Q^Q>0$.
  \end{itemize}
\end{lemma}
\begin{proof}
Suppose that $\mathcal{J}(u)<0$. Next, we claim that $\mathcal{P}(u)<0$. Indeed, if not, let $\mathcal{P}(u)\geq 0$. By using  $\theta>2Q-\beta$ and $G(s)\geq 0$, we obtain from \eqref{eq9.999} and \eqref{g3} that
\begin{align*}
   \mathcal{J}(u)
   &\geq \frac{1}{Q} \int_{\mathbb{H}^n}\frac{g(u)u-(2Q-\beta)G(u)}{r(\xi)^\beta}\,\mathrm{d}\xi\geq \left(\frac{\theta-2Q+\beta}{Q}\right)\int_{\mathbb{H}^n}\frac{G(u)}{r(\xi)^\beta} \,\mathrm{d}\xi\geq 0,
\end{align*}
which is a contradiction. This proves the claim and hence verifies $\textnormal{(i)}$. To prove $\textnormal{(ii)}$, assume that $\|D_H u\|_Q^Q\leq K$. In view of \eqref{eq1.6} and arguing similarly as in Lemma \ref{lem8.5}, we deduce that
\begin{align*}
\int_{\mathbb{H}^n}\frac{|g(u)u|}{r(\xi)^\beta}\,\mathrm{d}\xi&\leq C_1\big(\|D_Hu\|_Q^Q\big)^{\frac{\tau-Q+\beta}{Q}} a^{Q-\beta}+C_2\big(\|D_Hu\|_Q^Q\big)^{\frac{qt^\prime-Q+\beta}{Qt^\prime}} a^{\frac{Q-\beta}{t^\prime}},
\end{align*}
where $C_1,C_2>0$ are the same constants as in \eqref{eq8.3}. Moreover, by using \eqref{eq8.5} and proceeding exactly as in the proof of Corollary \ref{cor8.6}, we infer from the above estimate that
\begin{align*}
  \int_{\mathbb{H}^n}\frac{|g(u)u|}{r(\xi)^\beta}\,\mathrm{d}\xi\leq\frac{\|D_Hu\|_Q^Q}{4Q}
\end{align*}
for all $u\in S(a)$ with $\|D_H u\|_Q^Q\leq K$.
It follows from $a>0$ and $G(s)\geq 0$ that
\begin{align*}
    \mathcal{P}(u)\geq \|D_Hu\|_Q^Q- \int_{\mathbb{H}^n}\frac{g(u)u}{r(\xi)^\beta} \,\mathrm{d}\xi\geq\Big(1-\frac{1}{4Q}\Big)\|D_H u\|_Q^Q>0.
\end{align*}
This establishes $\textnormal{(ii)}$, which completes the proof.
\end{proof}
Next, the mountain pass level and the Poho\v{z}aev level are shown to coincide, as stated below.
\begin{lemma}\label{lem8.14}
 Assume that \eqref{g5} holds. Then the Poho\v{z}aev level $m(a)$ is characterized by $m(a)=\inf_{u\in \mathcal{P}(a)}\max_{s\in\mathbb{R}}\mathcal{J}(\mathcal{F}(u,s)).$
 Moreover, there holds $\gamma_{\mu}(a,\beta)=m(a)$.
\end{lemma}
\begin{proof}
Take $u\in \mathcal{P}(a)$. Then it follows from Lemma \ref{lem8.13} that $
\mathcal{J}(u)=\max_{s\in\mathbb R}\mathcal{J}(\mathcal{F}(u,s))
$, which in turn yields
\begin{align}\label{eq8.333}
 m(a)
=\inf_{u\in \mathcal{P}(a)}\mathcal{J}(u)
=\inf_{u\in \mathcal{P}(a)}\max_{s\in\mathbb R}\mathcal{J}(\mathcal{F}(u,s)).
\end{align}
Now, we fix $u\in S(a)$. By employing Lemma \ref{lem8.4}, we can find $s_1<0$ and $s_2>0$ such that $\|D_H \mathcal{F}(u,s_1)\|_Q^Q<\frac{K}{2}$ and $\mathcal{J}(\mathcal{F}(u,s_2))<0$.
Define the path $\psi\colon[0,1]\to S(a)$ by $\psi(t)=\mathcal{F}(u,(1-t)s_1+ts_2)$ for all $t\in[0,1]$. It follows that $\psi\in\Gamma$ and thus, $\Gamma\neq\emptyset$. Therefore, by the definition of  $ \gamma_{\mu}(a,\beta)$, one has
\begin{align*}
  \gamma_{\mu}(a,\beta)\leq \max_{t\in[0,1]} \mathcal{J}(\psi(t))\leq\max_{s\in\mathbb{R}} \mathcal{J}(\mathcal{F}(u,s)).
\end{align*}
In particular, for any $u\in \mathcal{P}(a)$, we obtain from \eqref{eq8.333} that
\begin{align*}
  \gamma_{\mu}(a,\beta)\leq\inf_{u\in \mathcal{P}(a)}\max_{s\in\mathbb{R}} \mathcal{J}(\mathcal{F}(u,s))= m(a).
\end{align*}
Conversely, for any $\psi\in\Gamma$, we have $\|D_H \psi(0)\|_Q^Q<\frac{K}{2}$ and $\mathcal{J}(\psi(1))<0$. By applying Lemma \ref{lem8.999}, there holds $\mathcal{P}(\psi(0))>0>\mathcal{P}(\psi(1))$. It is easy to see that $\mathcal{P}$ is continuous on $HW^{1,Q}(\mathbb{H}^n)$. Thus, there exists some $t_0\in(0,1)$ such that $\mathcal{P}(\psi(t_0))=0$. In particular, we have $\psi(t_0)\in\mathcal{P}(a)$ and
\begin{align*}
    \max_{t\in[0,1]} \mathcal{J}(\psi(t))\geq \mathcal{J}(\psi(t_0))\geq m(a).
\end{align*}
Hence, we get $\gamma_{\mu}(a,\beta)\geq m(a)$. Thus, the desired assertion follows, which completes the proof.
\end{proof}
Finally, we conclude the discussion by establishing the existence of nontrivial ground-state solutions to our main problem.
\subsection{{Proof of Theorem \ref{thm1.11}}}Fix $\mu>\mu^\ast$. Let $\{u_k\}_k\subset S(a)$ be a $\textnormal{(PS)}$-sequence at level $\gamma_\mu(a,\beta)$ for $\left.\mathcal{J}\right|_{S(a)}$ as described in the above analysis. The boundedness of the sequence $\{u_k\}_k$ in $HW^{1,Q}(\mathbb{H}^n)$ implies that there exists some $u_a\in HW^{1,Q}(\mathbb{H}^n)$ such that $u_k\rightharpoonup u_a$ in $HW^{1,Q}(\mathbb{H}^n)$ as $k\to\infty$. It follows from Lemma \ref{lem8.11} that $u_a\neq0$. Now, by using Lemmas \ref{lem8.3} and \ref{lem8.10}, we obtain from
\begin{align*}
 \limsup_{k\to\infty}\lambda_k= -\bigg(\frac{Q-\beta}{a^Q}\bigg)\liminf_{k\to\infty}\int_{\mathbb{H}^n}\frac{G(u_k)}{r(\xi)^\beta}\,\mathrm{d}\xi=-\bigg(\frac{Q-\beta}{a^Q}\bigg)\int_{\mathbb{H}^n}\frac{G(u_a)}{r(\xi)^\beta}\,\mathrm{d}\xi<0.
\end{align*}
This shows that, passing to a subsequence if necessary, we may assume  $\lambda_k\to\hat\lambda_a<0$ in $\mathbb{R}$ as $k\to\infty$. Since the sequences $\{u_k\}_k$ and $\{\lambda_k\}_k$ satisfy all the hypotheses of Proposition \ref{prop8.12}, we conclude that $u_a$ weakly solves
\begin{align*}
    -\Delta_{H,Q} u_a-\hat\lambda_a|u_a|^{Q-2}u_a
=\frac{g(u_a)}{r(\xi)^\beta}
\quad\text{in}\quad\left(HW^{1,Q}(\mathbb{H}^n)\right)^\ast.
\end{align*}
In particular, it holds
\begin{align}\label{eq8.16}
    \|D_H u_a\|_Q^Q-\hat\lambda_a \| u_a\|_Q^Q=\int_{\mathbb{H}^n}\frac{g(u_a)u_a}{r(\xi)^\beta}\,\mathrm{d}\xi.
\end{align}
Recall that
\begin{align*}
    \|D_H u_k\|_Q^Q-\lambda_k \| u_k\|_Q^Q=\int_{\mathbb{H}^n}\frac{g(u_k)u_k}{r(\xi)^\beta}\,\mathrm{d}\xi+o_k(1)\quad\text{as}\quad k\to\infty,
\end{align*}
which is equivalent to
\begin{align}\label{eq8.17}
    \|D_H u_k\|_Q^Q-\hat\lambda_a \| u_k\|_Q^Q=\int_{\mathbb{H}^n}\frac{g(u_k)u_k}{r(\xi)^\beta}\,\mathrm{d}\xi+o_k(1)\quad\text{as}\quad k\to\infty.
\end{align}
By applying Lemma \ref{lem8.3}, we obtain from \eqref{eq8.16} and \eqref{eq8.17} that
\begin{align*}
    \lim_{k\to\infty}\big(\|D_H u_k\|_Q^Q-\hat\lambda_a \| u_k\|_Q^Q\big)=\|D_H u_a\|_Q^Q-\hat\lambda_a \| u_a\|_Q^Q.
\end{align*}
The condition $\hat\lambda_a<0$ allows us to equip $HW^{1,Q}(\mathbb{H}^n)$ with the equivalent norm $\|u_a\|_{\lambda_a}=\big(\|D_H u_a\|_Q^Q-\hat\lambda_a \| u_a\|_Q^Q\big)^{\frac{1}{Q}}$, where $\lambda_a=-\hat\lambda_a$. The uniform convexity of $HW^{1,Q}(\mathbb{H}^n)$, together with the preceding convergence, gives $u_k\to u_a$ in $HW^{1,Q}(\mathbb{H}^n)$ as $k\to\infty$. This proves that $u_a\in S(a)$ is a nontrivial solution of \eqref{main problem} corresponding to the parameter  $\lambda=\lambda_a=-\hat\lambda_a>0$. It follows from $u_k\to u_a$ in $HW^{1,Q}(\mathbb{H}^n)$ as $k\to\infty$ and Lemma \ref{lem8.3} that $\mathcal{P}(u_a)=\lim_{k\to\infty}\mathcal{P}(u_k)=0$ and $ \mathcal{J}(u_a)=\lim_{k\to\infty}\mathcal{J}(u_k)=\gamma_{\mu}(a,\beta)$. Thus, we infer that $u_a\in\mathcal{P}(a)$. If, in addition, \eqref{g5} holds, then Lemma \ref{lem8.14} gives $\gamma_{\mu}(a,\beta)=m(a)$, and hence $\mathcal{J}(u_a)=m(a)$. This shows that $u_a$ is a nontrivial ground state solution of \eqref{main problem} on $\mathcal{P}(a)$. Hence, the proof is complete. \qed

\section{The Scope and Prospects of Open Questions}\label{sec9}
The results of this work raise several questions concerning alternative approaches to the existence of maximizers for weighted Moser--Trudinger inequalities, the attainment of the corresponding nonsingular inequality, and the role of singular weights in the symmetry and stability of extremals for Caffarelli--Kohn--Nirenberg interpolation-type inequalities in the borderline Sobolev regime on $\mathbb{H}^n$. It is also natural to ask whether a suitable growth condition on $g$ at infinity can replace the requirement $\mu>\mu^\ast$ in Theorem \ref{thm1.11}. This discussion thus leaves open the following questions:
\begin{itemize}
\item[$\diamond$] Is there an alternative approach to proving the existence of maximizers for the weighted Moser--Trudinger inequality on $\mathbb{H}^n$?
\item[$\diamond$] Does the Moser--Trudinger inequality on $\mathbb{H}^n$ admit a maximizer in the nonsingular  case ($\beta=0$)?
\item[$\diamond$] Do Caffarelli--Kohn--Nirenberg interpolation-type inequalities on $\mathbb{H}^n$ exhibit symmetry breaking and admit quantitative stability estimates for their extremals in the borderline Sobolev regime?
\item[$\diamond$] The condition $\mu>\mu^\ast$ in \eqref{g4} of Theorem \ref{thm1.11} is essential to our argument. In the Euclidean setting, assumptions of this type are often replaced by a suitable growth condition at infinity, namely, 
\[
\liminf_{|s|\to\infty}
\frac{sg(s)}{\exp\big(\alpha_Q(1-\frac{\beta}{Q})|s|^{Q^\prime}\big)}>0~\text{(large enough)},
\]
together with an estimate of the mountain pass level using Moser functions, as in \eqref{eq4.6}. Can this approach be adapted to replace the condition $\mu>\mu^\ast$ in \eqref{g4} by a suitable growth condition on $g$ at infinity?
\end{itemize}
\medskip
\noindent\textbf{Acknowledgements.} S. Kumar acknowledges support from the CSIR PhD Fellowship, Sanction No. 09/1125(25678)/2025-EMR-I, Government of India. N. Lam acknowledges support from the NSERC Discovery Grant RGPIN-2021-03584, Canada. D.K. Mahanta acknowledges support from the Academic Research Associate position in the School of Mathematics, Indian Institute of Science Education and Research Thiruvananthapuram, India. A. Sarkar acknowledges support from the DST-INSPIRE Grant No. DST/INSPIRE/04/2018/002208, Government of India.

\medskip
\noindent\textbf{Data availability.} No data was used for the research described in the article.

\medskip
\noindent\textbf{Conflict of interest.} The authors declare that they have no conflict of interest.



\begin{thebibliography}{99}
\bibitem{Adams-1988}
D.R. Adams,
{\it A sharp inequality of J. Moser for higher order derivatives},
Ann. of Math. (2) {\bf 128} (1988), no. 2, 385--398. DoI: \url{https://doi.org/10.2307/1971445}



\bibitem{Adi-Druet2004}
Adimurthi, O. Druet,
\emph{Blow-up analysis in dimension 2 and a sharp form of Trudinger-Moser inequality}, Comm. Partial Differential Equations \textbf{29} (2004), no. 1--2, 295--322. DoI: \url{https://doi.org/10.1081/PDE-120028854}



\bibitem{Adi-Sandeep-2007}
Adimurthi, K. Sandeep,
{\it A singular Moser-Trudinger embedding and its applications},
NoDEA Nonlinear Differential Equations Appl. {\bf 13} (2007), no. 5--6, 585--603. DoI: \url{https://doi.org/10.1007/s00030-006-4025-9}



\bibitem{Adi-Yang-2010}
Adimurthi, Y. Yang,
{\it An interpolation of Hardy inequality and Trudinger-Moser inequality in $\mathbb{R}^N$ and its applications},
Int. Math. Res. Not. IMRN {\bf 2010} (2010), no. 13, 2394--2426. DoI: \url{https://doi.org/10.1093/imrn/rnp194}



\bibitem{AJM-2022}
C.O. Alves, C. Ji, O.H. Miyagaki,
{\it Normalized solutions for a Schr\"{o}dinger equation with
critical growth in $\mathbb{R}^N$},
Calc. Var. Partial Differential Equations {\bf 61} (2022), no. 1, Paper No. 18. DoI: \url{https://doi.org/10.1007/s00526-021-02123-1}



\bibitem{BMT-2003}
Z.M. Balogh, J.J. Manfredi, J.T. Tyson,
{\it Fundamental solution for the $Q$-Laplacian and sharp Moser-Trudinger inequality in Carnot groups},
J. Funct. Anal. {\bf 204} (2003), no. 1, 35--49. DoI: \url{https://doi.org/10.1016/S0022-1236(02)00169-6}



\bibitem{BS-2017}
T. Bartsch, N. Soave,
{\it A natural constraint approach to normalized solutions of
nonlinear Schr\"{o}dinger equations and systems},
J. Funct. Anal. {\bf 272} (2017), no. 12, 4998--5037. DoI: \url{https://doi.org/10.1016/j.jfa.2017.01.025}



\bibitem{CKN-1984}
L. Caffarelli, R. Kohn, L. Nirenberg,
{\it First order interpolation inequalities with weights},
Compositio Math. {\bf 53} (1984), no. 3, 259--275. DoI: \url{http://www.numdam.org/item?id=CM_1984__53_3_259_0}



\bibitem{Cao-1992}
D.M. Cao,
{\it Nontrivial solution of semilinear elliptic equation with critical exponent in $\mathbb{R}^2$},
Comm. Partial Differential Equations {\bf 17} (1992), no. 3--4, 407--435. DoI: \url{https://doi.org/10.1080/03605309208820848}



\bibitem{CarlesonChang1986}
L. Carleson, S.-Y. A. Chang,
\emph{On the existence of an extremal function for an inequality of J. Moser}, Bull. Sci. Math. (2) \textbf{110} (1986), no. 2, 113--127.



\bibitem{CHTZ-2026}
D. Cassani, L. Huang, C. Tarsi, X. Zhong,
{\it The mass-mixed case for normalized solutions to NLS
equations in dimension two},
J. Anal. Math. {\bf 159} (2026), no. 1, 359--400. DoI: \url{https://doi.org/10.1007/s11854-026-0436-1}



\bibitem{Cassani-Sani-Tarsi-2014}
D. Cassani, F. Sani, C. Tarsi,
{\it Equivalent Moser type inequalities in $\mathbb{R}^2$ and the zero mass case},
J. Funct. Anal. {\bf 267} (2014), no. 11, 4236--4263. DoI: \url{https://doi.org/10.1016/j.jfa.2014.09.022}



\bibitem{Chang-Yang-1987}
S.-Y.A. Chang, P. Yang,
\emph{Prescribing Gaussian curvature on $S^2$}, Acta Math. \textbf{159} (1987), no. 3--4, 215--259. DoI: \url{https://doi.org/10.1007/BF02392560}



\bibitem{Chang-Yang2003}
S.-Y.A. Chang, P. Yang,
\emph{The inequality of Moser and Trudinger and applications to conformal geometry}, Comm. Pure Appl. Math. \textbf{56} (2003), no. 8, 1135--1150. DoI: \url{https://doi.org/10.1002/cpa.3029}



\bibitem{Chen-Lu-Zhu-2021}
L. Chen, G. Lu, M. Zhu,
{\it Sharpened Trudinger-Moser inequalities on the Euclidean
space and Heisenberg group},
J. Geom. Anal. {\bf 31} (2021), no. 12, 12155--12181. DoI: \url{https://doi.org/10.1007/s12220-021-00713-1}



\bibitem{Chen-Lu-Zhu-2023}
L. Chen, G. Lu, M. Zhu,
{\it Least energy solutions to quasilinear subelliptic equations with constant and degenerate potentials on the Heisenberg group},
Proc. Lond. Math. Soc. (3) {\bf 126} (2023), no. 2, 518--555. DoI: \url{https://doi.org/10.1112/plms.12495}





\bibitem{CRTY-2023}
S. Chen, V.D. R\u{a}dulescu, X. Tang, S. Yuan,
{\it Normalized solutions for Schr\"{o}dinger equations with
critical exponential growth in $\mathbb{R}^2$},
SIAM J. Math. Anal. {\bf 55} (2023), no. 6, 7704--7740. DoI: \url{https://doi.org/10.1137/22M1521675}





\bibitem{Cohn-Lam-Lu-Yang-2012}
W.S. Cohn, N. Lam, G. Lu, Y. Yang,
{\it The Moser-Trudinger inequality in unbounded domains of Heisenberg group and sub-elliptic equations},
Nonlinear Anal. {\bf 75} (2012), no. 12, 4483--4495. DoI: \url{https://doi.org/10.1016/j.na.2011.09.053}



\bibitem{Cohn-Lu-2001}
W.S. Cohn, G. Lu,
{\it Best constants for Moser-Trudinger inequalities on the Heisenberg group},
Indiana Univ. Math. J. {\bf 50} (2001), no. 4, 1567--1591. DoI: \url{https://doi.org/10.1512/iumj.2001.50.2138}



\bibitem{Cohn-Lu-2004}
W. Cohn, G. Lu,
\emph{Sharp constants for Moser-Trudinger inequalities on spheres in complex space $\mathbb{C}^n$}, Comm. Pure
Appl. Math. {\bf 57} (2004), no. 11, 1458--1493. DoI: \url{https://doi.org/10.1002/cpa.20043}



\bibitem{delPino-Musso-Ruf-2010}
M. del Pino, M. Musso, B. Ruf,
{\it New solutions for Trudinger-Moser critical equations in $\mathbb{R}^2$},
J. Funct. Anal. {\bf 258} (2010), no. 2, 421--457. DoI: \url{https://doi.org/10.1016/j.jfa.2009.06.018}

\bibitem{delPino-Musso-Ruf-2012}
M. del Pino, M. Musso, B. Ruf,
{\it Beyond the Trudinger-Moser supremum},
Calc. Var. Partial Differential Equations {\bf 44} (2012), no. 3-4, 543--576. DoI: \url{https://doi.org/10.1007/s00526-011-0444-5}

\bibitem{DMMT-2022}
F. De Marchis, A. Malchiodi, L. Martinazzi, P.-D. Thizy,
{\it Critical points of the Moser-Trudinger functional on
closed surfaces},
Invent. Math. {\bf 230} (2022), no. 3, 1165--1248. DoI: \url{https://doi.org/10.1007/s00222-022-01142-9}



\bibitem{DHZZ-2024}
S. Deng, L. Huang, J. Zhang, X. Zhong,
{\it Normalized solutions for a coupled Schr\"{o}dinger systems in
$\mathbb{R}^2$ with exponential critical growth: mass
super-critical case},
J. Differential Equations {\bf 411} (2024), 349--380. DoI: \url{https://doi.org/10.1016/j.jde.2024.08.015}



\bibitem{de Souza-do O-2011}
M. de Souza, J.M. do \'{O},
{\it On a singular and nonhomogeneous $N$-Laplacian equation involving critical growth},
J. Math. Anal. Appl. {\bf 380} (2011), no. 1, 241--263. DoI: \url{https://doi.org/10.1016/j.jmaa.2011.03.028}



\bibitem{do O-1997}
J.M. do \'O,
{\it $N$-Laplacian equations in $\mathbb{R}^N$ with critical growth},
Abstr. Appl. Anal. {\bf 2} (1997), no. 3--4, 301--315. DoI: \url{https://doi.org/10.1155/S1085337597000419}



\bibitem{Figueiredo-Mahanta-Sarkar-2027}
G.M. Figueiredo, D.K. Mahanta, A. Sarkar,
{\it On sharp singular Moser-Trudinger type inequalities and applications to zero mass $(p,Q)$-Laplace subelliptic equations in the Heisenberg group},
J. Differential Equations {\bf 489} (2027), Paper No. 114850. DoI: \url{https://doi.org/10.1016/j.jde.2026.114850}




\bibitem{Flucher1992}
M. Flucher,
\emph{Extremal functions for the Trudinger-Moser inequality in $2$ dimensions}, Comment. Math. Helv. \textbf{67} (1992), no. 3, 471--497. DoI: \url{https://doi.org/10.1007/BF02566514}



\bibitem{Folland-Stein}
G.B. Folland, E.M. Stein,
\emph{Estimates for the $\bar\partial_b$ complex and analysis on the Heisenberg group}, Comm. Pure Appl. Math. \textbf{27} (1974), 429--522. DoI: \url{https://doi.org/10.1002/cpa.3160270403}



\bibitem{FS-82}
G.B. Folland, E.M. Stein,
\emph{Hardy Spaces on Homogeneous Groups}, Mathematical Notes {\bf 28}, Princeton University Press, Princeton, NJ, 1982. DoI: \url{https://doi.org/10.2307/j.ctv17db3q0}



\bibitem{Changfeng-Moradifam-2018}
C. Gui, A. Moradifam,
{\it The sphere covering inequality and its applications},
Invent. Math. {\bf 214} (2018), no. 3, 1169--1204. DoI: \url{https://doi.org/10.1007/s00222-018-0820-2}



\bibitem{Hormander-1967}
L. H\"{o}rmander,
{\it Hypoelliptic second order differential equations},
Acta Math. {\bf 119} (1967), 147--171. DoI: \url{https://doi.org/10.1007/BF02392081}



\bibitem{Ikoma-Ishiwata-Wadade-2019}
N. Ikoma, M. Ishiwata, H. Wadade,
{\it Existence and non-existence of maximizers for the
Moser-Trudinger type inequalities under inhomogeneous
constraints},
Math. Ann. {\bf 373} (2019), no. 1--2, 831--851. DoI: \url{https://doi.org/10.1007/s00208-018-1709-5}



\bibitem{Ishiwata-2011}
M. Ishiwata,
{\it Existence and nonexistence of maximizers for variational problems associated with Trudinger-Moser type inequalities in $\mathbb{R}^N$},
Math. Ann. {\bf 351} (2011), no. 4, 781--804. DoI: \url{https://doi.org/10.1007/s00208-010-0618-z}



\bibitem{Ishiwata-Nakamura-Wadade-2014}
M. Ishiwata, M. Nakamura, H. Wadade,
{\it On the sharp constant for the weighted Trudinger-Moser type inequality of the scaling invariant form},
Ann. Inst. H. Poincar\'e{} C Anal. Non Lin\'eaire {\bf 31} (2014), no. 2, 297--314. DoI: \url{https://doi.org/10.1016/j.anihpc.2013.03.004}



\bibitem{Jeanjean-1997}
L. Jeanjean,
{\it Existence of solutions with prescribed norm for semilinear elliptic equations},
Nonlinear Anal. {\bf 28} (1997), no. 10, 1633--1659. DoI: \url{https://doi.org/10.1016/S0362-546X(96)00021-1}




\bibitem{JL1988}
D. Jerison, J.M. Lee,
{\it Extremals for the Sobolev inequality on the Heisenberg group and the CR Yamabe problem},
 J. Amer. Math. Soc. {\bf 1} (1988), no. 1, 1--13. DoI: \url{https://doi.org/10.2307/1990964}



\bibitem{Lam-2017}
N. Lam,
{\it Maximizers for the singular Trudinger-Moser inequalities in the subcritical cases},
 Proc. Amer. Math. Soc. \textbf{145} (2017), no. 11, 4885--4892. DoI: \url{https://doi.org/10.1090/proc/13624}



\bibitem{Lam-Lu-2012}
N. Lam, G. Lu,
{\it Sharp Moser-Trudinger inequality on the Heisenberg group at the critical case and applications},
Adv. Math. {\bf 231} (2012), no. 6, 3259--3287. DoI: \url{https://doi.org/10.1016/j.aim.2012.09.004}



\bibitem{Lam-Lu-2012@@}
N. Lam, G. Lu,
{\it Sharp Adams type inequalities in Sobolev spaces $W^{m,\frac{n}{m}}(\mathbb{R}^n)$ for arbitrary integer $m$},
J. Differential Equations {\bf 253} (2012), no. 4, 1143--1171. DoI: \url{https://doi.org/10.1016/j.jde.2012.04.025}



\bibitem{Lam-Lu-2013}
N. Lam, G. Lu,
\emph{A new approach to sharp Moser-Trudinger and Adams type inequalities: a rearrangement-free argument}, J. Differential Equations {\bf 255} (2013), no. 3, 298--325. DoI: \url{https://doi.org/10.1016/j.jde.2013.04.005}



\bibitem{Lam-Lu-Tang-2012}
N. Lam, G. Lu, H. Tang,
{\it On nonuniformly subelliptic equations of $Q$-sub-Laplacian type with critical growth in the Heisenberg group},
Adv. Nonlinear Stud. {\bf 12} (2012), no. 3, 659--681. DoI: \url{https://doi.org/10.1515/ans-2012-0312}



\bibitem{Lam-Lu-Tang-2014}
N. Lam, G. Lu, H. Tang,
{\it Sharp subcritical Moser-Trudinger inequalities on Heisenberg groups and subelliptic PDEs},
Nonlinear Anal. {\bf 95} (2014), 77--92. DoI: \url{https://doi.org/10.1016/j.na.2013.08.031}



\bibitem{Lam-Lu-Zhang-2017}
N. Lam, G. Lu, L. Zhang,
{\it Equivalence of critical and subcritical sharp
Trudinger-Moser-Adams inequalities},
Rev. Mat. Iberoam. {\bf 33} (2017), no. 4, 1219--1246. DoI: \url{https://doi.org/10.4171/RMI/969}



\bibitem{Lam-Lu-Zhang-2019}
N. Lam, G. Lu, L. Zhang,
\emph{Existence and nonexistence of extremal functions for sharp Trudinger-Moser inequalities},
Adv. Math. \textbf{352} (2019), 1253--1298. DoI: \url{https://doi.org/10.1016/j.aim.2019.06.020}



\bibitem{Lam-Tang-2013}
N. Lam, H. Tang,
{\it Sharp constants for weighted Moser-Trudinger inequalities on groups of Heisenberg type},
Nonlinear Anal. {\bf 89} (2013), 95--109. DoI: \url{https://doi.org/10.1016/j.na.2013.04.017}



\bibitem{Li2001}
Y. Li,
\emph{Moser-Trudinger inequality on compact Riemannian manifolds of dimension two}, J. Partial Differential Equations \textbf{14} (2001), no. 2, 163--192. DoI: \url{https://doi.org/10.4208/jpde.v14.n2.7}



\bibitem{Li-Lu-Zhu-2018}
J. Li, G. Lu, M. Zhu,
{\it Concentration-compactness principle for Trudinger-Moser inequalities on Heisenberg groups and existence of ground state solutions},
Calc. Var. Partial Differential Equations {\bf 57} (2018), no. 3, Paper No. 84. DoI: \url{https://doi.org/10.1007/s00526-018-1352-8}



\bibitem{Li-Lu-Zhu-2021}
J. Li, G. Lu, M. Zhu,
{\it Concentration-compactness principle for Trudinger-Moser’s inequalities on Riemannian manifolds and Heisenberg groups: a completely symmetrization-free argument},
Adv. Nonlinear Stud. {\bf 21} (2021), no. 4, 917--937. DoI: \url{https://doi.org/10.1515/ans-2021-2147}



\bibitem{Lions-1985}
P.-L. Lions,
{\it The Concentration-Compactness Principle in the Calculus of Variations. The limit case. I},
Rev. Mat. Iberoam. {\bf 1} (1985), no. 1, 145--201. DoI: \url{https://doi.org/10.4171/RMI/6}





\bibitem{DTAL-2026}
D.K. Mahanta, T. Mukherjee, A. Sarkar, L. Sharma,
{\it On the study of $(p,Q)$-Laplace Choquard equations
with critical Trudinger-Moser nonlinearity in $\mathbb{H}^N$},
Anal. Appl. (Singap.) {\bf 24} (2026), no. 7, 1807--1847. DoI: \url{https://doi.org/10.1142/S021953052650003X}


\bibitem{JJM2019}
J.J. Manfredi, V.N. Vera de Serio,
{\it Rearrangements in Carnot Groups},
Acta Math. Sin. (Engl. Ser.) {\bf 35} (2019), no. 7, 1115--1127. DoI: \url{https://doi.org/10.1007/s10114-019-8271-2}


\bibitem{Moser-1970/71}
J. Moser,
{\it A sharp form of an inequality by N. Trudinger},
Indiana Univ. Math. J. {\bf 20} (1970/71), 1077--1092. DoI: \url{https://www.jstor.org/stable/24890183}



\bibitem{Nguyen-2021}
V.H. Nguyen,
{\it The thresholds of the existence of maximizers for the critical sharp singular Moser-Trudinger inequality under constraints},
Math. Ann. {\bf 380} (2021), no. 3--4, 1933--1958. DoI: \url{https://doi.org/10.1007/s00208-020-02010-8}



\bibitem{Nguyen-2021@}
V.H. Nguyen,
{\it Extremal functions for sharp Moser-Trudinger type
inequalities in the whole space $\mathbb{R}^N$},
J. Funct. Anal. {\bf 280} (2021), no. 3, Paper No. 108833. DoI: \url{https://doi.org/10.1016/j.jfa.2020.108833}




\bibitem{Ozawa-Ruzhansky-Suragan-2019}
T. Ozawa, M. Ruzhansky, D. Suragan,
{\it $L^p$-Caffarelli-Kohn-Nirenberg type inequalities on homogeneous groups},
Q. J. Math. {\bf 70} (2019), no. 1, 305--318. DoI: \url{https://doi.org/10.1093/qmath/hay040}




\bibitem{Pokhozhaev-1965}
S.I. Poho\v{z}aev,
{\it The Sobolev embedding in the case $pl=n$},
in: Proceedings of the Technical Scientific Conference on Advances of Scientific Research 1964--1965, Mathematics Section, Moskov. Energet. Inst., Moscow, 1965, pp. 158--170.

\bibitem{PT-2022}
P. Pucci, L. Temperini,
{\it Existence for singular critical exponential $(p,Q)$
 equations in the Heisenberg group},
Adv. Calc. Var. {\bf 15} (2022), no. 3, 601--617. DoI: \url{https://doi.org/10.1515/acv-2020-0028}



\bibitem{Rockafellar-1970}
R.T. Rockafellar,
\emph{Convex analysis},
Princeton University Press, Princeton, 1970.



\bibitem{Royden-Fitzpatrick-1988}
H. L. Royden, P. Fitzpatrick, \emph{Real analysis},
Fourth Edition, Pearson, Boston, 2010.


\bibitem{Ruf-2005}
B. Ruf,
{\it A sharp Trudinger-Moser type inequality for unbounded domains in $\mathbb{R}^2$},
J. Funct. Anal. {\bf 219} (2005), no. 2, 340--367. DoI: \url{https://doi.org/10.1016/j.jfa.2004.06.013}



\bibitem{Ruf-Sani-2013}
B. Ruf, F. Sani,
{\it Sharp Adams-type inequalities in $\mathbb{R}^n$},
Trans. Amer. Math. Soc. {\bf 365} (2013), no. 2, 645--670. DoI: \url{https://doi.org/10.1090/S0002-9947-2012-05561-9}



\bibitem{Simon-1978}
J. Simon,
\emph{ R\'{e}gularit\'{e} de la solution d'une \'{e}quation non lin\'{e}aire dans $\mathbb{R}^N$},
in: Journ\'{e}es d'Analyse Non Lin\'{e}aire (Besan\c{c}on, 1977), Lecture Notes in Math. {\bf 665}, Springer, Berlin, 1978, pp. 205--227. DoI: \url{https://doi.org/10.1007/BFb0061807}





\bibitem{Soave-2020}
N. Soave,
{\it Normalized ground states for the NLS equation with combined nonlinearities: the Sobolev critical case},
J. Funct. Anal. {\bf 279} (2020), no. 6, Paper No. 108610. DoI: \url{https://doi.org/10.1016/j.jfa.2020.108610}



\bibitem{Struwe1988}
M. Struwe,
{\it Critical points of embeddings of $H_0^{1,n}$ into Orlicz spaces}, Ann. Inst. H. Poincar\'{e} Anal. Non Lin\'{e}aire \textbf{5} (1988), no. 5, 425--464. DoI: \url{https://www.numdam.org/item/AIHPC_1988__5_5_425_0/}



\bibitem{Trudinger-1967}
N. S. Trudinger,
{\it On imbeddings into Orlicz spaces and some applications}, J. Math. Mech. \textbf{17} (1967), 473--483. DoI: \url{https://doi.org/10.1512/iumj.1968.17.17028}



\bibitem{Williams-1991}
D. Williams,
\emph{Probability with martingales},
Cambridge University Press, Cambridge, 1991. DoI: \url{https://doi.org/10.1017/CBO9780511813658}



\bibitem{Yang-2006}
Y. Yang,
{\it A sharp form of Moser-Trudinger inequality in high dimension},
J. Funct. Anal. {\bf 239} (2006), no. 1, 100--126. DoI: \url{https://doi.org/10.1016/j.jfa.2006.06.002}



\bibitem{Yang-2012}
Y. Yang,
{\it Existence of positive solutions to quasilinear elliptic equations with exponential growth in the whole Euclidean space},
J. Funct. Anal. {\bf 262} (2012), no. 4, 1679--1704. DoI: \url{https://doi.org/10.1016/j.jfa.2011.11.018}



\bibitem{Yudovich-1961}
V.I. Yudovich,
{\it Some estimates connected with integral operators and with solutions of elliptic equations},
Dokl. Akad. Nauk SSSR {\bf 138} (1961), no. 4, 805--808. DoI: \url{https://www.mathnet.ru/eng/dan25090}

\end{thebibliography}
\end{document}